\documentclass[11pt, reqno]{amsart}
\usepackage{listings}
\usepackage{amsmath,amssymb,eucal}
\usepackage{graphicx,subfigure,epsfig}
\usepackage{wrapfig, lipsum,booktabs}
\usepackage{graphicx}
\usepackage{epstopdf}
\usepackage{verbatim}
\usepackage{bm}
\usepackage{multicol}
\usepackage{multirow}
\usepackage{float}
\usepackage{color}
\usepackage{url}
\usepackage{soul}
\usepackage{tikz}
\usepackage{tikz-3dplot}
\usepackage{pgfplots}

\usepackage{algorithm} 
\usepackage[noend]{algpseudocode} 

\usepackage[top=1in, bottom=1.in, left=1in, right=1in]{geometry}

\newcommand{\bld}[1]{\hbox{\boldmath$#1$}}
\newcommand{\Th}{\mathcal{T}_h}

\newcommand{\pol}{\mathcal{P}}

\begin{document}
\title[Variational Lagrangian Hydrodynamics]{High-order variational Lagrangian schemes for compressible fluids}
\author{Guosheng Fu}
\address{Department of Applied and Computational Mathematics and
Statistics, University of Notre Dame, USA.}
\email{gfu@nd.edu}
\author{Chun Liu}
\address{Department of Applied Mathematics, Illinois Institute of Technology, USA.}
\email{cliu124@iit.edu}
 \thanks{
 G. Fu's research is partially supported by NSF grant DMS-2012031. 
 C. Liu's research is partially supported by NSF grants
DMS-1950868, DMS-2153029 and DMS-2118181.
 }

\begin{abstract}
  We present high-order variational Lagrangian finite element methods for compressible fluids using a discrete energetic variational approach.
  Our spatial discretization is mass/momentum/energy conserving and entropy stable.   
  Fully implicit time stepping is used for the temporal discretization, which allows for a much larger time step size for stability compared to  explicit methods, especially for low-Mach number flows and/or on highly distorted meshes.  
  Ample numerical results are presented to showcase the good performance of our proposed scheme.
\end{abstract}
\maketitle

\section{Introduction}
\label{sec:intro}
We are interested in the numerical simulation of multidimensional compressible fluid flows.
There are two typical choices: a Lagrangian framework in which the computational mesh moves with the fluid velocity, and an Eulerian framework in which the fluid flows through a fixed spatial mesh. 
Lagrangian methods are widely used in many fields for multi-material flow simulations
such as astrophysics and inertial confinement fusion (ICF), due to their distinguished
advantage of being able to preserve the material interface sharply.
In this article, we are concerned exclusively with Lagrangian methods. 

To address problems in Lagrangian hydrodynamics, two fundamental approaches are commonly used.
The first one is the staggered grid hydrodynamics (SGH)  approach
which employs a spatial discretization where the {\it kinematic} variables such as position and velocity are {\it continuous} and defined at the mesh vertices, and 
the {\it thermodynamic} variables including density, pressure, and internal energy are {\it discontinuous} and  defined at cell centers; see \cite{VonNeumann50,Wilkins80}. 
Artificial viscosity, as originally proposed in \cite{VonNeumann50}, is used to generate entropy production and suppress numerical oscillations across shocks. 
The second approach is the cell-centered hydrodynamics (CCH), which treats all hydrodynamic variables as cell centered and use approximate Riemann solvers to reconstruct a continuous
velocity field; see \cite{Munz94,Despres05,Maire07,Cheng07}.
This process naturally introduces a
sufficient level of dissipation at shock boundaries without the need of artificial viscosity.

While the original SGH scheme \cite{VonNeumann50} does not preserve total energy conservation, it can be recovered using 
the notion of ``corner masses" and ``corner forces" as in the compatible hydro approach \cite{Caramana98}. 
Energy-conservative high-order finite element schemes were further developed in \cite{Dobrev12} which can be viewed as a high-order generalization of the SGH approach since the key idea of using a continuous kinematic approximation space and a discontinuous thermodynamic approximation space were adopted. Advantages of the high-order Lagrangian schemes in \cite{Dobrev12} over the low-order schemes in \cite{Caramana98} include 
the ability to more accurately capture the geometry of the flow and maintain robustness with respect to mesh motion using high-order curved meshes, sharper shock resolution and better symmetry preservation for symmetric flows.

Our spatial discretization follows the SGH approach using a (high-order) continuous finite element space for the kinematic variables, and a (high-order) discontinuous
finite element space for the thermodynamic variables, which turns out to be closely related to the high-order finite element scheme \cite{Dobrev12}. 
However, we use a completely different derivation. While most of the existing Lagrangian schemes
are obtained by directly discretizing the underling PDE system,  
the starting point of our spatial discretization is an {\it energy dissipation law}. 
In particular, we present a class of variational Lagrangian schemes for compressible flows using a discrete Energetic variational approach, which is an analogue to the energetic variational approach (EnVarA) \cite{liu2009introduction, Giga2017} in a semi-discrete level. 

For a given {energy-dissipation law} and the kinematic (transport) relation, the EnVarA 
\cite{liu2009introduction, Giga2017} provides a general framework to determine the dynamics of system in a unique and well-defined way, through two distinct variational processes: Least Action Principle (LAP) and Maximum Dissipation Principle (MDP). This approach is originated from  pioneering work of Onsager \cite{onsager1931reciprocal, onsager1931reciprocal2} and Rayleigh \cite{strutt1871some}, and has been successfully applied to build up many mathematical models \cite{liu2009introduction, sun2009energetic, eisenberg2010energy, Giga2017}.
Most of existing EnVarA literature focuses on the isothermal case where temperature variation is not allowed. 
Here we adopt the approach used in \cite{LiuSulzbach22} to model a thermodynamic system with temperature variation using EnVarA.
The obtained model is then discretized in space using the discrete EnVarA, which leads to a high-order variational Lagrangian finite element scheme. Main structures of the continuous model including mass/momentum/energy conservation and entropy stability are naturally preserved in the proposed spatial discretization.
The resulting ODE system is further discretized in time using high-order implicit time integrators. Due to the special structure of the ODE system, a nonlinear system only for the velocity degrees of freedom (DOFs) needs to be solved in each time step. The allowed time step size for stability is drastically improved over explicit time stepping, especially in the low Mach number regime, at the expense of a nonlinear system solve.
Ample numerical examples are used to show the good performance of the proposed scheme.

We summarize the main features of our scheme:
\begin{itemize}
\item Space-time high-order accuracy. 
\item Mass/momentum/energy conservation and entropy stability for the spatial discretization.
\item Implicit time stepping allows for large time step size especially for low Mach number flows and/or on highly distorted meshes.
\end{itemize}

The rest of the paper is organized as follows.
In Section \ref{sec:envar}, we present the
ideal gas model using EnVarA. In Section \ref{sec:space}, a variational Lagrangian scheme is constructed using discrete EnVarA. 
The implicit temporal discretization of the resulting ODE system from Section \ref{sec:space} is then presented in Section \ref{sec:time}. 
Numerical examples are presented in Section \ref{sec:num}, followed by a summary in Section \ref{sec:conclude}.
In the Appendix, we briefly discuss our approach to the compressible isothermal case where temperature is fixed.


\section{The energetic variational approach for ideal gas}
\label{sec:envar}
The EnVarA \cite{Giga2017} is a tool to {\it derive} the force balance equation starting from a total energy $E^{\text{total}}$
and an energy dissipation rate functional $\mathcal{D}$, 
where $E^{\text{total}}$ is the sum of 
the kinetic energy $\mathcal{K}=\int_{\Omega}\frac12\rho|\bld u|^2\mathrm{dx}$
and the Helmholtz free energy $\mathcal{F}=\int_{\Omega}\psi(\rho, \theta)\mathrm{dx}$, and 
\begin{align}
\label{dissip}
2\mathcal{D}=\int_{\Omega}\eta|\nabla_s\bld u|^2+(\xi-\frac23\eta)|\nabla\cdot\bld u|^2\mathrm{dx}
\end{align}
is the rate of energy dissipation with dynamic viscosity $\eta$ and bulk viscosity $\xi$, which may depend on both $\rho$ and $\theta$. 
Here $\nabla_s$ denotes the symmetric gradient operator, $\rho$ is the density, 
$\theta$ is the absolute temperature, and $\bld u$ is the fluid velocity.

Thermodynamics of idea gas is well studied in the literature \cite{Bai99, RRR00, Bir94, Daf79, Daf16, Eri98, McQ76, Sal01}.
In classical thermodynamics, the concept of free energy proves to be useful \cite{Bai99, Sal01}.
The internal energy and pressure of an ideal gas \cite{McQ76, RRR00} have a linear relationship with temperature and the product of temperature and density, respectively.
Using this observation, Liu and Sulzbach \cite{LiuSulzbach22} proposed a definition for the free energy density of an ideal gas:
\begin{align}
\label{ener-den}
\psi(\rho, \theta) = (c_p-c_v)\theta\rho\log(\rho)-c_v\rho\theta\log(\theta),
\end{align}
 which we utilize in our current work. 
 This definition includes the specific heat 
 at constant volume $c_v$ and  specific heat at constant pressure $c_p$. 
Associated with the free energy density \eqref{ener-den}, we define the three thermodynamic variables, namely  pressure $p$, internal energy (per unit mass) $e$, and  entropy (per unit mass) $s$: 
\begin{subequations}
\label{thermo}
\begin{align}
\label{pres}
p:=&\;\psi_\rho\rho-\psi = (c_p-c_v)\rho\theta,\\
\label{init}
e:=&\;(\psi-\psi_\theta\theta)/\rho = c_v\theta,\\
\label{entr}
s:=&\;-\psi_\theta/\rho = \log(\theta^{c_v}/\rho^{c_p-c_v})+c_v,
\end{align}
\end{subequations}
which are easily verified to satisfy the famous Gibbs equation (by chain rule) \cite{Daf16}:
\begin{align}
\label{gibbs}
\theta\partial s = \partial e + p\partial (\frac{1}{\rho}),
\end{align}
where $\partial$ represents any differentiation.

In Lagrangian coordinates, we introduce the flow map: $\bld x(\bld X, t): \Omega^0\rightarrow \Omega^t$,
which  satisfies the trajectory equation
\begin{align}
\label{motion}
\frac{d}{dt}\bld x(\bld X, t) = \bld u(\bld x(\bld X, t), t),
\end{align}
with initial condition $\bld x(\bld X, 0) = \bld X$.
Since we are concerning with conventional ideal gas, we postulate the kinematics of the temperature $\theta$ being transported along the trajectory: $\frac{d}{dt}\theta=\theta_t+\bld u\cdot\nabla \theta$.

\subsection{Kinematics: mass conservation}
Within an Lagrangian control volume
$V^t:=\{\bld x(\bld X, t): \; \forall \bld X\in V^0\}$, 
 mass does not change over time:
\begin{align}
\label{mass}
\frac{d}{dt}\int_{V^t}\rho \,\mathrm{dx} = 
\int_{V^0}\frac{d}{dt}(\rho J)\,\mathrm{dX} = 0,
\end{align}
where $J=\mathrm{Det}(\nabla_X\bld x)$ is the Jacobian determinant.
Since equality \eqref{mass} is valid for any control volume $V(t)$, there must hold
\begin{align}
\label{den}
\frac{d}{dt}(\rho J)=0,\quad \text{ or }\quad
\rho(\bld x(\bld X, t), t) = \rho_0(\bld X)/J(\bld X, t),
\end{align}
where $\rho_0:\Omega^0\rightarrow \mathbb{R}^+$ is the initial density. The equality \eqref{den} represents strong mass conservation.
Writing mass conservation \eqref{den} back to the Eulerian coordinates, we have $\rho_t+\nabla\cdot(\rho\bld u) = 0$, 
which is often referred to as the continuity equation.
By abuse of notation, from it we can also get 
$\delta\rho = -\nabla\cdot(\rho\delta \bld x)$, 
where $\delta$ represents the variational derivative.
Such relation will be used and made clear in the derivations later in this paper.


\subsection{Force balance: LAP and MDP}
Here we combine the Least Action Principle (LAP)  and Maximum Dissipation Principle (MDP) \cite{onsager1931reciprocal, onsager1931reciprocal2,strutt1871some} to derive the force balance equation.
The action functional for the system is
$$
\mathcal{A}=
\int_0^T\left(\mathcal{K}-\mathcal{F}\right)\mathrm{dt}
=\int_0^T\int_{\Omega^t}\left(\frac12\rho|\bld u|^2-\psi(\rho,\theta)\right)\mathrm{dx}\mathrm{dt}.
$$


The LAP performs variation on 
the kinetic and free energies to derive the 
inertial and conservative forces:
\[
\bld f_{\text{inertial}}=\frac{\delta \int_0^T\mathcal{K}\mathrm{dt}}{\delta \bld x},\quad 
\bld f_{\text{cons}}=\frac{\delta \int_0^T\mathcal{F}\mathrm{dt}}{\delta \bld x}.
\]
Taking variation on the kinetic energy and using mass conservation \eqref{den}, we get 
\begin{align}
\label{kinetic-V}
\delta \int_0^T\mathcal{K}\mathrm{dt}=&
\delta \int_o^T\int_{\Omega^0}\frac12\rho_0(X)|\bld x_t|^2\mathrm{dXdt}
= \int_o^T\int_{\Omega^0}\rho_0(X)\bld x_t\, \delta \bld x_t\,\mathrm{dXdt}\nonumber\\
=& -\int_o^T\int_{\Omega^0}\rho_0(X)\bld x_{tt}\,\delta \bld
x\,\mathrm{dXdt}\nonumber\\
=&
-\int_o^T\int_{\Omega^0}\rho_0(X)(\bld u_t+\bld u\cdot\nabla\bld u)\,\delta \bld
x\,\mathrm{dXdt}\nonumber\\
=&
-\int_o^T\int_{\Omega^t}\rho(\bld x, t)\dot{\bld u}\,\delta \bld
x\,\mathrm{dxdt}
\end{align}
where we used the short-hand  notation $\dot{\bld u}:=\frac{d}{dt}\bld u(\bld x(\bld X, t), t)$
for the material derivative.
This implies that 
$
\bld f_{\text{inertial}}=-\rho\dot{\bld u}.
$
Taking variation on the Helmholtz free energy $\mathcal{F}$ (Hamilton's principle of virtual work) leads to 
\begin{align}
\delta \int_0^T\mathcal{F}\mathrm{dt}=&\;
\delta \int_0^T\int_{\Omega^t}\psi(\rho, \theta)\,\mathrm{dxdt}
=\;
\int_0^T\int_{\Omega^t}
\psi_{\rho}\delta\rho+\psi_\theta\delta\theta
\,\mathrm{dxdt}\nonumber\\
=&\;
\int_0^T\int_{\Omega^t}
\psi_{\rho}(-\nabla\cdot(\rho\delta \bld x))+\psi_\theta(-\delta \bld x \cdot\nabla\theta)
\,\mathrm{dxdt}\nonumber\\
=&\;
\int_0^T\int_{\Omega^t}
(\rho\nabla\psi_{\rho}-\psi_\theta\nabla\theta)\cdot \delta\bld x
\,\mathrm{dxdt}\nonumber\\
=&\;
\int_0^T\int_{\Omega^t}
\nabla p\cdot \delta\bld x
\,\mathrm{dxdt}
\label{potential-V}
\end{align}
where we used  mass conservation and the kinematic transport of the temperature defined after formula \eqref{motion} in the second row, and
the definition of the pressure \eqref{pres} in the last row:
\[
\nabla p = \nabla (\psi_\rho\rho-\psi)
=\rho\nabla\psi_\rho+\underbrace{\psi_\rho\nabla \rho
-\psi_\rho\nabla\rho}_{=0}-\psi_\theta\nabla\theta.
\]
This gives the conservative force 
$
\bld f_{\text{cons}}=\nabla p.
$

The MDP performs variation on the energy dissipation rate \eqref{dissip} to get the dissipative force:
\[
\bld{f}_{\text{diss}} = \frac{\delta \mathcal{D}}{\delta \bld u},
\]
which implies 
\[
\bld{f}_{\text{diss}} =-\nabla\cdot\left(
\eta\nabla_s\bld u\cdot+(\xi-\frac23\eta)(\nabla\cdot\bld u)\bld{I}
\right),
\]
where $\bld{I}$ is the identity matrix.

Combining these, we get the force balance equation \cite{Giga2017}
\begin{align}
\label{force-b}
\frac{\delta\mathcal{A}}{\delta \bld x} = 
\frac{\delta\mathcal{D}}{\delta \bld u},
\end{align}
which takes the following form
\begin{align}
\label{force-b2}
\rho\dot{\bld u} + \nabla p -\nabla\cdot\left(\eta\nabla_s\bld u+(\xi-\frac23\eta)(\nabla\cdot\bld u)\bld I\right)=0,
\end{align}
This is the usual momentum equation for the Navier-Stokes equations, with its natural weak formulation
\begin{align}
\label{force-weak}
\int_{\Omega^t}
\left(\rho\dot{\bld u} \cdot\bld w
-p\nabla\cdot\bld w
+\eta\nabla_s\bld u\cdot\nabla_s\bld w+(\xi-\frac23\eta)(\nabla\cdot\bld u)(\nabla\cdot\bld w)\right)\mathrm{dx}
= 0,
\end{align}
for any test function $\bld w$ with homogeneous boundary conditions.

\subsection{Internal energy and entropy equations}
The Gibbs equation \eqref{gibbs} naturally gives an update equation for the  internal energy:
\begin{align}
\label{internal}
\rho\dot{e} =  \rho\theta\dot{s}-p\rho\dot{(1/\rho)}
= \rho\theta\dot{s}-p\nabla\cdot \bld u
\end{align}
The rate of change of the entropy $s$ can be contributed by the entropy flux $\bld j$ and entropy production $\Delta$:
\begin{align}
\label{entropy}
\rho\dot{s} = \nabla\cdot\bld j + \Delta.
\end{align}
The second law of thermodynamics states that entropy production is non-negative: $\Delta \ge 0$.
If we take the entropy flux $\bld j$
being given by Durhem relation 
\begin{align}
\bld j \theta = q,
\end{align}
and the heat flux 
\begin{align}
q = \kappa \nabla\theta
\end{align}
according to Fourier's law in which $\kappa$ is the heat conductivity.
From here, the explicit expression of $\Delta$ can be derived via conservation of total energy:
there holds 
\begin{align}
\label{ener-c}
0=\frac{d}{dt}\int_{\Omega}\left(\frac12\rho|\bld u|^2+\rho e\right)\mathrm{dx}
=&\int_{\Omega}\rho\dot{\bld u}\cdot \bld u+\rho \dot{e}\mathrm{dx}.
\end{align}
Using equations \eqref{force-weak} with test function $\bld w=\bld u$, \eqref{internal}, and \eqref{entropy}, we get:
\begin{align*}
0=&\;\int_{\Omega}
-(\eta|\nabla_s\bld u|^2+(\xi-\frac23\eta)(\nabla\cdot\bld u)^2
)
+\theta\left(\nabla\cdot(\frac{\kappa \nabla \theta}{\theta})+\Delta\right)
\mathrm{dx}\\
=&\;
\int_{\Omega}
\left(-(\eta|\nabla_s\bld u|^2+(\xi-\frac23\eta)(\nabla\cdot\bld u)^2
)
-\frac{\kappa |\nabla \theta|^2}{\theta}+\theta\Delta\right)
\mathrm{dx},
\end{align*}
where we applied the chain rule for the heat flux term and 
used the homogeneous boundary condition $\nabla\theta\cdot\bld n=0$ on $\partial\Omega$.
This implies that we can take the entropy dissipation rate as 
\begin{align}
\label{entropy-d}
\Delta = 
\left(\eta|\nabla_s\bld u|^2+(\xi-\frac23\eta)(\nabla\cdot\bld u)^2
+\frac{\kappa |\nabla \theta|^2}{\theta}\right)/\theta,
\end{align}
which satisfies the second law of thermodynamics as long as 
$\theta>0$.
This implies the following entropy equation:
\begin{align}
\label{entropy2}
\rho\theta\dot{s} = \nabla\cdot(\kappa\nabla \theta) + \eta|\nabla_s\bld u|^2+(\xi-\frac23\eta)(\nabla\cdot\bld u)^2.
\end{align}
Plugging \eqref{entropy} and \eqref{entropy-d} back to the internal energy equation \eqref{internal}, we obtain:
\begin{align}
\label{internal-e}
\rho \dot{e} = -p\nabla\cdot \bld u +\nabla\cdot(\kappa\nabla \theta)+
\eta|\nabla_s\bld u|^2+(\xi-\frac23\eta)(\nabla\cdot\bld u)^2.
\end{align}
By \eqref{init}, equation \eqref{internal-e} equivalently gives the dynamics of temperature $\theta$, which has the heat equation as the leading term.

\subsection{Summary}
Combining the above results, we finally obtain the model equations:
\begin{subequations}
\label{model}
\begin{align}
\dot{\bld x} = &\;\bld u,\\
\dot{(\rho J)} = &\;0,\\
\label{model-fb}
\rho\dot{\bld u} = &\;-\nabla p +\nabla\cdot\left(\eta\nabla_s\bld u+(\xi-\frac23\eta)(\nabla\cdot\bld u)\bld I\right),\\
\rho \dot{e} = &\;-p\nabla\cdot \bld u +\nabla\cdot(\kappa\nabla \theta)+
(\eta|\nabla_s\bld u|^2+(\xi-\frac23\eta)(\nabla\cdot\bld u)^2),
\end{align}
where 
\begin{equation}
e=c_v\theta,\quad p=(c_p-c_v)\rho\theta.
\end{equation}
\end{subequations}
This is nothing but the compressible Navier-Stokes equations of an ideal gas  \cite{Daf16}.
Moreover, this system further satisfy the entropy equation \eqref{entropy2}.
In the next section we derive a variational Lagrangian scheme
for this system, where the spatial derivatives are evaluated 
by pulling back to the reference configuration (Lagrange coordinates); see, e.g., \cite{Wilkins80, Cheng07,Dobrev12}.

\section{A discrete energetic variational approach: spatial discretization}
\label{sec:space}
In this section, we construct a variational Lagrangian scheme based on a discrete EnVarA. For simplicity, we ignore heat conduction in the model, i.e. we take $\kappa=0$ in this section. 



\subsection{Notation and the finite element spaces}
Our grid-based scheme starts with a conforming triangulation $\Th^0=\{T_\ell^0\}_{\ell=1}^{N_T}$ of the initial configuration $\Omega^0$ with $N_T$ elements, where we assume the element $T_\ell^0:=\Phi_{T_\ell}^0(\widehat{T})$
is obtained from a polynomial mapping $\Phi_{T_\ell}^0$
from the reference element $\widehat{T}$, which, for simplicity,  is a simplex or a hypercube. 

We denote  $\pol^k(\widehat{T})$ as the polynomial space of degree no greater than $k$ if $\widehat{T}$
is a reference simplex, or the 
tensor-product polynomial space of degree no greater than $k$ in each direction if $\widehat{T}$
is a reference hypercube, for $k\ge 1$.
The mapped polynomial space on a spatial physical element
$T^0\in\Th^0$ 
is denoted as 
\[
\pol^k(T^0):=
\{\widehat{v}\circ(\Phi_{T}^0)^{-1}:\; \forall 
\widehat{v}\in\pol^k(\widehat{T}))
\}.
\]

We denote $\{\widehat{\bm \xi}_i\}_{i=1}^{N_k}$
as a set of quadrature points with positive weights 
$\{\widehat{\omega}_i\}_{i=1}^{N_k}$ that is accurate for polynomials of degree up to $2k+1$ on the reference element $\widehat T$, i.e.,
\begin{align}
\label{quad-2k}
\int_{\widehat T}\widehat f\, \mathrm{dx}
=\sum_{i=1}^{N_k}\widehat{\omega}_i\widehat f(\widehat{\bm \xi}_i), \quad\forall \widehat f\in\pol^{2k+1}(\widehat T).
\end{align}
Note that when $\widehat T$ is a reference square, we simply use the Gauss-Legendre quadrature rule with $N_k = (k+1)^2$, which is optimal.
On the other hand, when $\widehat T$ is a reference simplex, the optimal choice of quadrature rule is more complicated; see, e.g., 
\cite{Zhang09,Witherden15} and references cited therein.
Table \ref{tab:quad} list the number $N_k$ for $0\le k\le 6$ of the symmetric quadrature rules provided in \cite{Zhang09}.
\begin{table}[ht!]
\begin{tabular}{cccccccc}
 & $k=0$ & $k=1$ & $k=2$ & $k=3$ & $k=4$& $k=5$ & $k=6$\\\hline
$N_k$ on Triangle 
& 1 & 6 & 7 & 15 & 19&28&37\\
$N_k$ on Tetrahedron
& 1 & 8 & 14 & 36 & 61 & 109 & 171\\\hline
\end{tabular}
\caption{Number of quadrature points $N_k$ for the quadrature rule 
on a simplex that is accurate up to degree $2k+1$ for $0\le k\le 6$.
}
\label{tab:quad}
\end{table}
The integration points and weights on a physical element $T_\ell^0$ are simply obtained via mapping: $\{\bm \xi_i^{\ell}:=\Phi_{T_\ell^0}(\widehat{\bm \xi}_i)\}_{i=1}^{N_k}$, and 
$\{\omega_i^{\ell}:=|\nabla \Phi_{T_\ell^0}(\widehat{\bm \xi}_i)|\widehat \omega_i\}_{i=1}^{N_k}$.
To simplify the notation, we denote the set of physical integration points and weights
\begin{subequations}
\label{quad-rule}
\begin{align}
\Xi_h^k:=&\;\{\bm\xi_i^{\ell}:\;\;1\le i\le N_k, \, 1\le \ell\le N_T\},\\
\Omega_h^k:=&\;\{\omega_i^{\ell}:\;\;1\le i\le N_k, \, 1\le \ell\le N_T\},
\end{align}
\end{subequations}
and denote $(\cdot, \cdot)_h$ as the discrete inner-product on the mesh $\Th^0$ using the quadrature points $\Xi_h^k$ and weights $\Omega_h^k$:
\[
(\alpha, \beta)_h:=\sum_{\ell=1}^{N_T}\sum_{i=1}^{N_k}\alpha(\bm\xi_i^{\ell})
\beta(\bm\xi_i^{\ell})\omega_i^{\ell}.
\]

We are now ready to present our 
continuous and discontinuous finite element spaces: 
\begin{align}
\label{fes-V}
\bld V_h^k:=&\;\{\bld v\in [H^1(\Omega^0)]^d: \;\; \bld v|_{T_\ell^0}
\in[\pol^k(T_\ell^0)]^d,\;\;
\forall T_\ell^0\in\Th^0\},\\
\label{fes-W}
W_h^k:=&\;\{w\in L^2(\Omega^0): \;\; w|_{T_\ell^0}
\in W^{k}(T_\ell^0),\;\;
\forall T_\ell^0\in\Th^0\},
\end{align}
where the local space 
\[
W^k(T_\ell^0):=\pol^{k}(T_\ell^0)\oplus \delta W_k(T_\ell^0),
\]
is associated with the integration rule in \eqref{quad-2k}
such that $\dim W^k(T_\ell^0) = N_k$, and the nodal conditions \begin{align}
\label{collocation}
\varphi_i^{\ell}({\bm \xi}_j^{\ell})=
\delta_{ij},\quad \forall 1\le j\le N_k,
\end{align}
in which  $\delta_{ij}$ is the Kronecker delta function determines a unique solution $\varphi_i^\ell\in W^k(T_\ell^0)$. 
This implies that $\{\varphi_i^\ell\}_{i=1}^{N_k}$ is a set of nodal bases for the space ${W}^k(T_\ell^0)$, i.e., 
\begin{align}
\label{nodal-basis}
{W}^k(T_\ell^0)=\mathrm{span}_{1\le i\le N_k}\{\varphi_i^\ell\}.
\end{align}
When $T_\ell^0$ is a mapped hypercube, 
we have $N_k=(k+1)^2$, hence $\delta W_k(T_\ell^0)=\emptyset$ and
$W^k(T_\ell)$ is simply the (mapped) tensor product polynomial space $\pol^{k}(T_\ell^0)$.
On the other hand, when $T_\ell^0$ is a mapped simplex, 
we have $\dim \delta W_k(T_\ell^0)= N_k-\dim \pol^k(T_\ell^0)>0$ for $k\ge 1$.
However, we emphasize that the explicit expression of 
$\delta W_k(T_\ell^0)$ or
the basis function $\phi_i^\ell$ does not matter in our discretization, as only their nodal degrees of freedom (DOFs) on the quadrature nodes will enter into the numerical integration. 
Any function $\alpha(\bld X, t)$ in $W_h^k$ (for fixed $t$) can be expressed as 
\[
\alpha(\bld X, t) = \sum_{\ell=1}^{N_T}\sum_{i=1}^{N_k}{\sf a}_{i}^\ell(t)
\phi_i^\ell(\bld X),
\]
where $\{{\sf a}_{i}^\ell(t)\}$
are the unknown coefficients.
We refer to $W_h^k$ as the (discontinuous) {\it integration rule space}, which {\it only } contains the $N_T\times N_k$ quadrature points and weights \eqref{quad-rule}, and is easy to implement in practice. 

We further denote a set of basis functions for $\bld V_h^k$ as 
$\{\bld \varphi_i\}_{i=1}^{N_V^k}$, where $N_V^k$ is the dimension of $\bld V_h^k$.
We use the continuous space $\bld V_h^k$ to approximate the flow map $\bld x_h$ and velocity $\bld u_h$, and 
the discontinuous space $W_h^k$ to approximate the 
density $\rho_h$, pressure $p_h$, internal energy $e_h$,
temperature $\theta_h$, and entropy $s_h$. 
More specifically, we have 
\begin{subequations}
\label{approx}
\begin{alignat}{2}
\bld x_h(\bld X, t) =&\; \sum_{i=1}^{N_V^k}{\sf x}_i(t)\bld\varphi_i(\bld  X), \quad &&\bld u_h(\bld X, t) =\; \sum_{i=1}^{N_V^k}{\sf u}_i(t)\bld\varphi_i(\bld X),\\
\rho_h(\bld X, t) =&\; \sum_{\ell=1}^{N_T}\sum_{i=1}^{N_k}{\rho}_i^\ell(t)\phi_i^\ell(\bld X),\quad
&&
p_h(\bld X, t) =\; \sum_{\ell=1}^{N_T}\sum_{i=1}^{N_k}{\sf p}_i^\ell(t)\phi_i^\ell(\bld X),\\
e_h(\bld X, t) =&\; \sum_{\ell=1}^{N_T}\sum_{i=1}^{N_k}{e}_i^\ell(t)\phi_i^\ell(\bld X),\quad
&&
\theta_h(\bld X, t) =\; \sum_{\ell=1}^{N_T}\sum_{i=1}^{N_k}{\theta}_i^\ell(t)\phi_i^\ell(\bld X),\\
s_h(\bld X, t) =&\; \sum_{\ell=1}^{N_T}\sum_{i=1}^{N_k}{\sf s}_i^\ell(t)\phi_i^\ell(\bld X),
\end{alignat}
\end{subequations}
where ${\sf X}_h=[{\sf x}_1, \cdots, {\sf x}_{N_V^k}]^T$, 
${\sf U}_h=[{\sf u}_1, \cdots, {\sf u}_{N_V^k}]^T$, 
${\sf R}_h=[{\sf \rho}_1^1, \cdots, {\sf \rho}_{N_k}^{N_T}]^T$, 
${\sf P}_h=[{\sf p}_1^1, \cdots, {\sf p}_{N_k}^{N_T}]^T$, 
${\sf E}_h=[{\sf e}_1^1, \cdots, {\sf e}_{N_k}^{N_T}]^T$, 
${\sf \Theta}_h=[{\sf \theta}_1^1, \cdots, {\sf \theta}_{N_k}^{N_T}]^T$,
and
${\sf S}_h=[{\sf s}_1^1, \cdots, {\sf s}_{N_k}^{N_T}]^T$
are the 
time dependent coefficient vectors for $\bld x_h$, $\bld u_h$, 
$\rho_h$, $p_h$, $e_h$, $\theta_h$, and $s_h$, respectively.
Note that by the thermodynamic relations \eqref{thermo}, there are only two independent thermodynamic variables. 
Here we take density $\rho$ and temperature $\theta$
as the independent variables. The other variables will be updated through the discrete formulas \eqref{thermo-h} below.

\subsection{Trajectory equation and mass conservation}
With the notation given in \eqref{approx}, the trajectory equation \eqref{motion} simply implies that 
\begin{align}
\label{traj}
{\sf X}_h'(t) = {\sf U}_h(t).
\end{align}
We require mass conservation to be satisfied pointwise at the quadrature nodes level, specifically, \eqref{den} implies that 
\begin{align}
\label{mass-c}
{\sf\rho}_i^\ell(t) = 
\rho_0(\bld \xi_i^\ell)/J_h(\bld \xi_i^\ell, t), 
\quad \forall 1\le \ell\le N_T, \;1\le i\le N_k,
\end{align}
where the discrete Jacobian on the quadrature point $\bld \xi_i^\ell$ is 
$J_h(\bld \xi_i^\ell, t):=
\mathrm{Det}(\nabla_X\bld x_h(\bld \xi_i^\ell, t))
$.

\subsection{Thermodynamic relations}
We require the relations in \eqref{thermo} for the thermodynamic variables be satisfied on the quadrature points level, which implies 
\begin{subequations}
\label{thermo-h}
\begin{align}
\label{pres-h}
{\sf p}_i^\ell = &\;(c_p-c_v)\rho_i^\ell\theta_i^\ell,\\
\label{init-h}
{\sf e}_i^\ell = &\; c_v\theta_i^\ell,\\
\label{entr-h}
{\sf s}_i^\ell = &\;
c_v\log(\theta_i^\ell)
-(c_p-c_v)\log(\rho_i^\ell)+c_v,
\end{align}
\end{subequations}
for all $1\le i\le N_k$ and $1\le \ell\le N_T$.
It is easy to see that the (pointwise) Gibbs equation \eqref{gibbs} is satisfied:
\begin{align*}
\rho_i^\ell\theta_i^\ell({\sf s}_i^\ell)'= 
\rho_i^\ell ({\sf e}_i^\ell)'
-\frac{{\sf p}_i^\ell}{\rho_i^\ell}(\rho_i^\ell)'.
\end{align*}
By the density  definition \eqref{mass-c} and Jacobi's formula, we have 
\[
(\rho_i^\ell)'=-\rho_i^\ell\nabla\cdot\bld u_h(\bld \xi_i^\ell),
\]
which implies that
\begin{align}
\label{gibbs-h}
\rho_i^\ell\theta_i^\ell({\sf s}_i^\ell)'= 
\rho_i^\ell ({\sf e}_i^\ell)'
+{{\sf p}_i^\ell}\nabla\cdot\bld u_h(\bld \xi_i^\ell).
\end{align}
We note that the above pointwise relations also hold for the classical low-order SGH schemes \cite{VonNeumann50,Wilkins80,Caramana98} where the thermodynamic variables are approximated via piecewise constants, which, however, does not hold in general  
for the high order scheme \cite{Dobrev12} due to the use of a different high-order thermodynamic finite element space.

\subsection{Discrete EnVarA and velocity equation}
Instead of discretizing the force balance equation \eqref{model-fb}, 
here we discretize the energy law and use 
the EnVarA to {\it derive} the discrete force balance equation directly.
We denote the discrete action functional
\begin{align}
\label{disc-A}
\mathcal{A}_h:=\int_0^T(\mathcal{K}_h-\mathcal{F}_h)\mathrm{dt},
\end{align}
where the discrete kinetic energy $\mathcal{K}_h$ and 
the discrete Helmholtz free energy $\mathcal{F}_h$ are given as 
\[
\mathcal{K}_h=\frac12(\rho_hJ_h\bld u_h, \bld u_h)_h=
\sum_{\ell=1}^{N_T}\sum_{i=1}^{N_k}\frac12\rho_0(\bld \xi_i^\ell)|\bld u_h(\bld \xi_i^\ell, t)|^2\omega_i^\ell,
\]
and 
\[
\mathcal{F}_h=
(\psi(\rho_h, \theta_h)J_h, 1)_h=
\sum_{\ell=1}^{N_T}\sum_{i=1}^{N_k}\psi\left(\rho_h(\bld \xi_i^\ell, t), \theta_h(\bld \xi_i^\ell,t)\right)J_h(\bld \xi_i^\ell,t)\omega_i^\ell.
\]
Moreover, the discrete dissipation rate is given as 
\begin{align}
\label{disc-D}
\mathcal{D}_h:=\frac12(\eta J_h\nabla_s\bld u_h, \nabla_s\bld u_h)_h
+\frac12\left((\xi-\frac23\eta)J_h\nabla\cdot\bld u_h, \nabla\cdot\bld u_h\right)_h.
\end{align}
The discrete force balance equation \eqref{force-b} is then 
\begin{align}
\label{disc-fb}
\frac{\delta\mathcal{A}_h}{\delta {\sf x}_j}
= 
\frac{\delta\mathcal{D}_h}{\delta {\sf u}_j},\quad \forall 1\le j\le N_V^k.
\end{align}
Elementary  calculation, using \eqref{traj}, \eqref{mass-c}
and Jacobi's formula,
yields that 
\begin{align}
\label{var-A}
\frac{\delta\mathcal{A}_h}{\delta {\sf x}_j}
=-(\rho_0\dot{\bld u}_h, \bld \varphi_j)_h
+(p_hJ_h, \nabla\cdot\bld \varphi_j)_h,
\end{align}
where $\dot{\bld u}_h=\sum_{i=1}^{N_V^k}{\sf u}_i(t)'\bld\varphi_i$, and
the pressure $p_h\in W_h^k$ satisfies 
\begin{align*}
{\sf p}_i^\ell = \psi_{\rho}(\rho_i^\ell, \theta_i^\ell)\rho_i^\ell-\psi(\rho_i^\ell, \theta_i^\ell)
= (c_p-c_v)\rho_i^\ell\theta_i^\ell,
\end{align*}
according to \eqref{pres-h}.
We also have
\begin{align}
\label{var-D}
\frac{\delta\mathcal{D}_h}{\delta {\sf u}_j}
=(\eta J_h\nabla_s{\bld u}_h, \nabla_s\bld \varphi_j)_h
+((\xi-\frac23\eta) J_h\nabla\cdot{\bld u}_h, \nabla\cdot\bld \varphi_j)_h, \quad \forall 1\le j\le N_V^k.
\end{align}
Plugging \eqref{var-A} and \eqref{var-D} back to \eqref{disc-fb}, and using the definition of the pressure, 
we get  the semi-discrete force balance equation:
\begin{align}
\label{disc-F}
(\rho_0\dot{\bld u}_h, \bld \varphi_j)_h
-\left((c_p-c_v)\rho_0\theta_h, \nabla\cdot\bld \varphi_j\right)_h
+(\sigma_h, \nabla\bld \varphi_j)_h
= 0 , \quad \forall 1\le j\le N_V^k,
\end{align}
where $\sigma_h$ is the viscous stress defined as 
\begin{align}
\label{stress}
\sigma_h:=\eta J_h\nabla_s{\bld u}_h+ (\xi-\frac23\eta) J_h\nabla\cdot{\bld u}_h\bld I.
\end{align}
Here the evaluation of spatial derivative terms shall be pulled back to the initial configuration $\Omega^0$. In particular, 
\[
\nabla\bld u_h = (F_h)^{-1}\nabla_X\bld u_h,
\]
where $F_h = \nabla_X\bld x_h$ is the deformation tensor.
Equation \eqref{disc-F} provides an ODE system for the velocity coefficient vector ${\sf U}_h(t)$.
Taking test function $\bld \varphi_j$ as a constant vector, 
we immediately obtain global momentum conservation:
\[
\frac{d}{dt}(\rho_0 \bld u_h, 1)_h=(\rho_0\dot{\bld u}_h, 1)_h=0.
\]

\subsection{Energy conservation  and temperature equation}
We use energy conservation to get an update equation 
for the temperature coefficient vector ${\sf \Theta}_h(t)$.
In the absence of heat conduction ($\kappa=0$), the spatial discretization of the internal energy equation \eqref{internal-e} 
leads to
\begin{align}
\label{int-eh}
(\rho_h\dot{e}_h, J_h\phi_j^\ell)_h = 
-(p_h J_h\nabla\cdot\bld u_h,\phi_j^\ell)_h
+(\sigma_h:\nabla\bld u_h, \phi_j^\ell)_h, \quad\forall 
1\le \ell\le N_T,\, 1\le j\le N_k.
\end{align}
Equivalently, the equation \eqref{int-eh} has
the following pointwise form for the coefficient vector ${\sf E}_h$:
\begin{align}
\label{int-eC}
\rho_0(\bld \xi_j^\ell)({\sf e}_j^\ell)' = 
-{\sf p}_j^\ell J_h(\bld \xi_j^\ell)
(\nabla\cdot\bld u_h)(\bld \xi_j^\ell)
+\sigma_h(\bld \xi_j^\ell):\nabla\bld u_h(\bld \xi_j^\ell), 
\end{align}
for all $1\le \ell\le N_T,\,1\le j\le N_k$.
Plugging in the relations \eqref{pres-h} and \eqref{init-h}
back to \eqref{int-eC}, 
we get the following ODE system for the coefficent vector $\Theta_h(t)$:
\begin{align}
\label{int-eC2}
c_v\rho_0(\bld \xi_j^\ell)({\sf \theta}_j^\ell)' = 
-(c_p-c_v)\rho_0(\bld \xi_j^\ell)
(\nabla\cdot\bld u_h)(\bld \xi_j^\ell) {\sf \theta}_j^\ell
+\sigma_h(\bld \xi_j^\ell):\nabla\bld u_h(\bld \xi_j^\ell), 
\end{align}
One key observation is that \eqref{int-eC2} is a {\it linear}
ODE system for ${\sf \Theta}_h$.
Moreover, taking $\bld\varphi_j=\bld u_h$
in \eqref{disc-F} and $\phi_j^\ell=1$ in \eqref{int-eh}
and adding, we obtain the total energy conservation:
\[
\frac{d}{dt}\left(\frac12\rho_0|\bld u_h|^2+\rho_0e_h, 1\right)_h 
=
(\rho_0\dot{\bld u}_h, \bld u_h)_h+(\rho_0\dot{e}_h, 1)_h 
=0.
\]

\subsection{The entropy equation}
Combining the Gibbs equation \eqref{gibbs-h}
with the internal energy equation \eqref{int-eC} and simplifying, 
we obtain the ODE system satisfied by the entropy:
\begin{align}
\label{ent-X2}
J_h(\bld \xi_j^\ell)\rho_i^\ell\theta_i^\ell({\sf s}_i^\ell)'= 
\sigma_h(\bld \xi_j^\ell):\nabla\bld u_h(\bld \xi_j^\ell), 
\end{align}
for all $1\le \ell\le N_T$, $1\le j\le N_k$.
This implies that 
\begin{align}
\label{entro-sb}
\frac{d}{dt}(\rho_0 s_h, 1)_h=(\rho_0\dot{s}_h, 1)_h = (\frac{\sigma_h:\nabla\bld u_h}{\theta_h}, 1)_h,
\end{align}
where positivity of the right hand side, i.e., semi-discrete entropy stability, is guaranteed 
as long as the temperature $\Theta_h >0$.
We remark that the entropy stability \eqref{entro-sb}
is a direct consequence of our special choice of (nodal) thermodynamic finite element space \eqref{fes-W} and \eqref{nodal-basis}.

\subsection{Artificial viscosity}
To make the scheme robust even in the case of zero {\it physical viscosities} with $\eta=\xi=0$, we add artificial viscosity \cite{VonNeumann50} to the system so that shocks can be dissipated.
Specifically, we add to the stress term \eqref{stress} an artificial stress tensor $\sigma_h^{av}$ of the following form:
\begin{align}
\label{av-1}
\sigma_h\leftarrow \sigma_h+\sigma_h^{av},\quad\text{ where }
\sigma_h^{av}:= \mu_{av}J_h\nabla_s\bld u_h,
\end{align}
in which, following \cite{Dobrev12}, the artificial viscosity coefficient $\mu_{av}$ is: 
\begin{align}
\label{av}
\mu_{av}=\rho_h\left(
q_2\ell_{s_1}^2|\Delta_{s_1}\bld u_h|+q_1\psi_0\psi_1\ell_{s_1}c_s
\right)
\end{align}
where $q_1$ and $q_2$ are linear and quadratic scaling coefficients, $c_s=\sqrt{\gamma p_h/\rho_h}$ is the speed of sound with $\gamma=c_p/c_v$ being the adiabatic constant, $\Delta_{s_1}\bld u_h:={s_1\cdot\nabla\bld u_h\cdot s_1}$
is the {directional measure of compression}
and $\ell_{s_1}=\ell_0 |J_h s_1|$
is the {directional length scale} along the 
 direction $s_1$, 
 and the two linear switches are
$\phi_0 = \frac{|\nabla\cdot \bld u_h|}{\|\nabla\bld u_h\|}$,
and 
$
\phi_1 = \begin{cases}
1, &\text{if } \Delta_s \bld u_h<0,\\
0, &\text{if } \Delta_s \bld u_h\ge0.\\
\end{cases}
$
Here the direction $s_1$ is the unit eigenvector 
of the symmetric tensor $\nabla_s\bld u_h$ with the smallest eigenvalue $\lambda_1$, i.e.,  
\[
(\nabla_s\bld u_h) s_1=\lambda_1 s_1, \quad |s_1| = 1, \text{ and $\lambda_1$ is the smallest eigenvalue.}
\]
With this notation, we have $\Delta_{s_1}\bld u_h=\lambda_1$.
Moreover $\ell_0 = h_0/k$ is the mesh size of the initial domain divided by the polynomial degree $k$. 
We refer interested reader to \cite{Dobrev12} and references cited therein for more discussion on the choice of the
artificial viscosity coefficient.

\subsection{Summary}
The final form of the semi-discrete scheme is 
summarized in Algorithm \ref{alg:1} below.
This spatial discretization is high-order, mass/momentum/energy conserving, and entropy stable.
\begin{algorithm}[H]
\caption{Spatial discretization for model \eqref{model}.}
\label{alg:1}
\begin{algorithmic}
\State 
$\bullet$ Find $\bld x_h, \bld u_h\in \bld V_h^k$, 
and $\theta_h\in W_h^k$ such that 
the ODE system 
\begin{align*}
{\sf X}_h'(t) =&\; {\sf U}_h(t), \\
(\rho_0\dot{\bld u}_h, \bld \varphi_h)_h
-\left((c_p-c_v)\rho_0\theta_h, \nabla\cdot\bld \varphi_h\right)_h
+(\sigma_h, \nabla\bld \varphi_h)_h
= &\;0 , \quad\forall \bld\varphi_h\in \bld V_h^k,\\
({\sf \theta}_j^\ell)' +(\gamma-1)
(\nabla\cdot\bld u_h)(\bld \xi_j^\ell) {\sf \theta}_j^\ell
-\frac{\sigma_h(\bld \xi_j^\ell)}{c_v\rho_0(\bld \xi_j^\ell)}:\nabla\bld u_h(\bld \xi_j^\ell)
=&\;0,\quad\forall 1\le \ell\le N_T, 1\le j\le N_k,
\end{align*}
holds for the coefficient vectors ${\sf X}_h$, ${\sf U}_h$, and ${\sf \Theta}_h$, where the numerical stress 
\[
\sigma_h=(\eta+\mu_{av}) J_h\nabla_s{\bld u}_h+ (\xi-\frac23\eta)J_h\nabla\cdot{\bld u}_h\bld I,
\]
in which the artificial viscosity $\mu_{av}$ is given in \eqref{av}.
Here the notation \eqref{approx} for the finite element approximations is used. 
\State 
$\bullet$ The density approximation $\rho_h\in W_h^k$ satisfies mass conservation
\eqref{mass-c}, and 
the pressure, internal energy, and entropy approximations
$p_h, e_h, s_h\in W_h^k$ satisfy the thermodynamic relations \eqref{thermo-h}.
\end{algorithmic}
\end{algorithm}




\section{Temporal discretization}
\label{sec:time}
In this section, we focus on the discretization of the ODE system in Algorithm \ref{alg:1}.
We use fully implicit time discretizations so that the fully discrete scheme is robust for all mach numbers. 
We refer to \cite{Dobrev12,Adrian21} for alternative conservative explicit schemes.

\subsection{First order energy dissipative scheme}
\label{sec:BE}
Using implicit Euler for the time derivative terms, we arrive at the following first order scheme:
Given data ${\bld x}_h^{n-1}, {\bld u}_h^{n-1}\in \bld V_h^k$ and ${\theta}_h^{n-1}\in W_h^k$ at time $t^{n-1}$, and time step size $\delta t$, find solution
${\bld x}_h^{n}, {\bld u}_h^{n}\in \bld V_h^k$ and ${\theta}_h^{n}\in W_h^k$ such that
\begin{subequations}
\label{BE}
\begin{align}
\label{x-Th1}
\frac{{\bld x}_h^{n}-{\bld x}_h^{n-1}}{\delta t} = &\;{\bld  u}_h^n, \\
\label{u-Th1}
(\rho_0\frac{{\bld u}_h^n-{\bld u}_h^{n-1}}{\delta t}, \bld \varphi_h)_h
-\left((c_p-c_v)\rho_0\theta_h^n, \nabla\cdot\bld \varphi_h\right)_h
+(\sigma_h^{n}, \nabla\bld \varphi_h)_h
= &\;0 , \quad \forall \varphi_h\in\bld V_h^k,\\
\label{theta-Th1}
\frac{{\sf \theta}_j^{\ell, n}-
{\sf \theta}_j^{\ell, n-1}}{\delta t} + (\gamma-1)
\nabla\cdot\bld u_h^n(\bld \xi_j^{\ell}) {\sf \theta}_j^{\ell,n}
-\frac{\sigma_h^{n}(\bld \xi_j^\ell):\nabla\bld u_h^n(\bld \xi_j^\ell)}{c_v\rho_0(\bld \xi_j^\ell)} =&\; 0, \quad \forall j, \ell,
\end{align}
\end{subequations}
where the stress
\[
\sigma_h^n=(\eta+\mu_{av}^{n-1}) J_h^n\nabla_s{\bld u}_h^n+ (\xi-\frac23\eta)J_h^n\nabla\cdot{\bld u}_h^n\bld I,
\]
in which the artificial viscosity coefficient $\mu_{av}^{n}$ is evaluated at time $t^{n}$, and 
the Jacobian determinant $J_h^n=|\nabla_X\bld x_h^n|$.
The above system can be solved by first expressing $\bld x_h^n$ and $\theta_h^n$ in terms of $\bld u_h^n$
using \eqref{x-Th1} and \eqref{theta-Th1}:
\begin{subequations}
\label{x-theta}
\begin{align}
\bld x_h^n = &\;\bld x_h^{n-1} + \delta t \bld u_h^n,\\
\label{sss}
{\sf \theta}_j^{\ell, n}=&\;
\frac{{\sf \theta}_j^{\ell, n-1}+\delta t\frac{\sigma_h^n(\bld \xi_i^\ell) : \nabla \bld u_h^{n}(\bld \xi_i^\ell)}{c_v\rho_0(\bld \xi_i^\ell)}}{1+\delta t(\gamma-1)\nabla\cdot\bld u_h^n(\bld \xi_i^\ell)},
\end{align}
\end{subequations}
and then solve the nonlinear system for $\bld u_h^n$
in \eqref{u-Th1} using \eqref{x-theta}. 
We use Newton's method to solve this nonlinear system for the velocity DOFs.

For the scheme \eqref{BE}, positivity of density on the quadrature points is guaranteed as long as the Jacobian $J_h^n>0$ on these quadrature points. And positivity of the temperature (hence positivity of pressure and internal energy) is satisfied 
as long as the denominator of right hand side of \eqref{sss}
stays positive, i.e., 
\[
1+\delta t(\gamma-1)\nabla\cdot\bld u_h^n(\bld \xi_i^\ell)>0,
\quad \forall i,\ell.
\]
Moreover, strong mass conservation is satisfied due to \eqref{mass-c}, and global momentum conservation is satisfied by taking $\bld \varphi_j$ to be a global constant in \eqref{u-Th1}. 
Finally, taking $\bld\varphi_j=\bld u_h^n$ in \eqref{u-Th1}
and combining with \eqref{theta-Th1}, we get 
\begin{align*}
(\rho_0(\bld u_h^{n}-\bld u_h^{n-1}), \bld u_h^n)_h
+(\rho_0(e_h^{n}-e_h^{n-1}), 1)_h = 0, 
\end{align*}
which implies that 
\begin{align*}
\left(\rho_0(\frac12|\bld u_h^n|^2+e_h^n),1\right)_h
-
\left(\rho_0(\frac12|\bld u_h^{n-1}|^2+e_h^{n-1}),1\right)_h
=
-\left(\rho_0\frac12|\bld u_h^n-\bld u_h^{n-1}|^2,1\right)_h
\le 0
\end{align*}
Hence, the total energy is {\it dissipated} over time for the scheme 
\eqref{BE}.

\subsection{Second order energy-conservative scheme}
\label{sec:CN}
Energy conservation can be recovered from the first order scheme \eqref{BE} by applying a time filter, which is the same as the midpoint rule; see \cite{Burkardt20}.
This time stepping algorithm has two steps and is recorded in Algorithm \ref{alg:2} for reference.
\begin{algorithm}[H]
\caption{The midpoint rule with Backward Euler -- Forward Euler implementation.}
\label{alg:2}
\begin{algorithmic}
\State 
$\bullet$ Apply the backward Euler scheme \eqref{BE} with half time step $\delta t/2$ to get approximations at midpoint $t^{n-\frac12}=t^{n-1}+\delta t/2$, and denote the solutions as 
$\bld x_h^{n-\frac12}$,
$\bld u_h^{n-\frac12}$, and 
$\theta_h^{n-\frac12}$.
\State 
$\bullet$ 
Apply a time filter (forward Euler) step to 
approximation the solutions at time $t^{n}=t^{n-1}+\delta t$:
\begin{align}
\label{FE}
\bld x_h^{n} = 2\bld x_h^{n-\frac12}-\bld x_h^{n-1},\;\;
\bld u_h^{n} = 2\bld u_h^{n-\frac12}-\bld u_h^{n-1},\;\;
\theta_h^{n} = 2\theta_h^{n-\frac12}-\theta_h^{n-1}.
\end{align}
\State
$\bullet$ The density, pressure, internal energy, and entropy approximations
$\rho_h, p_h, e_h, s_h\in W_h^k$ are then recovered through \eqref{mass-c} and \eqref{thermo-h}.
\end{algorithmic}
\end{algorithm}
We notice that the BE step implies
\begin{align*}
(\rho_0(\bld u_h^{n-\frac12}-\bld u_h^{n-1}), \bld u_h^{n-\frac12})_h
+(\rho_0(e_h^{n-\frac12}-e_h^{n-1}), 1)_h = 0,
\end{align*}
By the extrapolation relations \eqref{FE}, we have 
\[
(\bld u_h^{n-\frac12}-\bld u_h^{n-1})\cdot \bld u_h^{n-\frac12}
=\frac14(|\bld u_h^n|^2-|\bld u_h^{n-1}|^2), 
\text{ and }
e_h^{n-\frac12}-e_h^{n-1}
=\frac12(e_h^{n-1}-e_h^{n-1}).
\]
Combining these equations, we get total energy conservation for this two-step method:
\[
\left(\rho_0(\frac12|\bld u_h^n|^2+e_h^n),1\right)_h
=
\left(\rho_0(\frac12|\bld u_h^{n-1}|^2+e_h^{n-1}),1\right)_h.
\]

\subsection{High-order BDF schemes}
\label{sec:BDF}
It is natural to extend the implicit Euler scheme 
\eqref{BE} to higher order by replacing the backward Euler time difference terms in \eqref{BE} using higher-order backward difference formulas (BDFs). 
For completeness, we record these high-order schemes with uniform time stepping in Algoritm~\ref{alg:3} below.
\begin{algorithm}[H]
\caption{The BDF$[m]$ scheme.}
\label{alg:3}
\begin{algorithmic}
\State 
$\bullet$
Given data ${\bld x}_h^{n-j}, {\bld u}_h^{n-j}\in \bld V_h^k$ and ${\theta}_h^{n-j}\in W_h^k$ at time $t^{n-j}=(n-j)\delta t$ for $j=1\cdots, m$, find 
solution
${\bld x}_h^{n}, {\bld u}_h^{n}\in \bld V_h^k$ and ${\theta}_h^{n}\in W_h^k$ 
at time $t^n=n\delta t$ such that
\begin{subequations}
\label{BDF}
\begin{align}
\label{x-Th3}
{D}_{\delta t}^m(\bld x_h^n) = &\;{\bld  u}_h^n, \\
\label{u-Th3}
(\rho_0{D_{\delta t}^m({\bld u}_h^n)}, \bld \varphi_h)_h
-\left((c_p-c_v)\rho_0\theta_h^n, \nabla\cdot\bld \varphi_h\right)_h
+(\sigma_h^{n}, \nabla\bld \varphi_h)_h
= &\;0 , \quad \forall \varphi_h\in\bld V_h^k,\\
\label{theta-Th3}
D_{\delta t}^m({\sf \theta}_j^{\ell, n}) + (\gamma-1)
\nabla\cdot\bld u_h^n(\bld \xi_j^{\ell}) {\sf \theta}_j^{\ell,n}
-\frac{\sigma_h^{n}(\bld \xi_j^\ell):\nabla\bld u_h^n(\bld \xi_j^\ell)}{c_v\rho_0(\bld \xi_j^\ell)} =&\; 0, \quad \forall j, \ell,
\end{align}
where $D_{\delta t}^m(\alpha)$ is the 
approximation to the time derivative term $\alpha(t)'$ using BDF$[m]$, e.g., 
\[
D_{\delta t}^2(\alpha^n)=\frac{3\alpha^n-4\alpha^{n-1}+\alpha^{n-2}}{2\,\delta t},\quad
D_{\delta t}^3(\alpha^n)=
\frac{11\alpha^n-18\alpha^{n-1}+9\alpha^{n-2}-2\alpha^{n-3}}{6\,\delta t}.
\]
\end{subequations}
\State
$\bullet$ The density, pressure, internal energy, and entropy approximations
$\rho_h, p_h, e_h, s_h\in W_h^k$ are then recovered through \eqref{mass-c} and \eqref{thermo-h}.
\end{algorithmic}
\end{algorithm}

\section{Numerical results}
\label{sec:num}
We present numerical results in this section using the open-source finite-element software
{\sf NGSolve} \cite{Schoberl16}, \url{https://ngsolve.org/}.
For all the simulation results, we consider invisid models with 
{\it zero} physical viscosities. 

We use a (variable time step size) BDF2 time stepping  with high-order spatial discretizations in Algorithm \ref{alg:3} for all the examples, except for Example \ref{sec:ex1} where higher order BDF time steppings are also used to verify the space/time high-order accuracy of the proposed methods. 
For problems with shocks, we take $q_1=0.5$ and $q_2=2$ as the default choice of the artificial viscosity parameters in \eqref{av} unless otherwise stated.

We take the time step size as 
\begin{align}
\label{cfl}
\delta t =\min\{\text{CFL} \frac{h_{min}}{|\bld u_h|+c_s}\}, 
\end{align}
where the length scale $h_{min}=h_0\alpha_0/k$ with $\alpha_0$ being the minimal singular value of the Jacobian matrix $\nabla_X\bld x_h$.
Here the default choice of the CFL constant is taken to be  $\text{CFL}=1$ unless otherwise stated.
The automatic time-step control detailed in \cite[7.3]{Dobrev12} is also used for the examples with shocks.
The Newton's method is used to solve the nonlinear system for velocity DOFs in each time step, where the average iteration counts for all cases are observed to be around 4-8.
Most of the linearized systems in each Newton iteration are solved using the sparse Cholesky factorization, with the only exception of 
the low-Mach number cases in Example \ref{sec:gr} where a direct pardiso solver is used as the matrix failed to be positive definite therein. 

\subsection{Accuracy test: 2D Taylor-Green Vortex}
\label{sec:ex1}
We consider the invisid 2D Taylor-Green Vortex problem 
proposed in \cite{Dobrev12} to check the high-order 
convergence of our proposed algorithm on deforming domains.
Following \cite{Dobrev12}, we take 
$\gamma=5/3$ and 
add a source term 
\[
e_{\sf src}(\bld x, t) = 
\frac{3\pi}{8}(\cos (3\pi x)\cos(\pi y)-\cos (\pi x)\cos(3\pi y))
\]
to the internal energy equation \eqref{internal-e}.
The computational domain is a unit square with wall boundary conditions, 
and the initial conditions are taken such that the exact solutions are:
\begin{align*}
\rho(\bld x,t) =&\;1, \\
\bld u(\bld x, t)=&\;(\sin(\pi x)\cos(\pi y), -\cos(\pi x)\sin(\pi y)),\\
p(\bld x, t)=&\;\frac{1}{4}(\cos (2\pi x)+\cos(2\pi y))+1.
\end{align*}

We turn off artificial viscosity in the numerical simulations, and perform mesh convergence studies for the $L^2$-errors of velocity and internal energy at final time $t=0.5$ on a sequence of four consecutive uniform rectangular meshes with size $2^{3+l}\times 2^{3+l}$ for $l=0,1,2,3$.
We consider the  BDF$[m]$ scheme Algorithm~\ref{alg:3} with polynomials of degree $m$ used for the spatial discretization for $m=1,2,3,4$. The time step size is taken be $\delta t=0.05/2^l$ for $l=0,1,2,3$.
The midpoint rule Algorithm~\ref{alg:2} (with smaller time steps) is used to generate the starting values for high order BDF schemes.
History of convergence for the two $L^2$ errors 
at $t=0.5$ are recorded in Table~\ref{tab:1}.
We observe $m$-th order convergence for both variables for all $1\le m\le 4$. Hence, spatial and temporal high order convergence is achieved.
This convergence rate is optimal for the BDF time stepping, but suboptimal by one order for the spatial discretization. We find the convergence behavior of
the BDF$[m+1]$-$P^m$ scheme (not reported here for simplicity) is similar to BDF$[m]$-$P^m$ for $1\le m\le 3$, which suggests the one-order spatial convergence rate reduction is unavoidable for our scheme. 

We remark that our spatial discretization can be made slightly more efficient by taking polynomial degree one order lower for the thermodynamic variables than the 
flow map approximation, while still maintaining a similar convergence behavior. However, since the major computational cost for our scheme is in the nonlinear system solver in \eqref{u-Th1}. Such efficiency gain is not significant, and will not be investigated further in this work. 

\begin{table}[tb]
\centering
\begin{tabular}{ll|ll|ll}
\toprule
 $m$&  mesh   & $\|\bld u_h-\bld u\|_{\Omega}$  & order  &
 $\|e_h-e\|_{\Omega}$& order   \\ \midrule
& $8\times 8$ &  1.052e-01 & --&1.337e-01&--\\
& $16\times 16$ & 4.131e-02 & 1.35 & 7.284e-02 & 0.88\\
1& $32\times 32$ & 1.949e-02 & 1.08 & 3.697e-02 & 0.98\\
& $64\times 64$ & 9.710e-03 & 1.00 & 1.852e-02 & 1.00
\\\midrule
& $8\times 8$ & 1.077e-02 & --&1.263e-02&--\\
& $16\times 16$ & 3.250e-03 & 1.73 & 3.578e-03 & 1.82\\
2& $32\times 32$ & 7.933e-04 & 2.03 & 8.764e-04 & 2.03\\
& $64\times 64$ & 1.964e-04 & 2.01 & 2.192e-04 & 2.00\\
\midrule
& $8\times 8$ &  5.590e-03 & --&2.947e-03&--\\
& $16\times 16$ &  5.809e-04 & 3.27 & 3.175e-04 & 3.21\\
3& $32\times 32$ &  7.070e-05 & 3.04 & 4.181e-05 & 2.92\\
& $64\times 64$ &  8.766e-06 & 3.01 & 5.404e-06 & 2.95\\
\midrule
& $8\times 8$ &  3.986e-03 & --&2.510e-03&--\\
& $16\times 16$ &  1.361e-04 & 4.87 & 7.305e-05 & 5.10\\
4& $32\times 32$ &  5.690e-06 & 4.58 & 2.928e-06 & 4.64\\
& $64\times 64$ &  4.755e-07 & 3.58 & 3.304e-07 & 3.15\\
\bottomrule
\end{tabular}
\vspace{.3ex}
\caption{Example \ref{sec:ex1}: History of convergence 
for the BDF$[m]$-$P^m$ scheme for $1\le m\le 4$.}
\label{tab:1}
\end{table}

\subsection{1D Shock Tube}
\label{sec:ex2}
\begin{figure}[ht!]
    \centering
    \includegraphics[width = .9\textwidth]{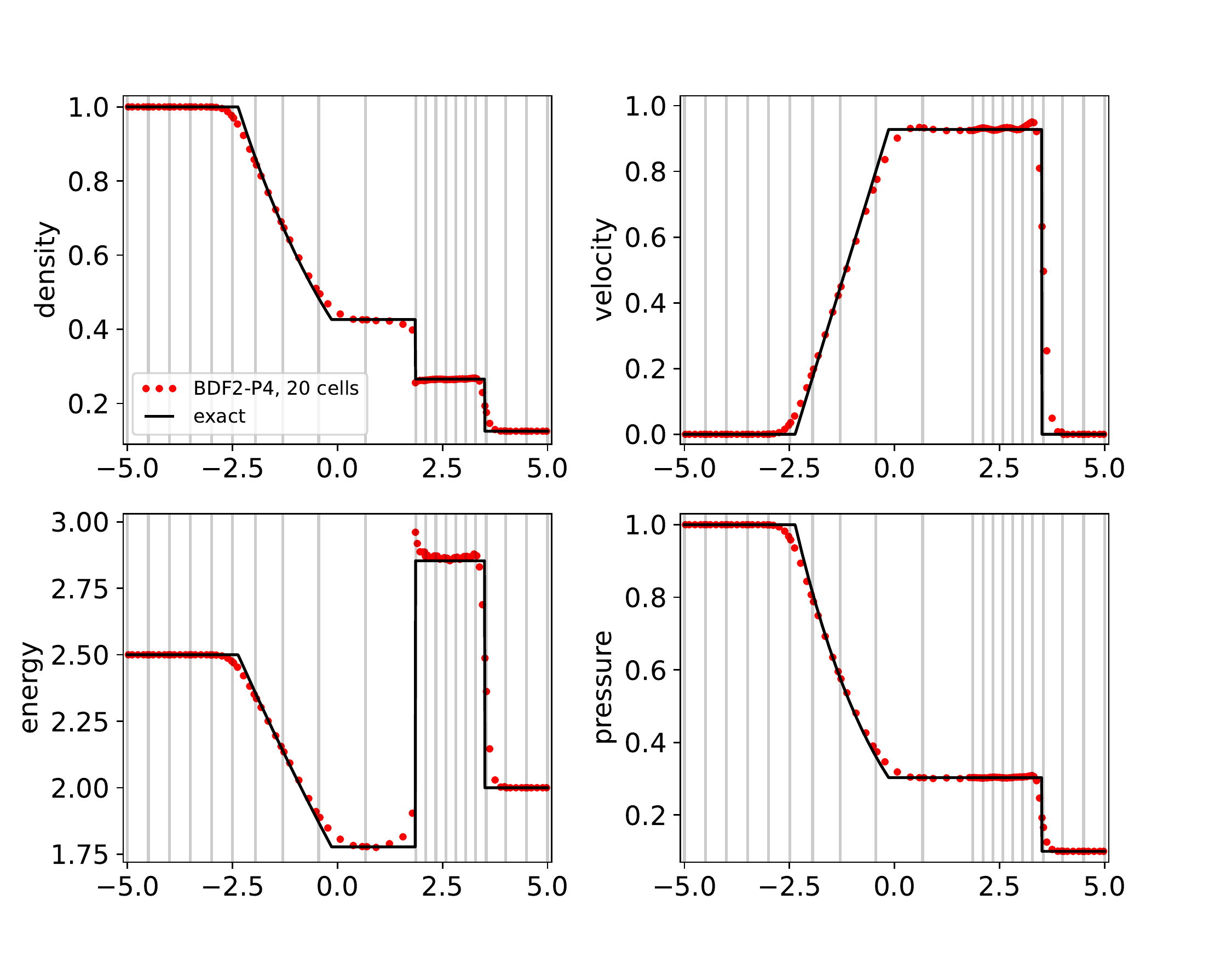}
    \vspace{-5ex}
    \caption{Example \ref{sec:ex2}: Results 
     at final time $t=0.2$ sampled on 100 quadrature points.
    The gray lines are the deformed cell boundaries.
    }
    \label{fig:1}
\end{figure}

We consider a simple 1D Riemann problem, the Sod shock tube on the domain $\Omega=[-5,5]$ with initial condition
\[
(\rho, \bld u, p)=(1,0,1) \text{ for } x\in[-5,0], 
\quad
(\rho, \bld u, p)=(0.125,0,0.1) \text{ for } x\in[0,5].
\]
Here $\gamma=1.4$.
We apply the BDF2 scheme with polynomial degree $k=4$ on an initial mesh with 20 uniform cells. 
The results at final time $t=2.0$ 
on all $20\times (4+1) =100$ quadrature points
are shown in Figure~\ref{fig:1}, along with the deformed cells. 
The numerical approximation agrees with the exact Riemann solution quite well even on this coarse mesh, where the shock is resolved within 2 cells. 
As typical of Lagrangian schemes, the contact discontinuity  is captured without dissipation.
We also observe the ``wall heating" phenomenon in the internal energy at the contact.

\subsection{2D Sedov explosion}
\label{sec:sedov}
The Sedov explosion \cite{Sedov} models the expanding wave by an intense explosion in a perfect gas. It is a standard problem to test the ability of codes to preserve the radial symmetry of shocks. 
The domain is a square $\Omega=[0, 1.2]\times [0, 1.2]$.
The initial condition is set to have unit density and zero velocity, and also to have zero internal energy except at the left bottom corner cell $T_0$, where its is a bilinear function whos
value is $\frac{0.2448\times 4}{\mathrm{area}(T_0)}$ on the left/bottom corner vertex and zero on the other three vertices. 
So the total initial internal energy is 
$
\int_{\Omega_0}\rho e\mathrm{dX} = 0.2448.
$
Symmetry boundary conditions are imposed on the left and bottom boundaries, while free boundary conditions are used on the top and right boundaries.
The analytic solution at time $t=1$ gives a shock at radius $r=1$ with a  peak density of $6$.

\begin{figure}[ht!]
    \centering
\subfigure[$k=4$, $N=16$]{\includegraphics[width = .27\textwidth]{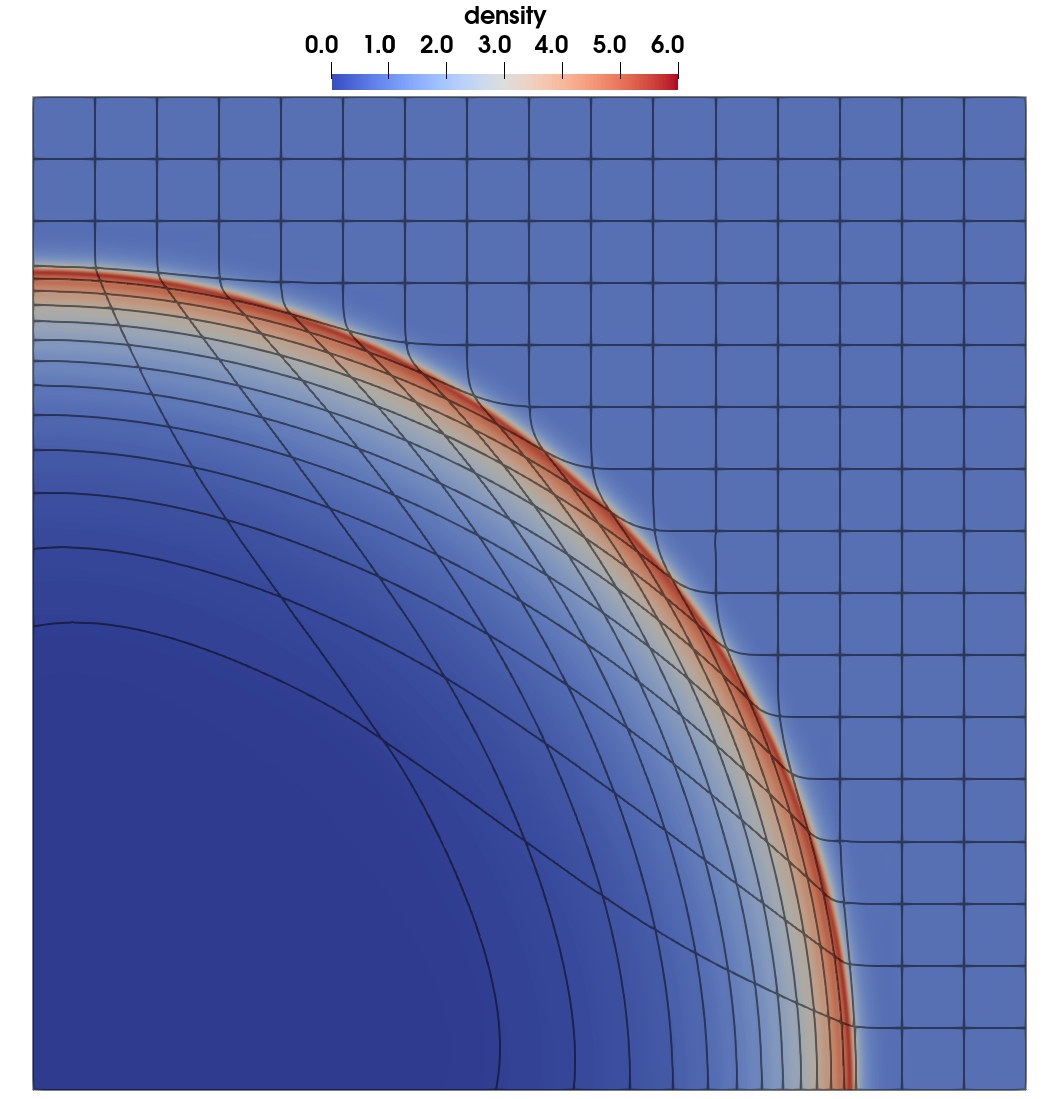}}
\subfigure[$k=4$, $N=32$]{\includegraphics[width = .27\textwidth]{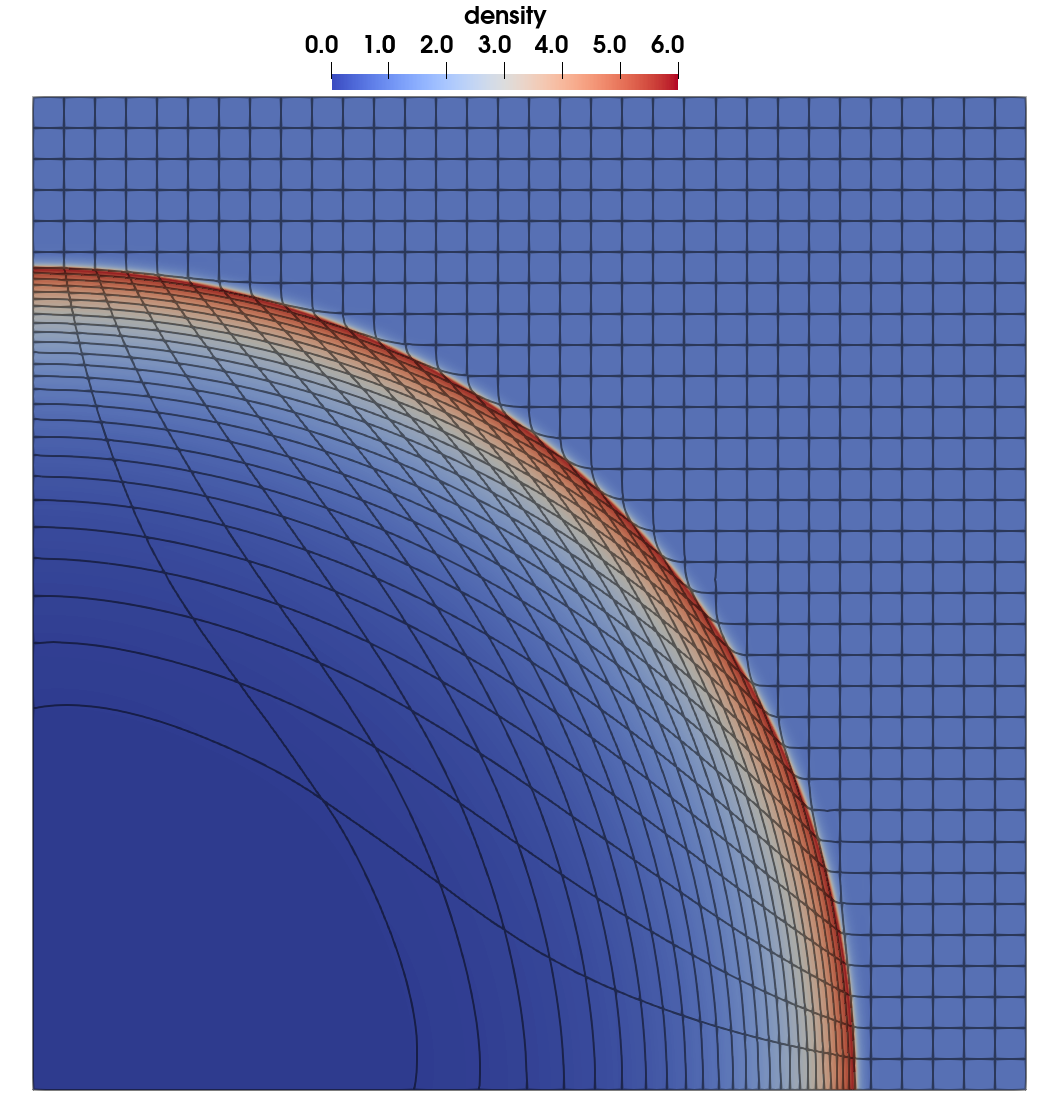}}
\subfigure[den. v.s. rad.]{\includegraphics[width = .4\textwidth]{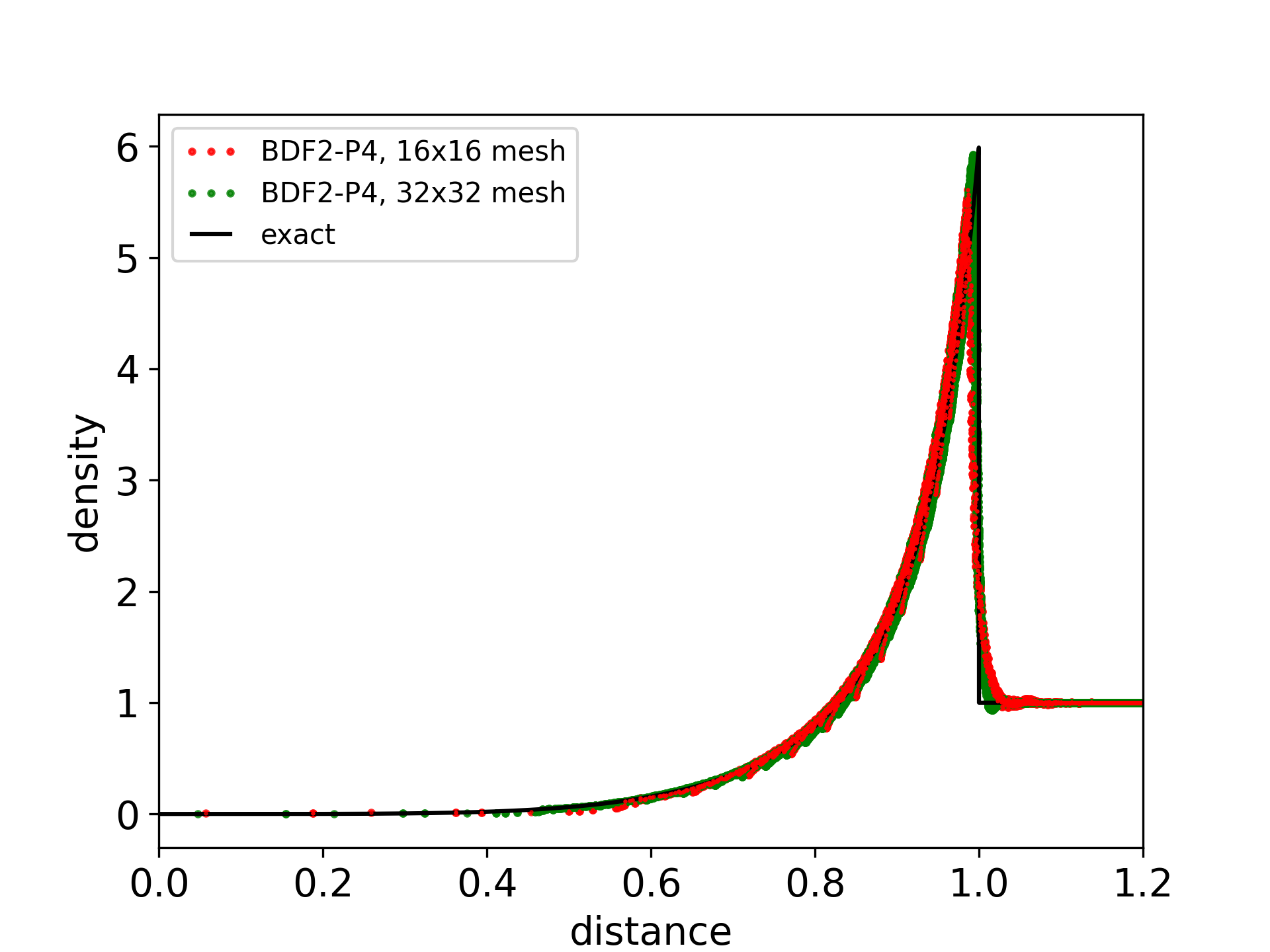}}
    \caption{Example \ref{sec:sedov}:  (a)-(b): Density contour on deformed domain at final time $t=1.0$ using BDF2 time stepping with polynomial degree $k=4$ on rectangular meshes with size $N\times N$.
(c): Scattered plot of density v.s. radius on all quadrature points at final time $t=1$.
     }
    \label{fig:sedov}
\end{figure}

We apply the BDF2 scheme with polynomial degree
$k=4$ on initial uniform rectangular meshes of size $N\times N$ with $N=16$ and $N=32$.
The density field on deformed meshes, 
and also the scattered plot of density
versus radius $r=\sqrt{x^2+y^2}$ on all quadrature points 
are shown in Figure~\ref{fig:sedov}.
The radial symmetry of the solution is preserved, and the numerical solution agrees quite well with the analytic solution in Figure~\ref{fig:sedov}(c). 


\subsection{2D Noh explosion}
\label{sec:noh}
The Noh explosion problem \cite{Noh} 
consists of an ideal gas with $\gamma =5/3$, initial density $\rho_0=1$, initial internal energy $e_0=0$, and initial velocity 
$\bld u_0=(-\frac{x}{\sqrt{x^2+y^2}},-\frac{y}{\sqrt{x^2+y^2}})$. 
The computational domain is $\Omega=[0,1]\times[0,1]$. We use symmetry boundary conditions on left and bottom boundaries, and 
free boundary conditions on top and right boundaries.
Similar to the previous case, this problem has a radial symmetry, and the analytic solution at time $t=0.6$ gives a shock at radius $r=0.2$ with a  peak density of $16$.

\begin{figure}[ht!]
    \centering
\subfigure[$k=4$, $N=16$]{\includegraphics[width = .27\textwidth]{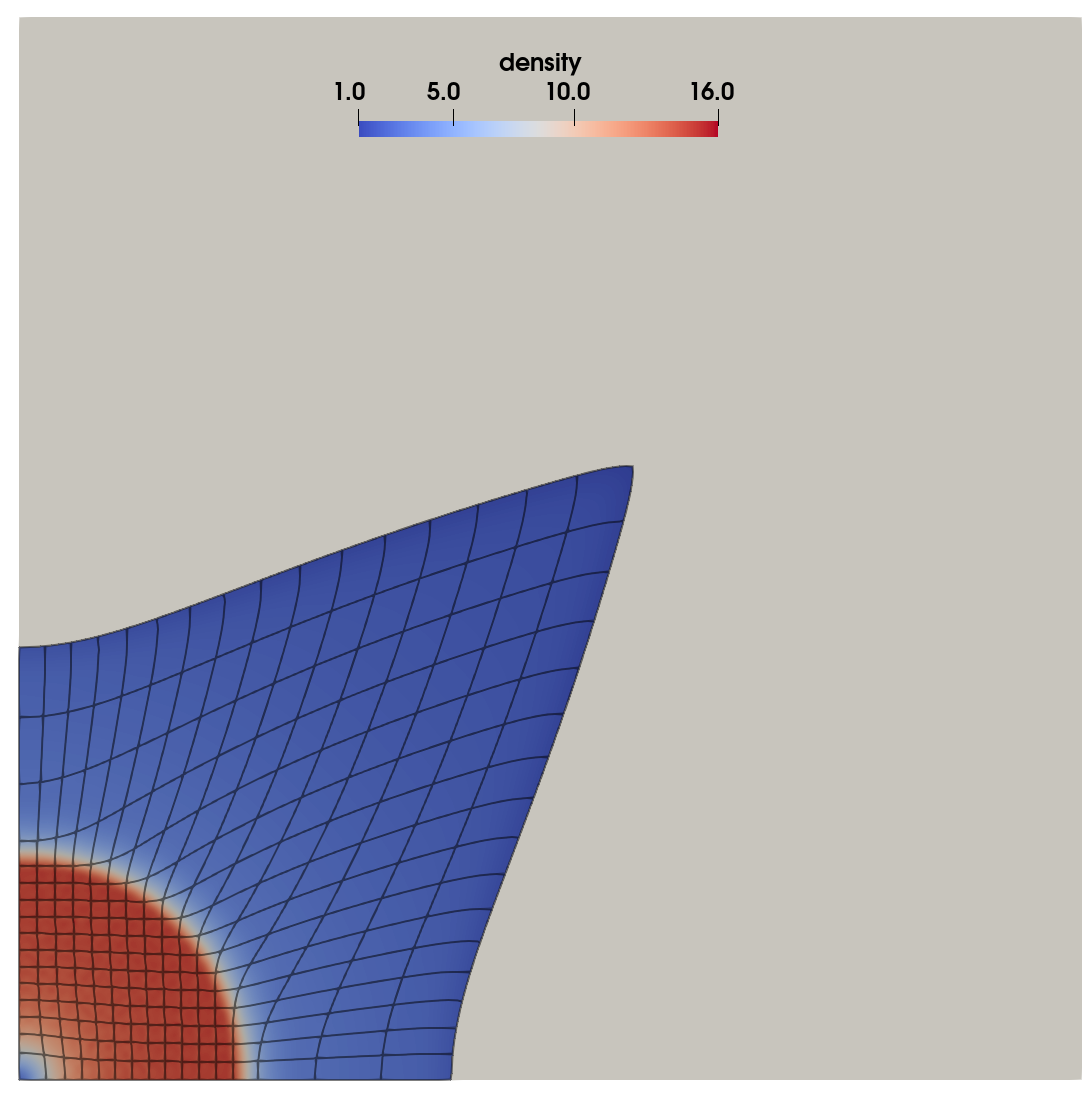}}
\subfigure[$k=4$, $N=32$]{\includegraphics[width = .27\textwidth]{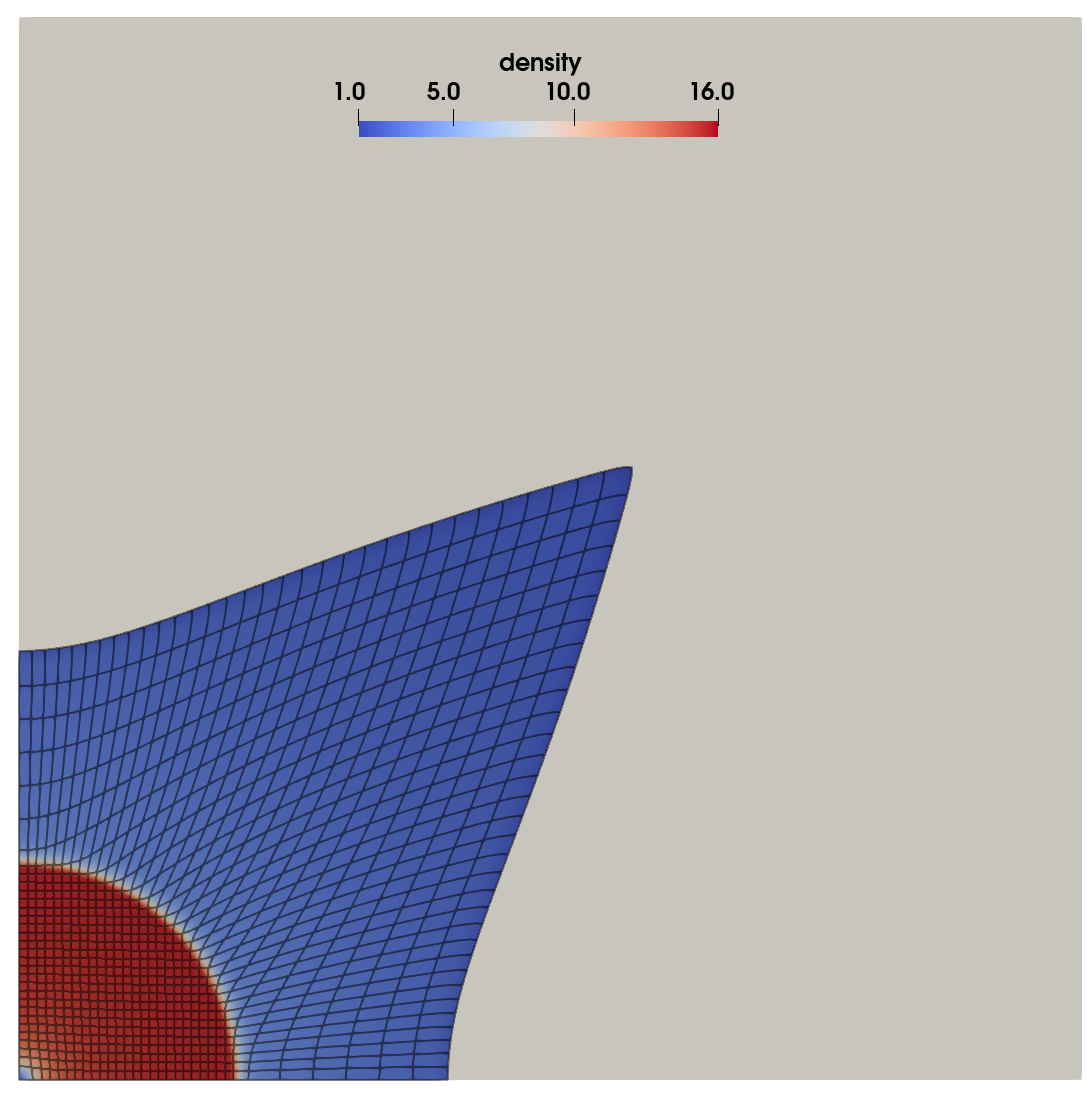}}
\subfigure[den. v.s. rad.]{\includegraphics[width = .4\textwidth]{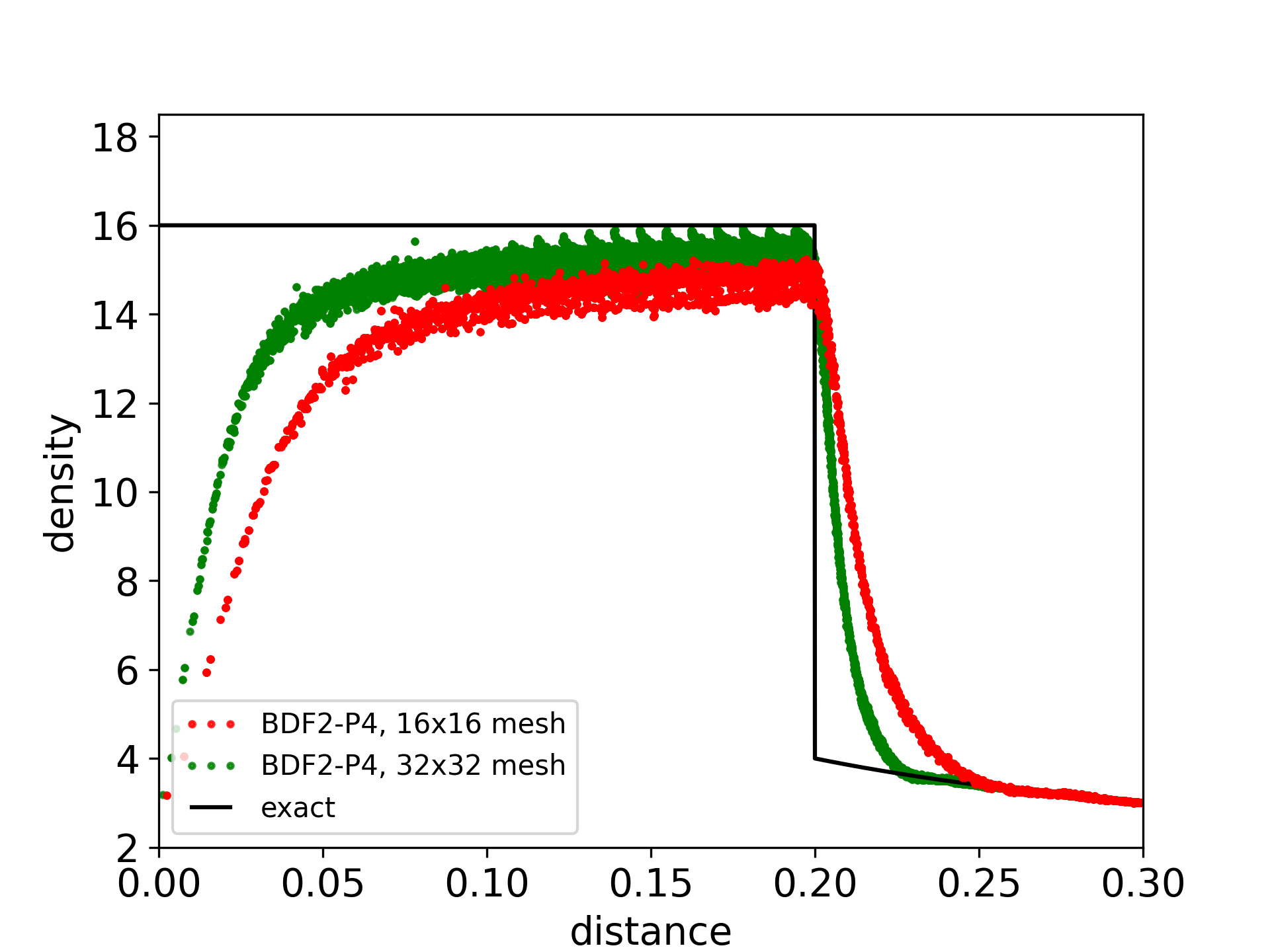}}
\caption{Example \ref{sec:noh}:  (a)-(b): Density contour on deformed domain at final time $t=0.6$ using BDF2 time stepping with polynomial degree $k=4$ on rectangular meshes with size $N\times N$.
Here the gray square is the initial domain.
(c): Scattered plot of density v.s. radius on all quadrature points at final time $t=0.6$.
     }
    \label{fig:noh}
\end{figure}

We apply the BDF2 scheme with polynomial degree
$k=4$ on initial uniform rectangular meshes of size $N\times N$ with $N=16$ and $N=32$.
For this problem, we observe the default choice of artificial viscosity coefficients with $q_1=0.5$ and $q_2 = 2$ leads to quite large post-shock oscillations.
So we increase these coefficients to $q_1=1$ and $q_2 = 4$ in our numerical experiments reported here.
The radial symmetry of the solution is preserved as seen in (a) and (b) of Figure~\ref{fig:noh}, and the numerical solution has a good agreement with the analytic solution for $0.1<r<0.3$ in Figure~\ref{fig:sedov}(c), although it is slightly oscillatory.

\subsection{Triple point problem}
\label{sec:t3}
\begin{figure}[ht!]
    \centering
    \includegraphics[width = .45\textwidth]{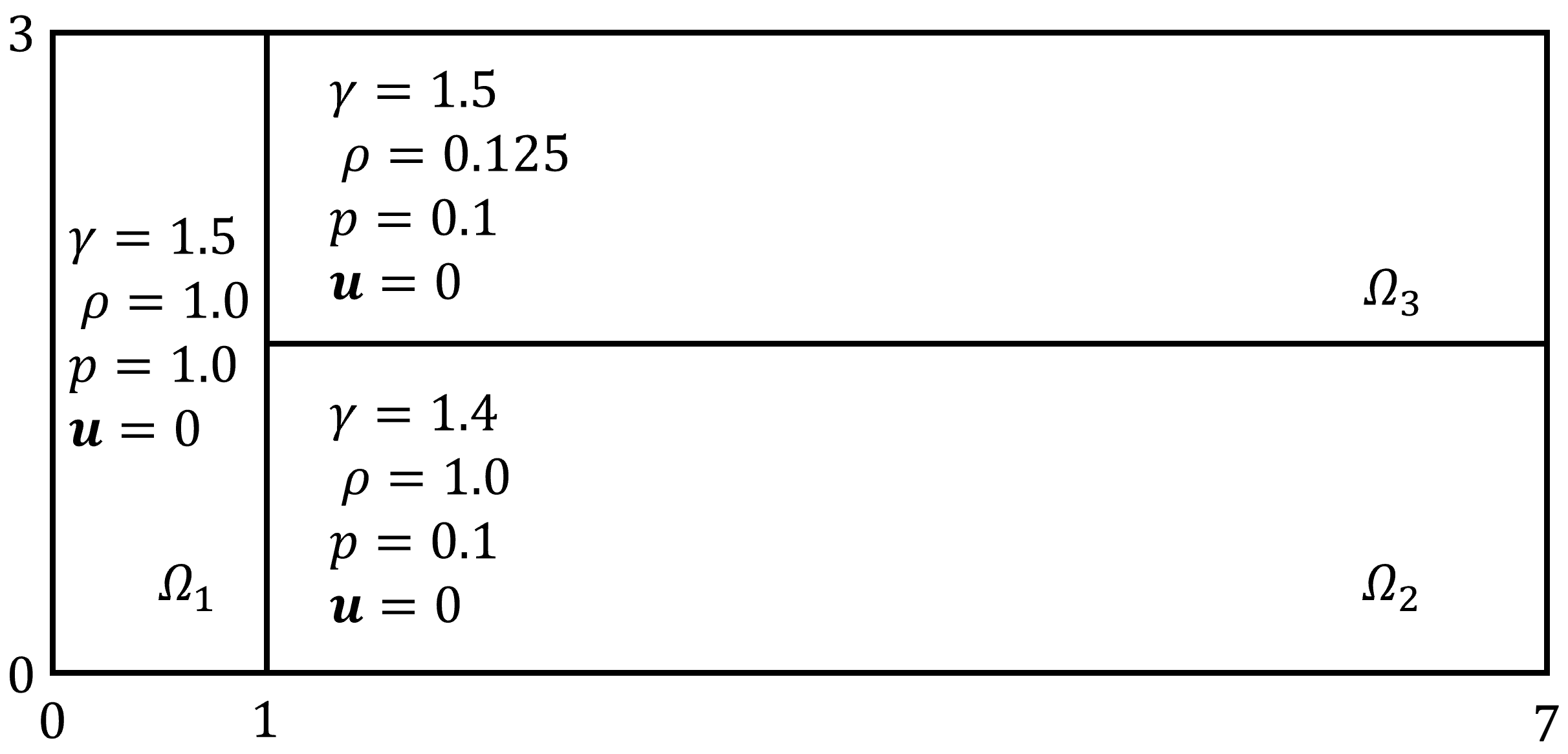}
    \caption{Example \ref{sec:t3}: Initial data for the triple point problem.
     }
    \label{fig:t3}
\end{figure}

The triple point problem is a multimaterial test case proposed in \cite{Kucharik10}; see also \cite{Galera10}. 
 The initial data is shown in Fig.~\ref{fig:t3}.
The computational domain has a rectangular shape with $7\times 3$ edge ratio. 
It includes three materials  at rest located in 
$\Omega_1=[0,1]\times[0,3]$, 
$\Omega_2=[1,7]\times[0,1.5]$, and
$\Omega_3=[1,7]\times[1.5,3]$, initially forming a T-junction. The high-pressure material in $\Omega_1$ creates a shock wave moving to the right. Due to different material properties in $\Omega_2$ and $\Omega_3$, the shock wave moves faster in 
$\Omega_3$, which leads to vortex formation around the triple point.
For pure Lagrangian methods, there is a limit to how long this problem can be run due to vortex generation. 
Similar to \cite{Dobrev12}, we run simulation till time $t=3.3$.

We apply the BDF2 scheme with polynomial degree
$k=4$ and on initial uniform rectangular meshes of size $28\times 12$ and $56\times 24$.
Here we reduce the artificial viscosity coefficients to be $q_1=0.25$ and $q_2=1$, and take a larger CFL number with $\text{CFL}=3$. 
A total of 238 time steps is used to drive the solution to final time $t=3.3$ on the fine mesh with $56\times 28$ cells. 
We note that an explicit scheme would require about two order of magnitude more time steps for a stable simulation. For example, the Laghos code freely available in the github repository \url{https://github.com/CEED/Laghos}, which implements the high-order Lagrangian finite element scheme in \cite{Dobrev12},
requires about $90,000$ time steps for $k=4$ on the fine mesh with RK4 time stepping and CFL=1.
We further note that the default choice of parameters also work for this problem; we make the modifications to illustrate that the high-order method is still robust with a larger time step size and a smaller artificial viscosity. 

\begin{figure}[ht!]
    \centering
    \includegraphics[width = .45\textwidth]{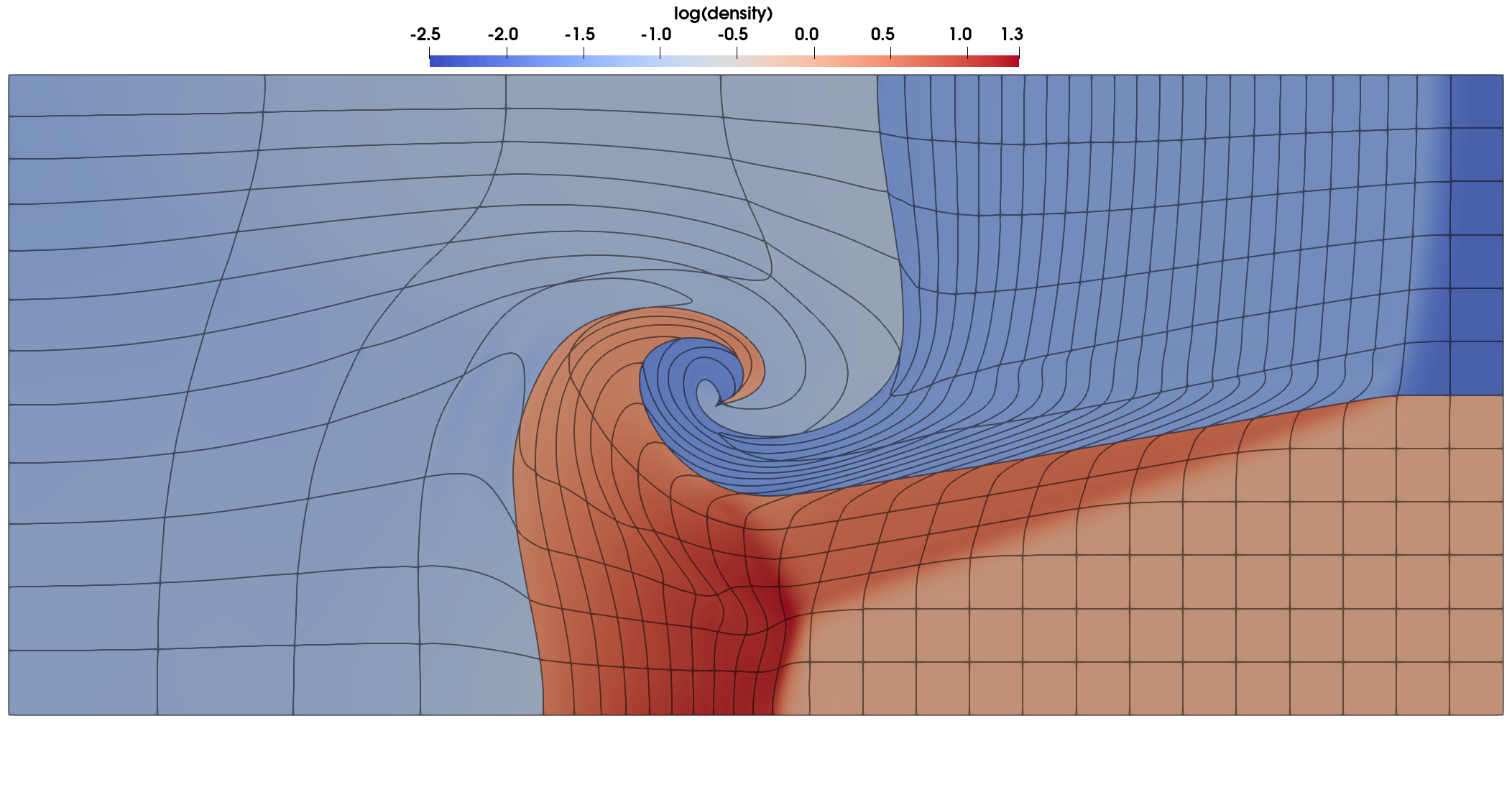}
    \includegraphics[width = .45\textwidth]{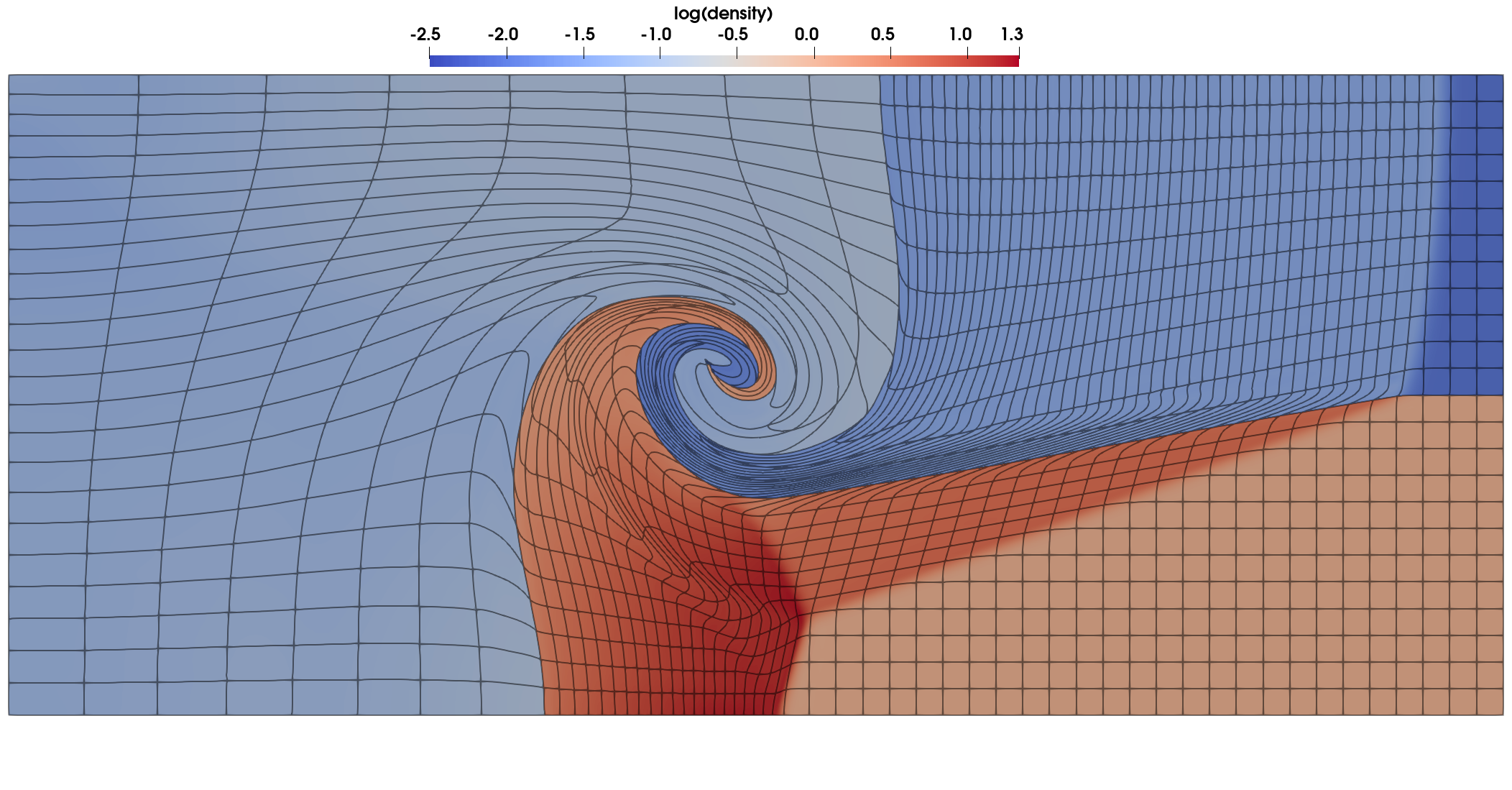}
    \caption{Example \ref{sec:t3}: Results 
     at final time $t=3.3$.
     }
    \label{fig:t32}
\end{figure}

The density plot in log-scale on the deformed meshes
are shown in Figure~\ref{fig:t32}. We observe the material interfaces are sharply preserved as typical of Lagrangian schemes, and the the shock locations are essentially the same, but the total amount of "roll-up" at the triple point increases as the mesh refines.

\subsection{Gresho vortex}
\label{sec:gr}
The Gresho vortex problem \cite{Gresho88} is an example of a stationary, incompressible rotating flow around the origin in two spatial dimensions, where centrifugal forces are exactly balanced by pressure gradients. It was first applied to the compressible Euler equations in \cite{Liska03}.
Here we use the low-Mach setup given in \cite{Miczek}.
The domain is a unit square $\Omega=(0,1)\times(0,1)$ with wall boundary conditions, $\gamma=1.4$, 
and the initial conditions are
\begin{align*}
\rho_0 = 1, \bld u_0=\;
u_{\phi}\bld e_{\phi},
p_0 = \;\frac{1}{\gamma M_{\max}^2}-\frac12
+ &\begin{cases}
12.5r^2& \text{if } r<0.2,\\
4\mathrm{ln}(5r)+4-20r+12.5r^2& \text{if } 0.2\ge r\ge 0.4,\\
4\mathrm{ln}2-2&\text{if } r>0.4,
\end{cases}
\end{align*}
where  the radius $r = \sqrt{(x-0.5)^2+(y-0.5)^2}$, 
angular velocity $u_{\phi}$ is 
\[
u_{\phi}(r)=
\begin{cases}
5r& \text{if } r<0.2,\\
2-5r& \text{if } 0.2\ge r\ge 0.4,\\
0&\text{if } r>0.4,\\
\end{cases}
\]
the unit vector $\bld e_{\phi}=(-(y-0.5)/r, (x-0.5)/r)$, and $M_{\max}$ is the parameter used to adjust the maximum Mach number of the problem, where $
M(r) = \frac{|\bld u|}{\sqrt{\gamma p/\rho}}
$
is the Mach number. The maximum number of the initial condition is achieved at $r=0.2$ with $M(0.2) = M_{\max}$.

\begin{figure}[ht!]
    \centering
\subfigure[$k=2$, $N=40$, $M_{\max}=0.1$]{\includegraphics[width = .3\textwidth]{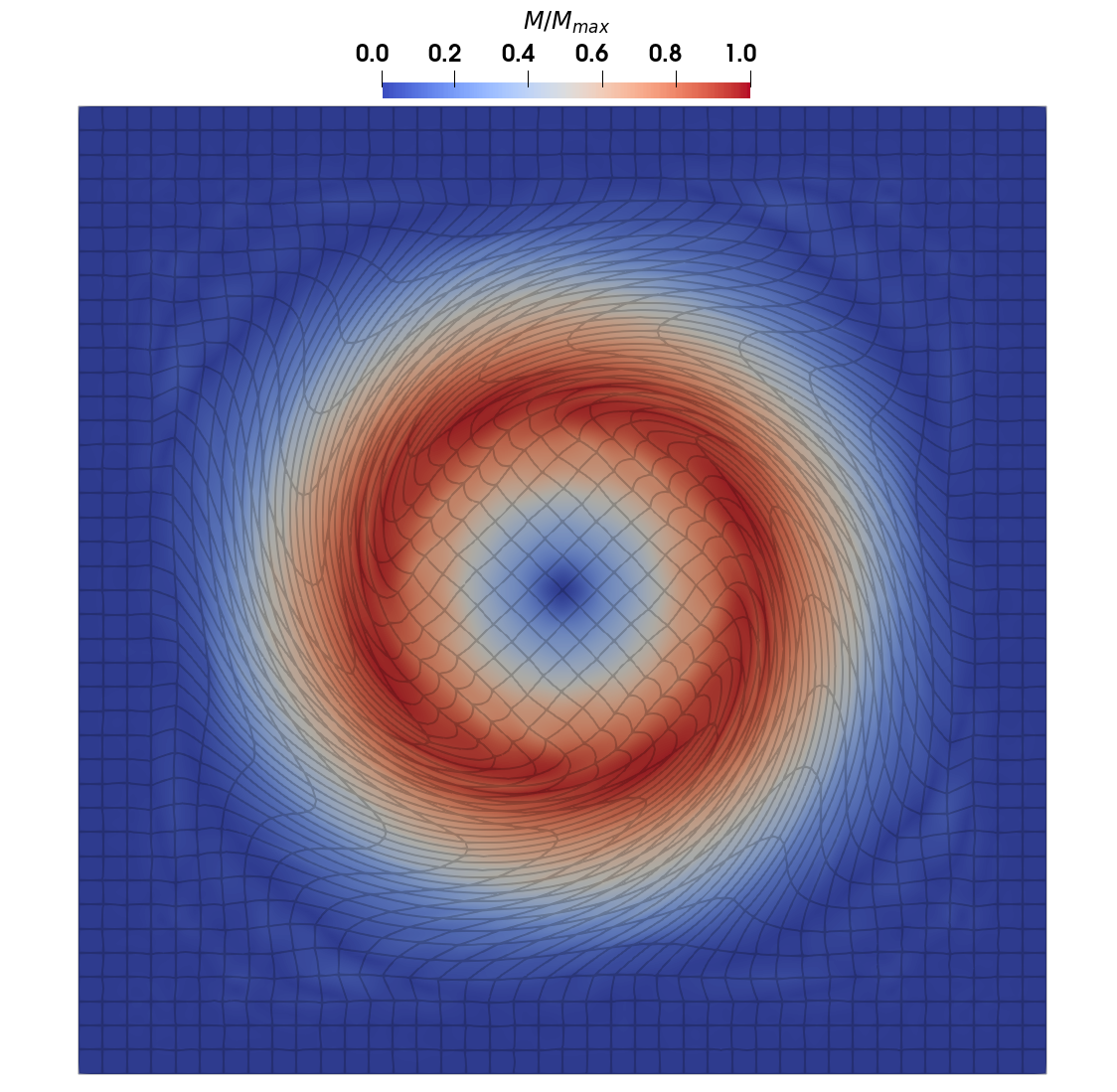}}
\subfigure[$k=2$, $N=40$, $M_{\max}=0.01$]{\includegraphics[width = .3\textwidth]{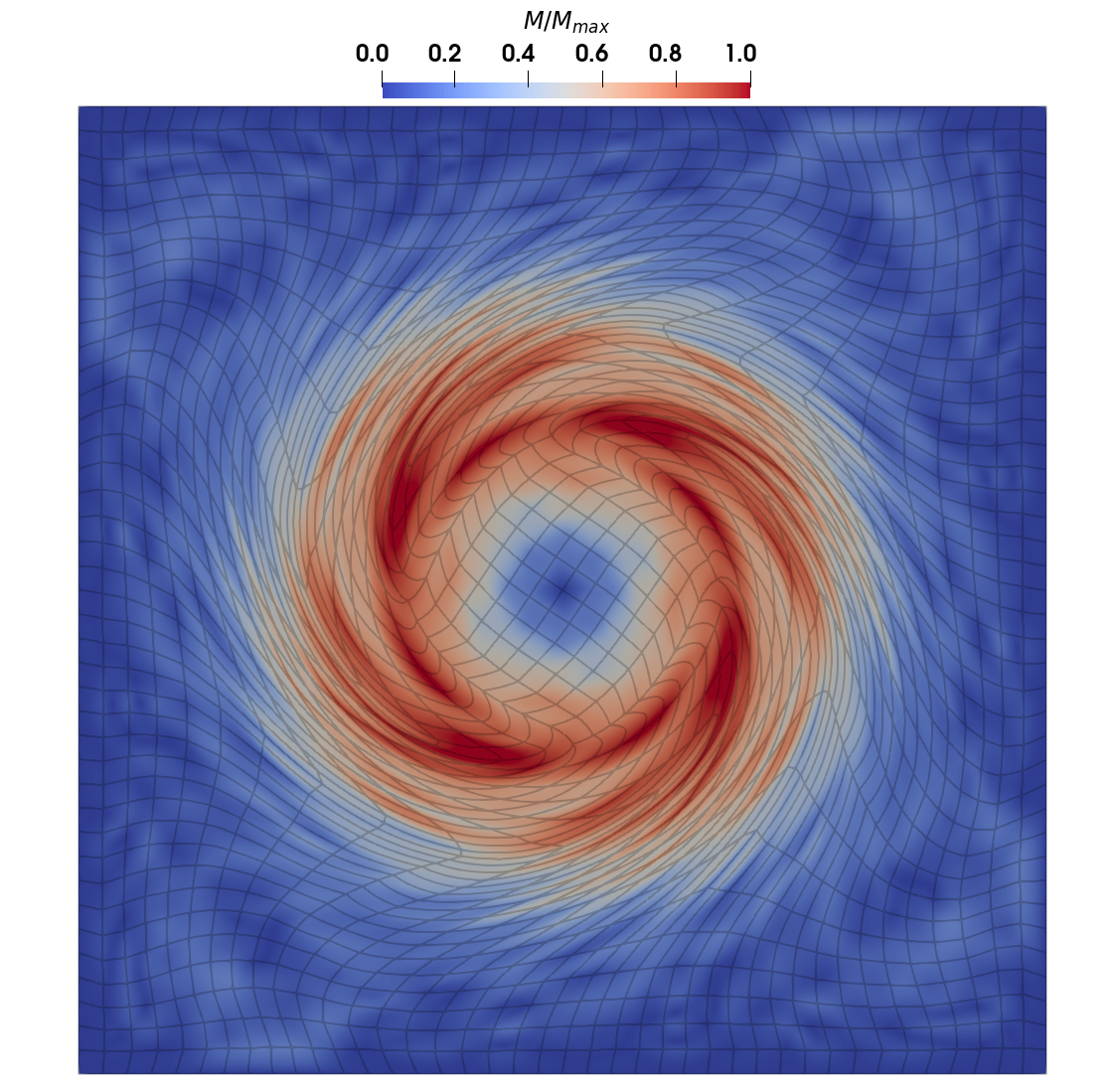}}
\subfigure[$k=2$, $N=40$, $M_{\max}=0.001$]{\includegraphics[width = .3\textwidth]{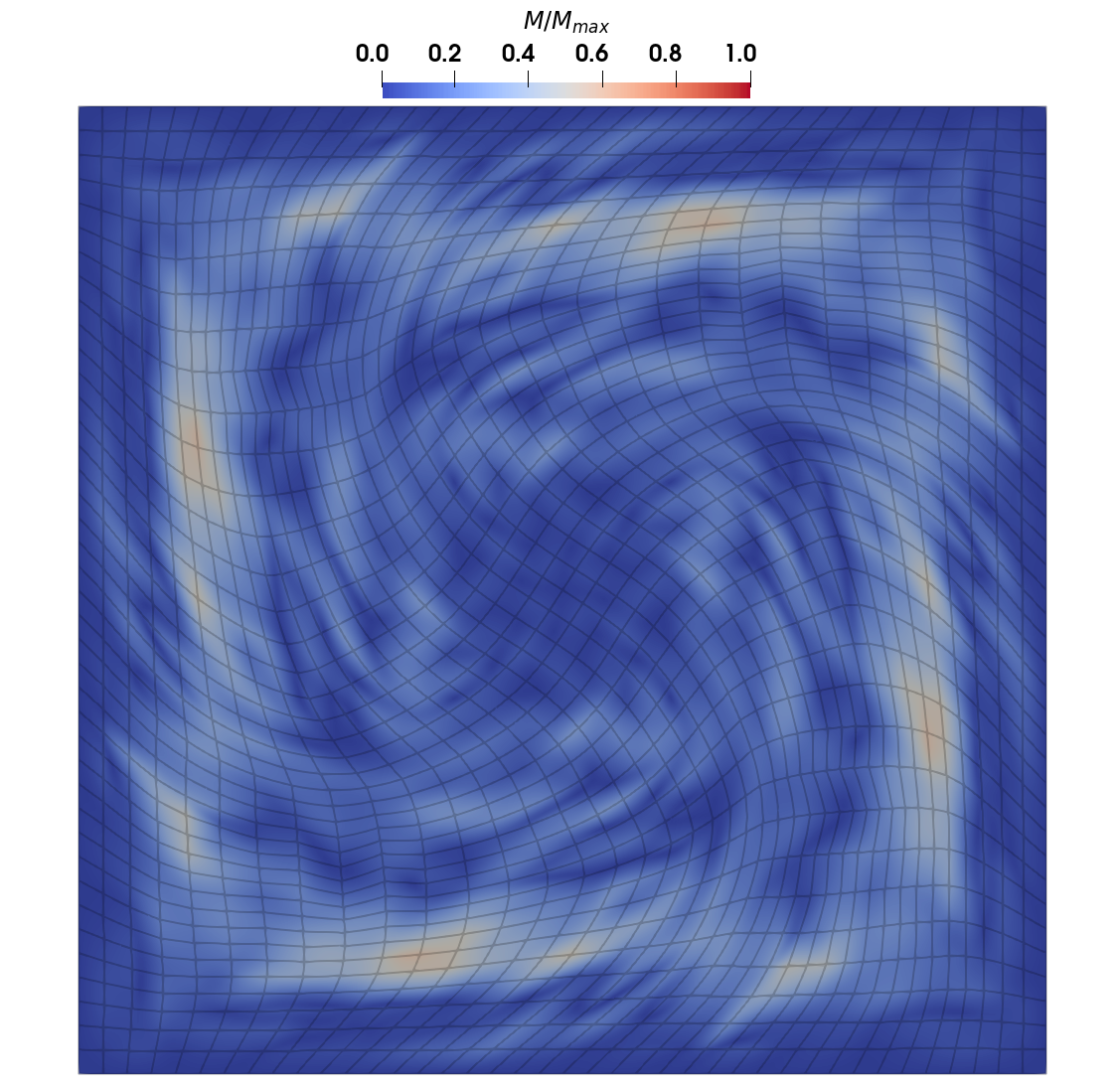}}
\subfigure[$k=4$, $N=20$, $M_{\max}=0.1$]{\includegraphics[width = .3\textwidth]{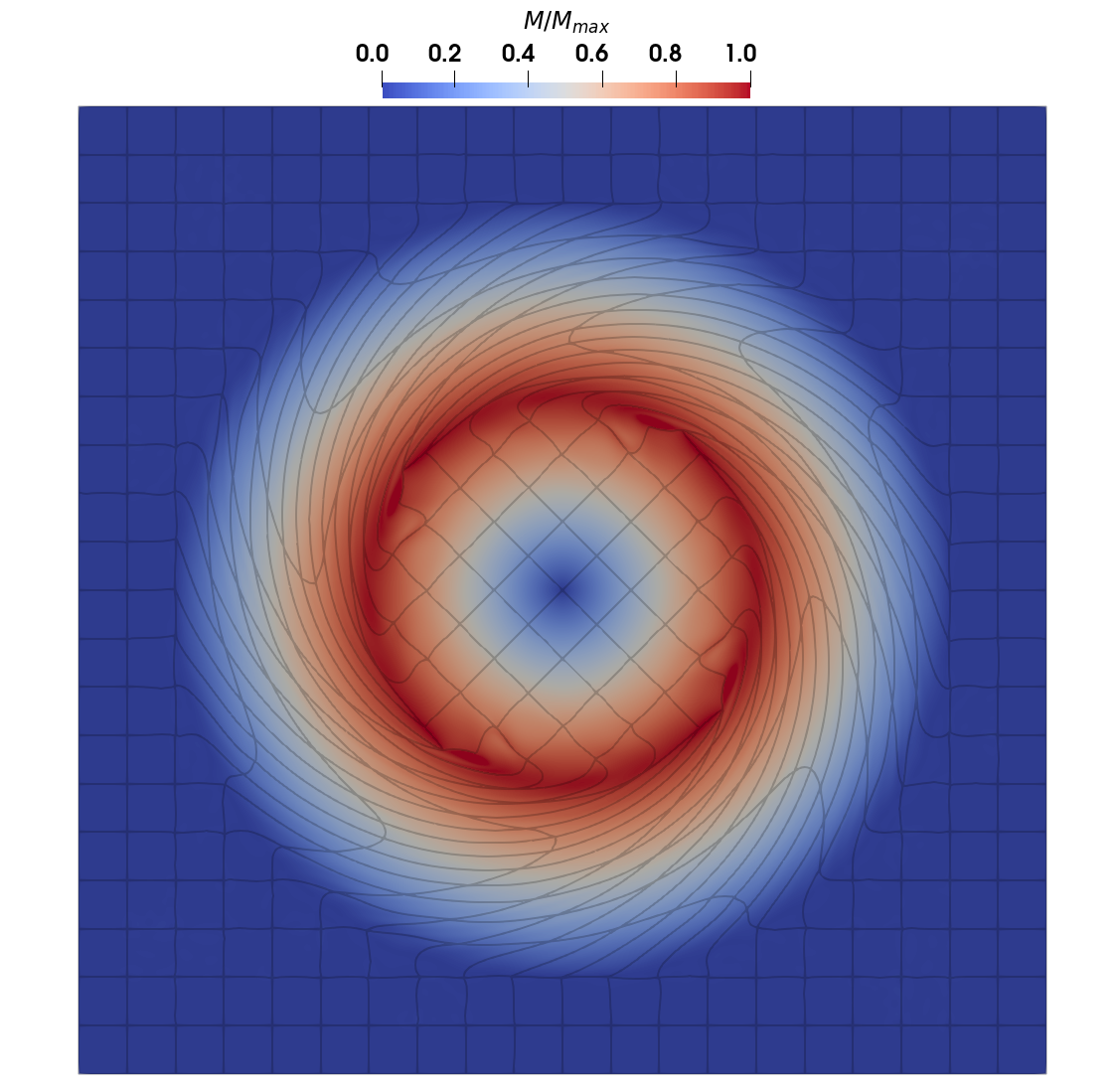}}
\subfigure[$k=4$, $N=20$, $M_{\max}=0.01$]{\includegraphics[width = .3\textwidth]{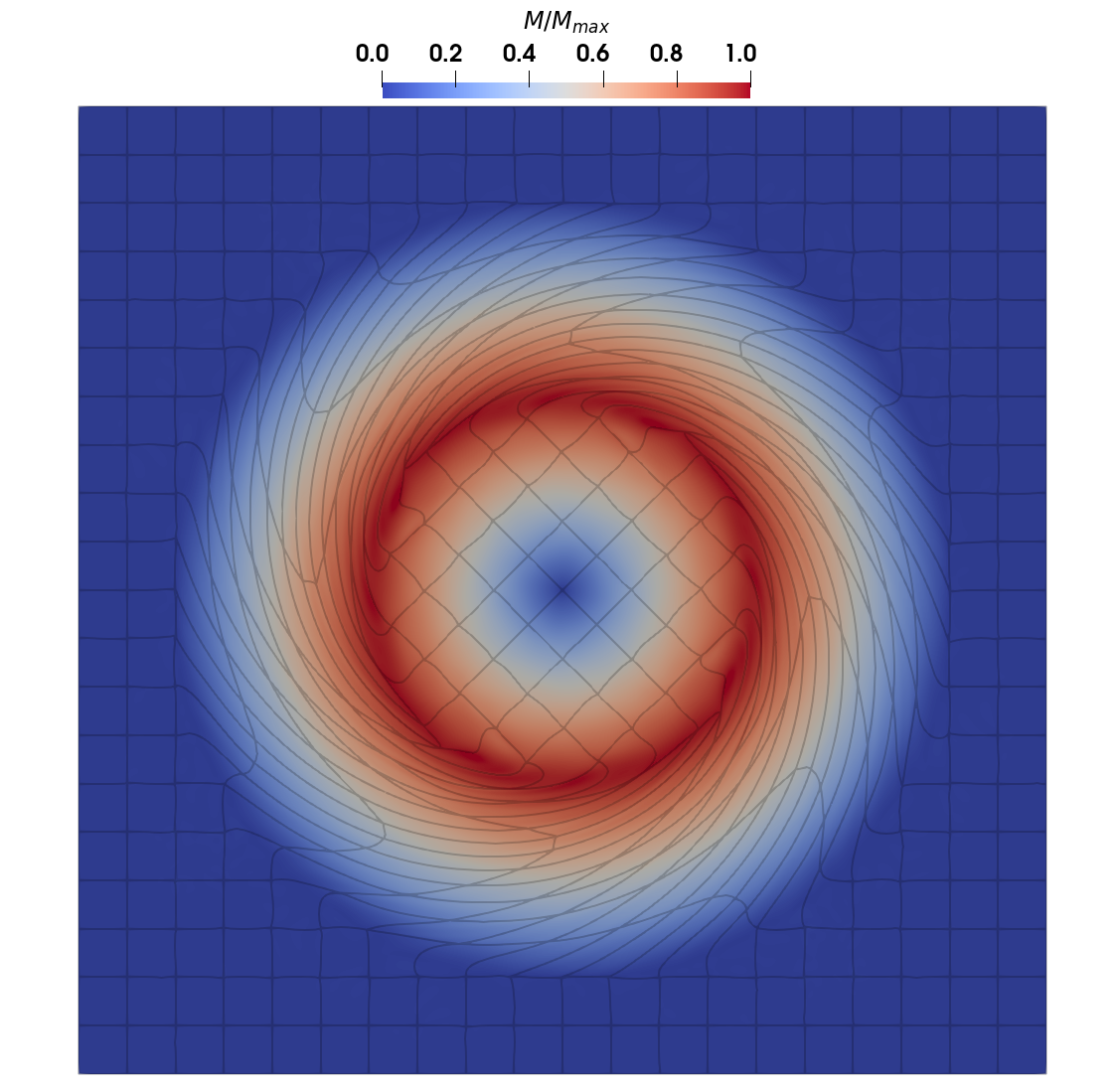}}
\subfigure[$k=4$, $N=20$, $M_{\max}=0.001$]{\includegraphics[width = .3\textwidth]{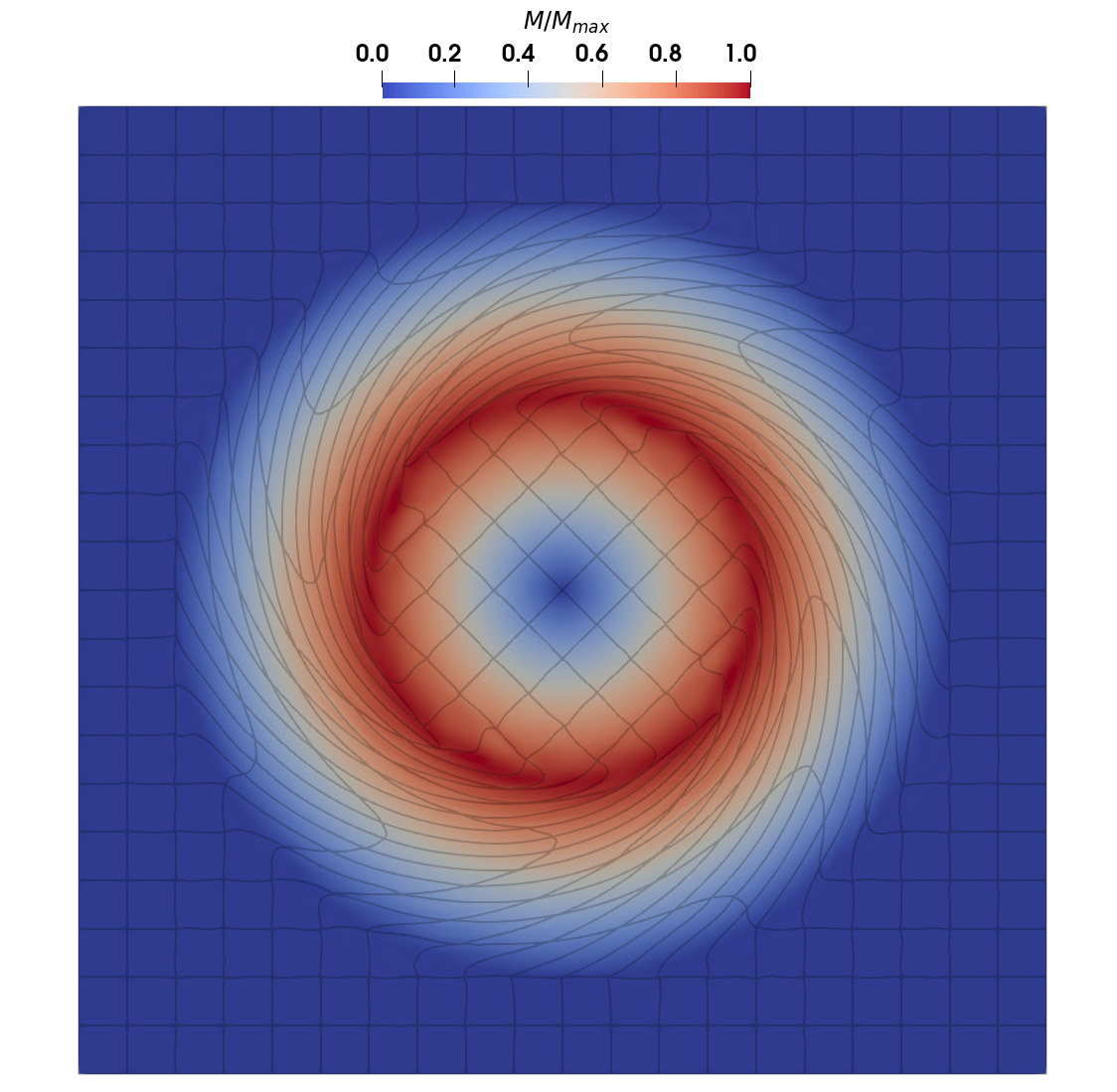}}
    \caption{Example \ref{sec:gr}: Density contour on deformed domain 
     at final time $t=1.0$ using BDF2 time stepping and 
     Algorithm \ref{alg:1}  with polynomial degree $k$ on a rectangular mesh with size $N\times N$ for different choices of $k$ and $N$.
     (parameters: CFL = 1.0, $q_1=0.5$, $q_2=2.0$)
     }
    \label{fig:gr}
\end{figure}

We apply the BDF2 scheme with polynomial degree $k=2$
on a $40\times 40$ mesh,
and $k=4$ on a $20\times 20$ mesh. 
The total number of velocity DOFs are the same for the two cases.
Since the problem is smooth, we turn off the artificial viscosity. 
The period of one rotation for $r=0.2$ is $2\pi r=0.4\pi$. We run simulation till final time $t = \frac34\times {0.4\pi}$ so that the internal flow has rotated $\frac{3}{4}\times 180= 135$ degrees.
The default choice of time step size \eqref{cfl}
is linear proportional to the  Mach number  since 
$\max c_s \approx 1/M_{\max}$. 
We remove this Mach-number dependency on time step size by multiplying the sound speed in \eqref{cfl} with the maximum Mach number, i.e., 
\[
\delta t =\min\{\text{CFL} \frac{h_{min}}{|\bld u_h|+M_{\max}c_s}\}.
\]
We use $M_{\max}=0.1, 0.01, 0.001$, and take the CFL number to be CFL=0.25
 for all cases.
The relative Mach number $M/M_{\max}$ at final time 
$t=\frac34\times 0.4\pi$ for all cases are shown in Figure~\ref{fig:gr}. 
The first row of Figure~\ref{fig:gr} show results for the $k=2$ simulations, where we clearly observe a locking phenomena as the Mach number decreases.
On the other hand, the higher order simulations with $k=4$ leads to almost identical results
for all three Mach numbers. This example illustrates the advantage of using a higher order scheme over a low-order scheme in the low-Mach number regime.

Moreover, the total number of time steps
for $k=4$ are between 900 to 1000 for all three cases. If an explicit scheme (e.g., RK3) were to be used to solve this problem, the total number of time steps will be about three orders of magnitude larger when $M_{\max}=0.001$ due to sound speed based CFL constraints.  
Finally, we note that the nonlinear system in each time step is harder to solve as the Mach number
$M_{\max}$ decreases.
For example, Newton's method with sparse Cholesky direct solver works for $M_{\max}=0.1$ with CFL=1, but it fails for $M_{\max}=0.001$, where we have to reduce CFL to be 0.25 and replace the Cholesky solver by a pardiso solver, which indicates the linear system for $M_{\max}=0.001$ in the Newton iteration is no longer positive definite.
The linear system solver issue in the low Mach number regime will be further investigate in our future work.

\subsection{Shock-bubble interaction}
\label{sec:sb}
\begin{figure}[ht!]
    \centering
\subfigure[Geometry]{\includegraphics[width = .55\textwidth]{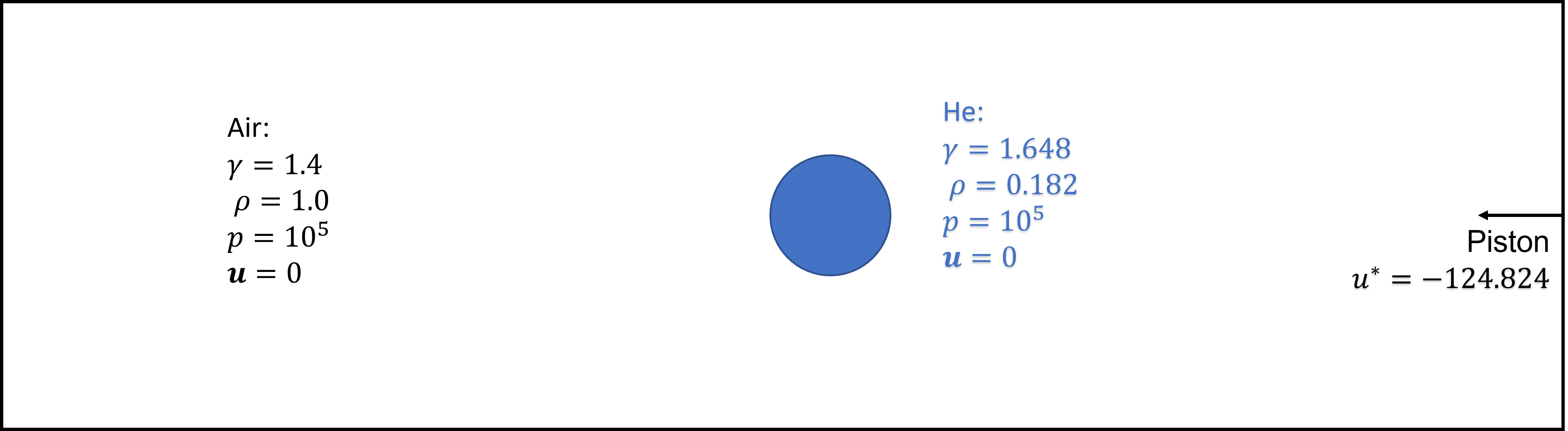}}
\subfigure[Coarse mesh]{\includegraphics[width = .18\textwidth]{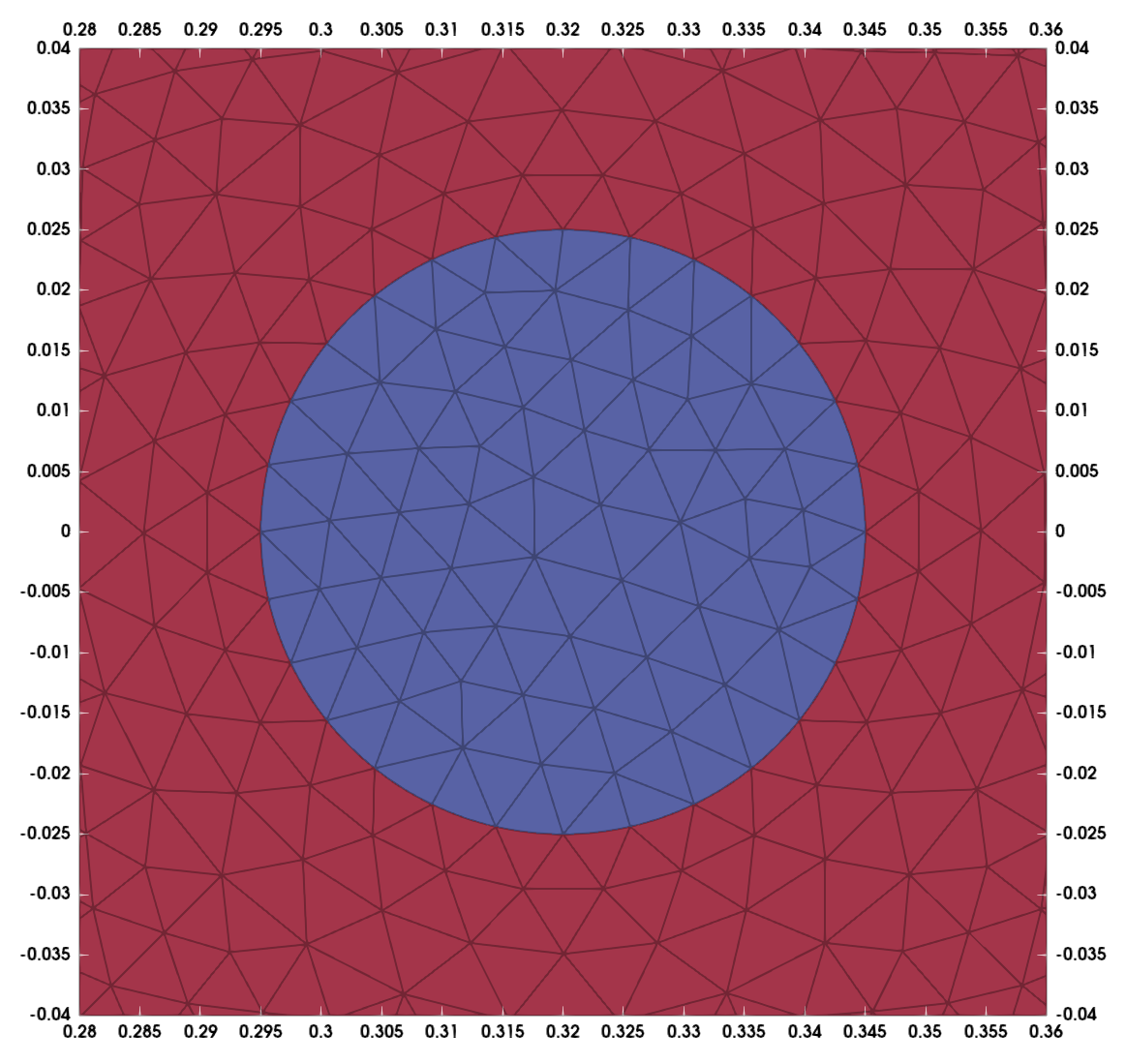}}
\subfigure[Fine mesh]{\includegraphics[width = .18\textwidth]{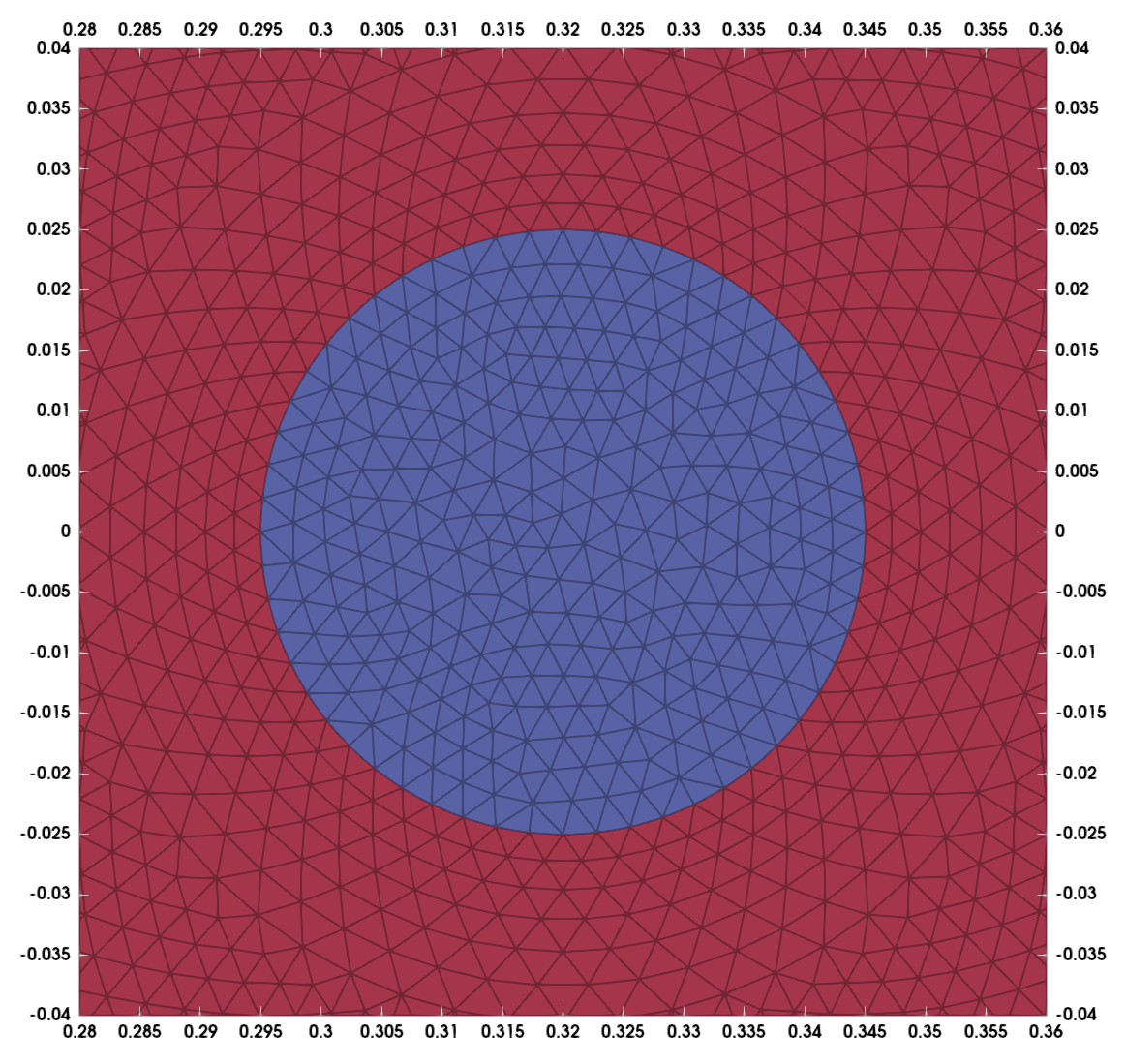}}
    \caption{Example \ref{sec:sb}: Geometry setup and zoomed-in meshes around the bubble.
  }
    \label{fig:sb0}
\end{figure}

This test case corresponds to the interaction of shock wave with a cylindrical Helium bubble surrounded by air at rest \cite{Quirk96}.
We use the same steup as in \cite[Section 8.4]{Galera10}.
The initial domain is a rectangular box $\Omega_0 = (0, L)\times (-H/2, H/2) = (0, 0.650)\times (-0.089, 0.089)$, which includes a  circular bubble with center $(0.320, 0)$ and radius $r_b=0.025$. Initial data are shown in Figure~\ref{fig:sb0}(a). Wall boundary conditions at each boundary is prescribed except at the right boundary, where we impose a piston-like boundary condition with inward velocity $(-124.824, 0)$.
The left going shock wave hits the bubble at time 
$t_i = 668.153\times 10^{-6}$. The final time of simulation is $t_{f}=t_i+674\times 10^{-6}=1342.153\times 10^{-6}$, which corresponds to the time where experimental shadow-graph extracted from \cite{Haas87} is displayed in \cite{Quirk96}.

\begin{figure}[ht!]
    \centering
    \includegraphics[width = .85\textwidth]{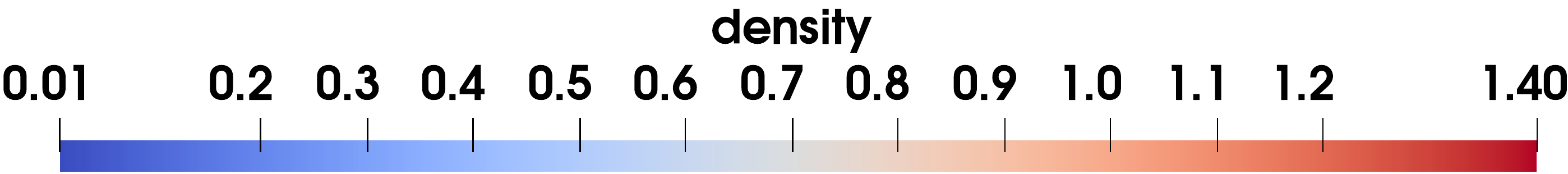}
\subfigure[$k$=2, coarse mesh]{\includegraphics[width = .19\textwidth]{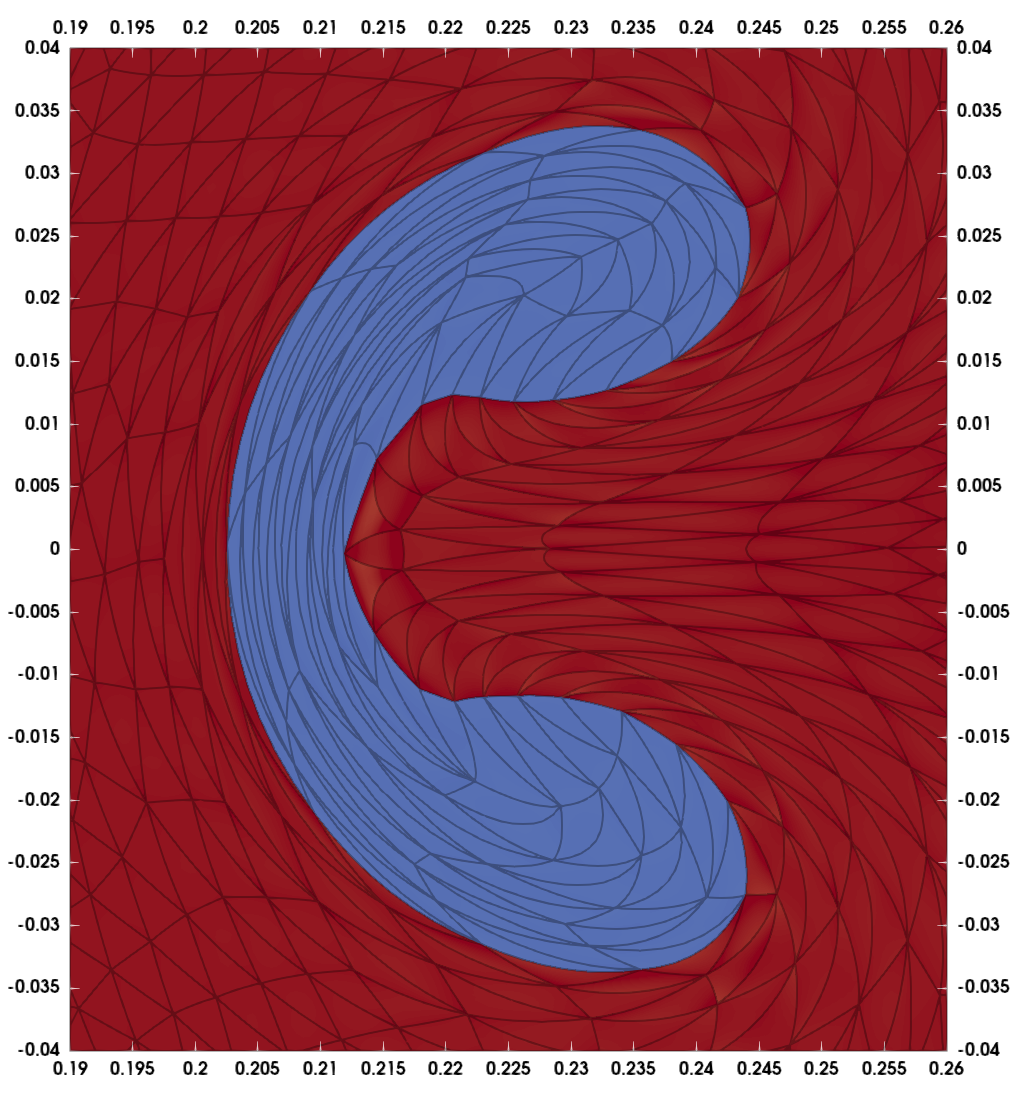}}
\subfigure[$k$=2, fine mesh]{\includegraphics[width = .19\textwidth]{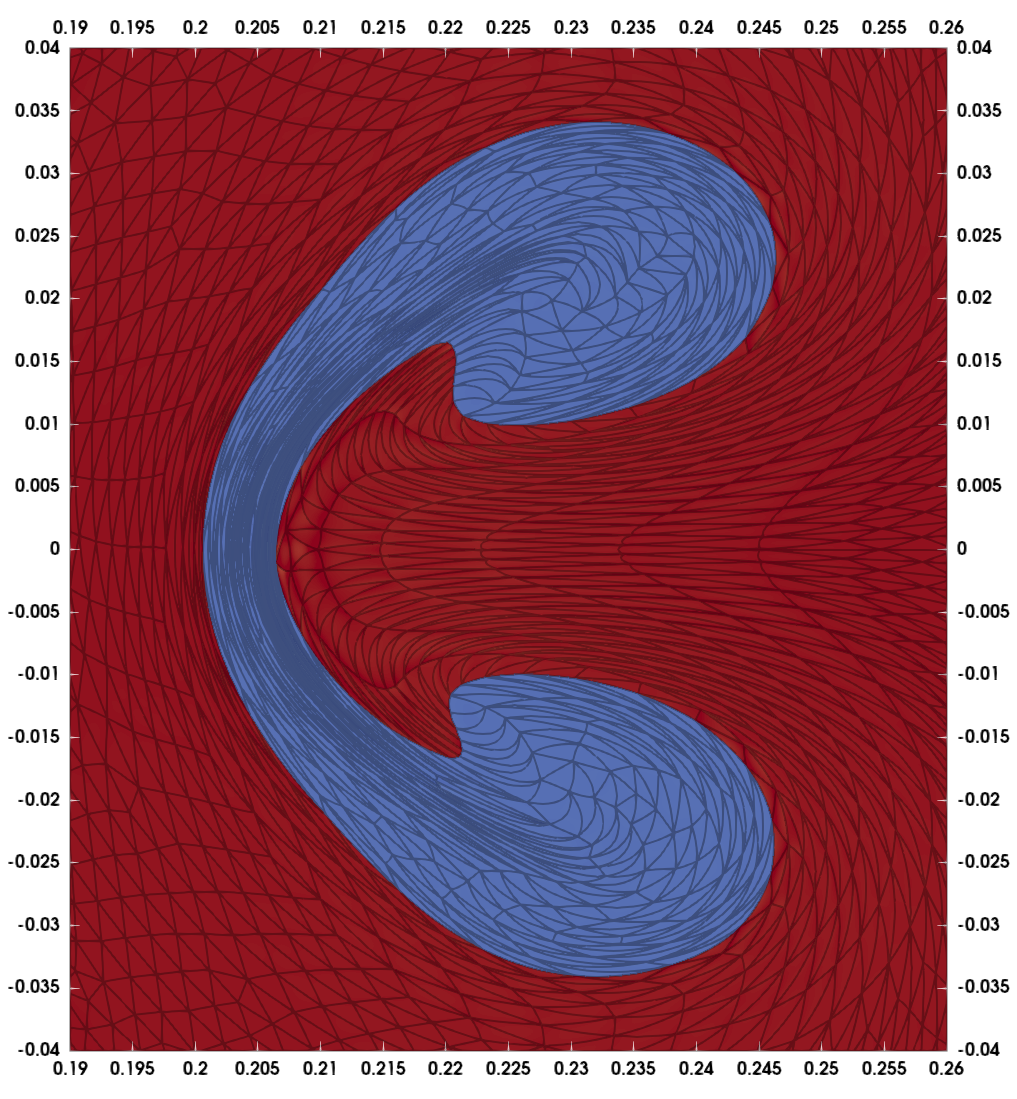}}
\subfigure[$k$=4, coarse mesh]{\includegraphics[width = .19\textwidth]{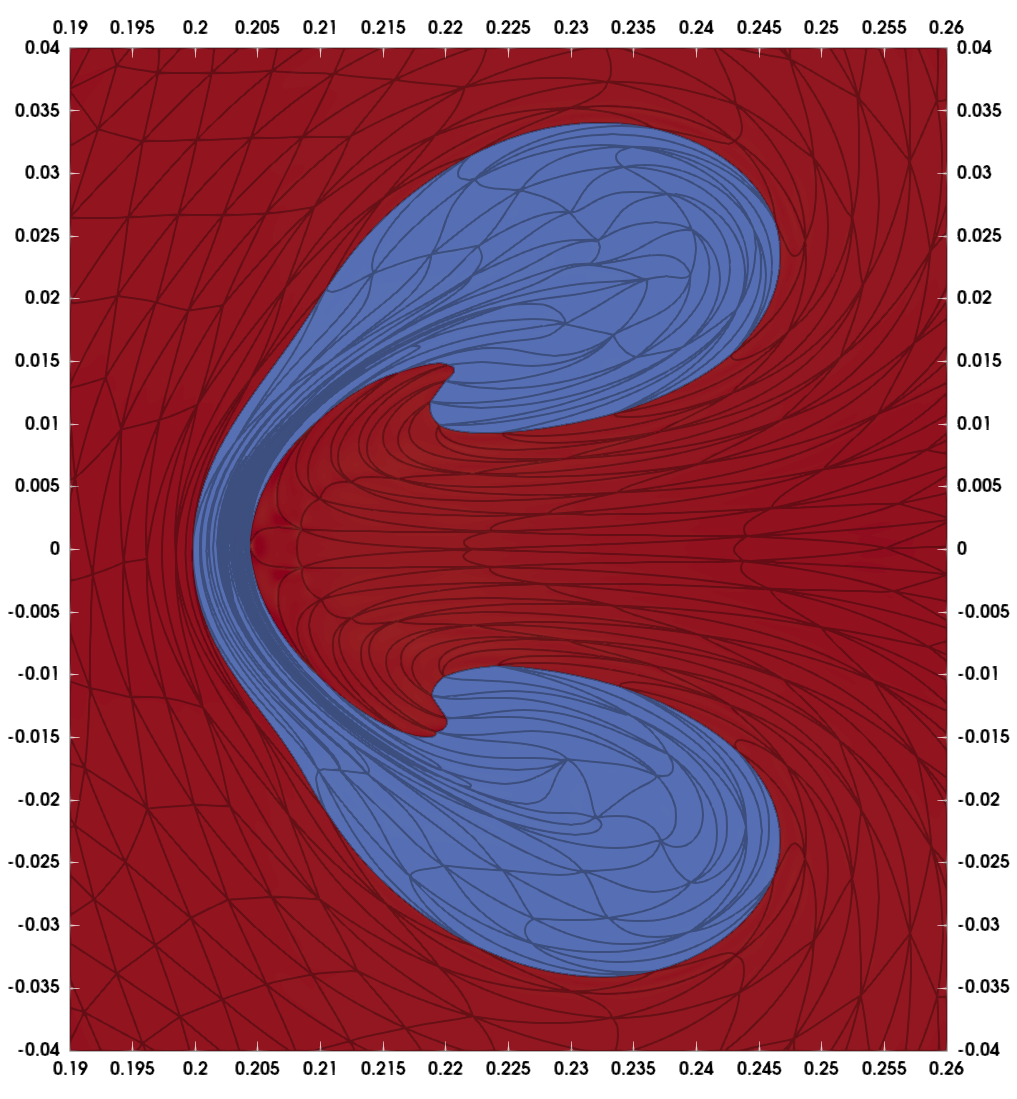}}
\subfigure[$k$=4, fine mesh]{\includegraphics[width = .19\textwidth]{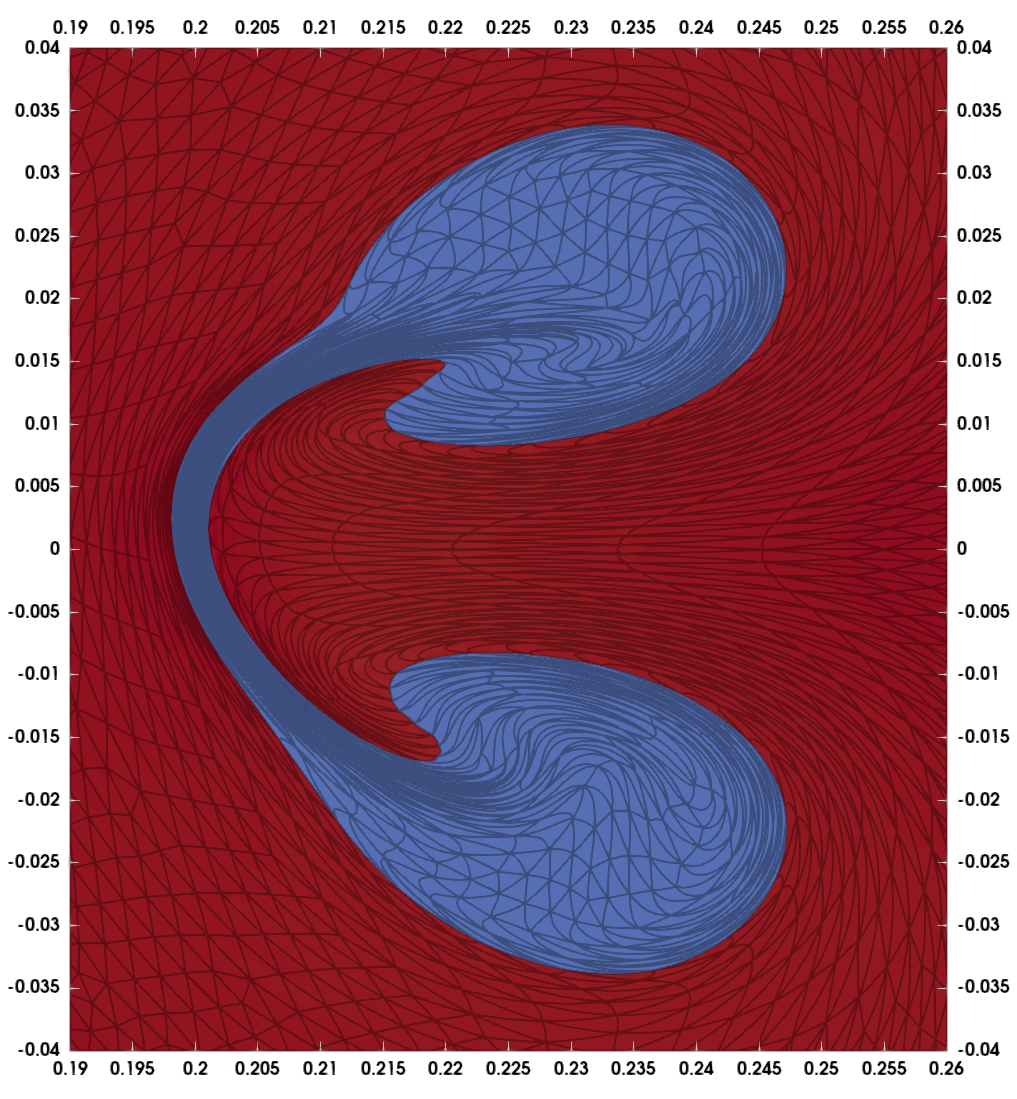}}
\subfigure[Data from \cite{Haas87}]{\includegraphics[width = .19\textwidth]{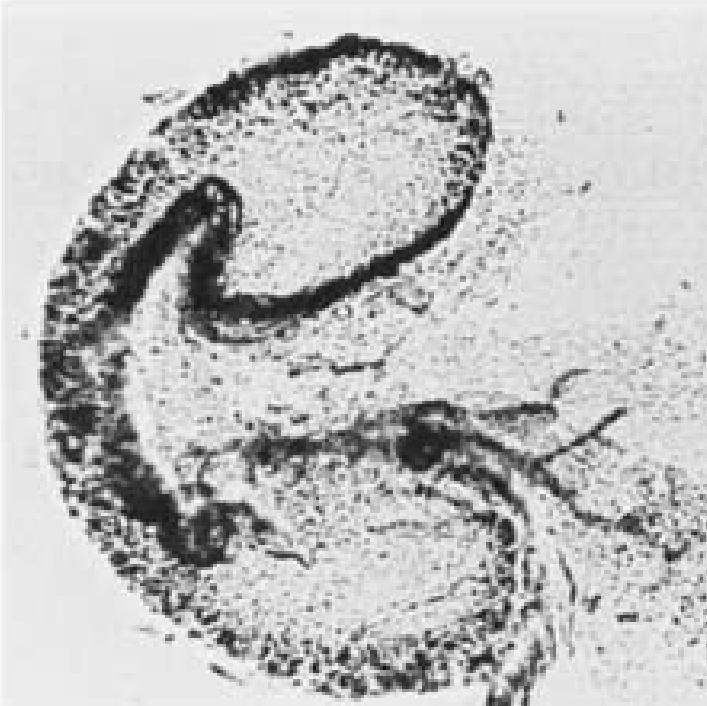}}
    \caption{Example \ref{sec:sb}: Zoom on the deformed bubble at final time $t$= 1342.153e-6.
   Here (e) is the Schlieren image from experimental data \cite{Haas87}.
  }
    \label{fig:sb}
\end{figure}

We use two unstructured triangular meshes that is fitted to the bubble boundary for this problem. The coarse mesh has 8324 triangular cells whose mesh size is $h=H/32$, while the fine mesh 
has 33526 cells whose mesh size is $h = H/64$; see Figure \ref{fig:sb0}(b)-(c) for the zoom-in view of the two meshes around the bubble.
The BDF2 time stepping is used in combination with polynomial degree $k=2$ and $k=4$ on these two meshes. The zoom-in view around the deformed bubble at final time $t_f$ are shown in Figure~\ref{fig:sb}. We observe that the location and shape of the deformed bubble is similar for each simulation, with a better resolution being obtained on a finer mesh with a higher order polynomial degree. These shapes are also qualitatively similar to the experimental Schlieren image 
in Figure~\ref{fig:sb}(e) obtained from  \cite{Haas87}.

\begin{figure}[ht!]
    \centering
    \includegraphics[width = .85\textwidth]{figs/denRangeX.pdf}
\subfigure[$k=4$, coarse mesh, $t$=800e-6]{\includegraphics[width = .45\textwidth]{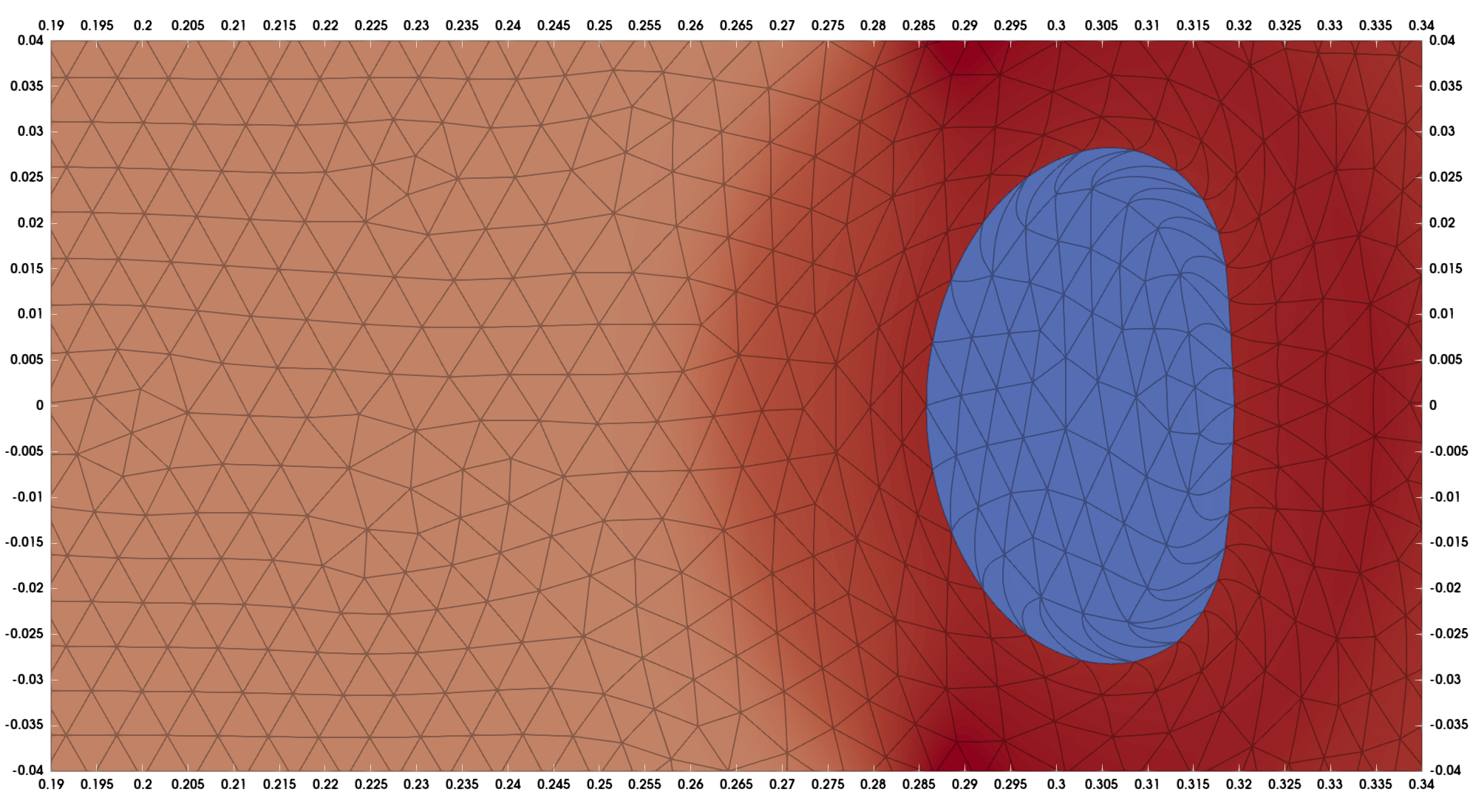}}
\subfigure[$k=4$, fine mesh, $t$=800e-6]{\includegraphics[width = .45\textwidth]{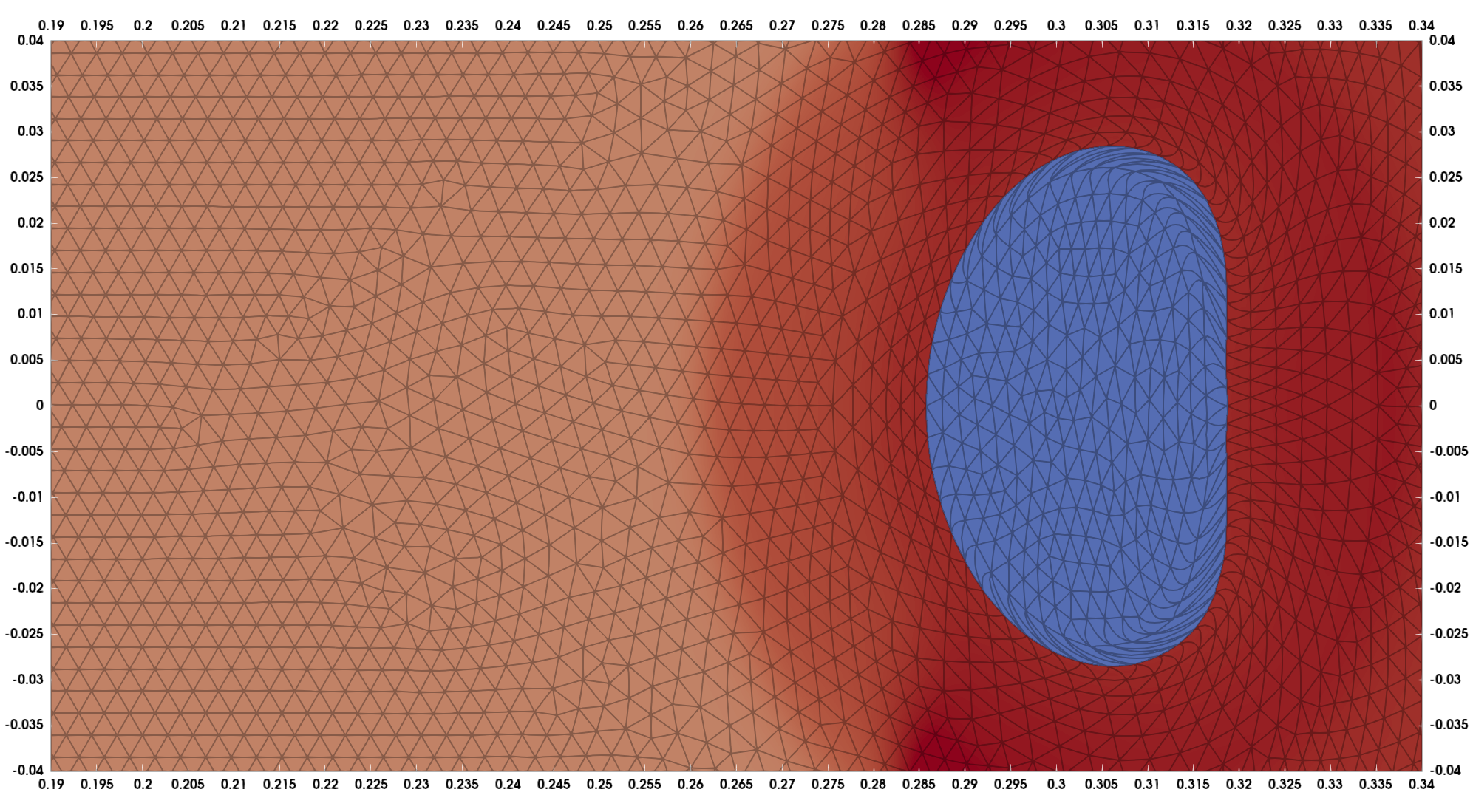}}
\subfigure[$k=4$, coarse mesh, $t$=1100e-6]{\includegraphics[width = .45\textwidth]{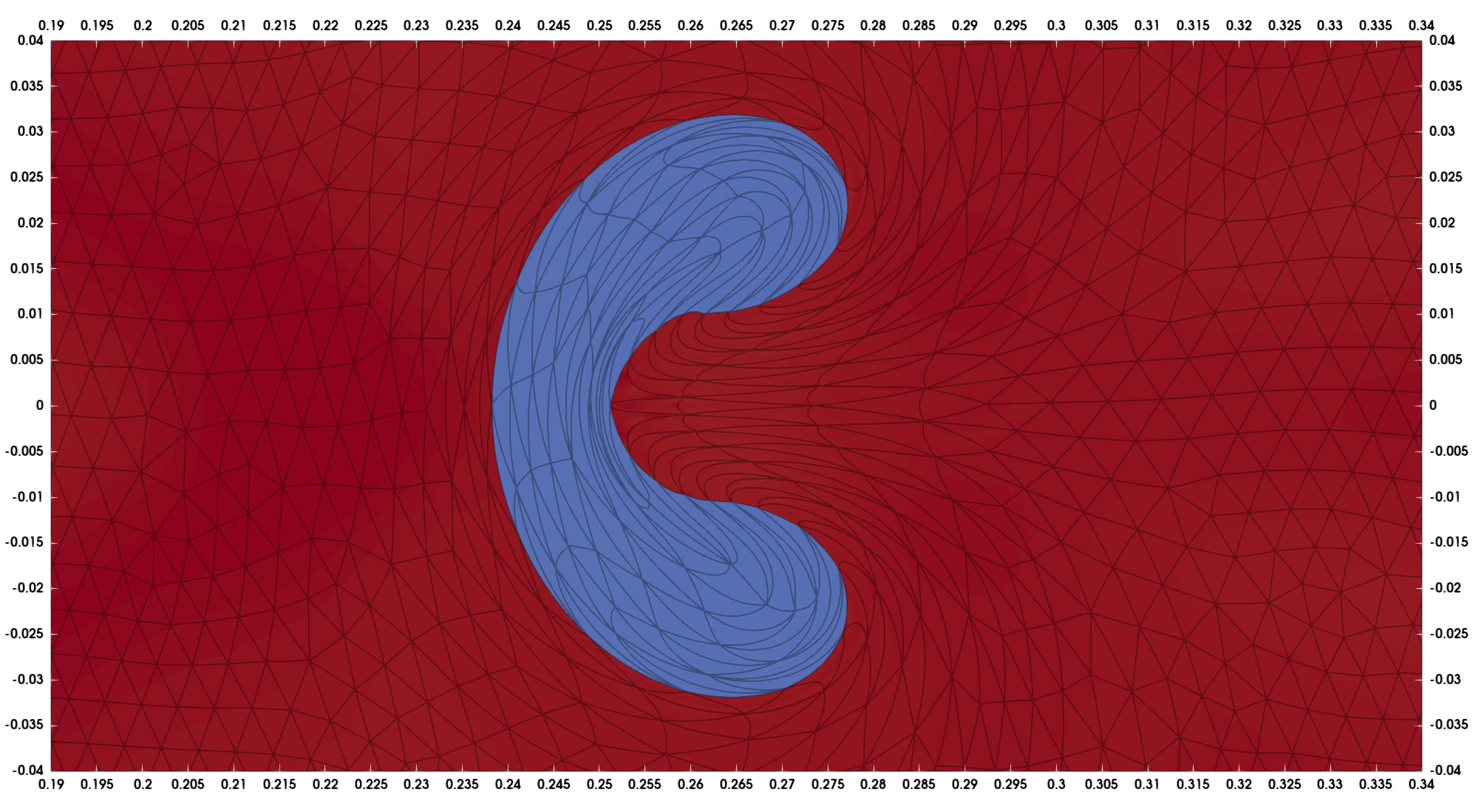}}
\subfigure[$k=4$, fine mesh, $t$=1100e-6]{\includegraphics[width = .45\textwidth]{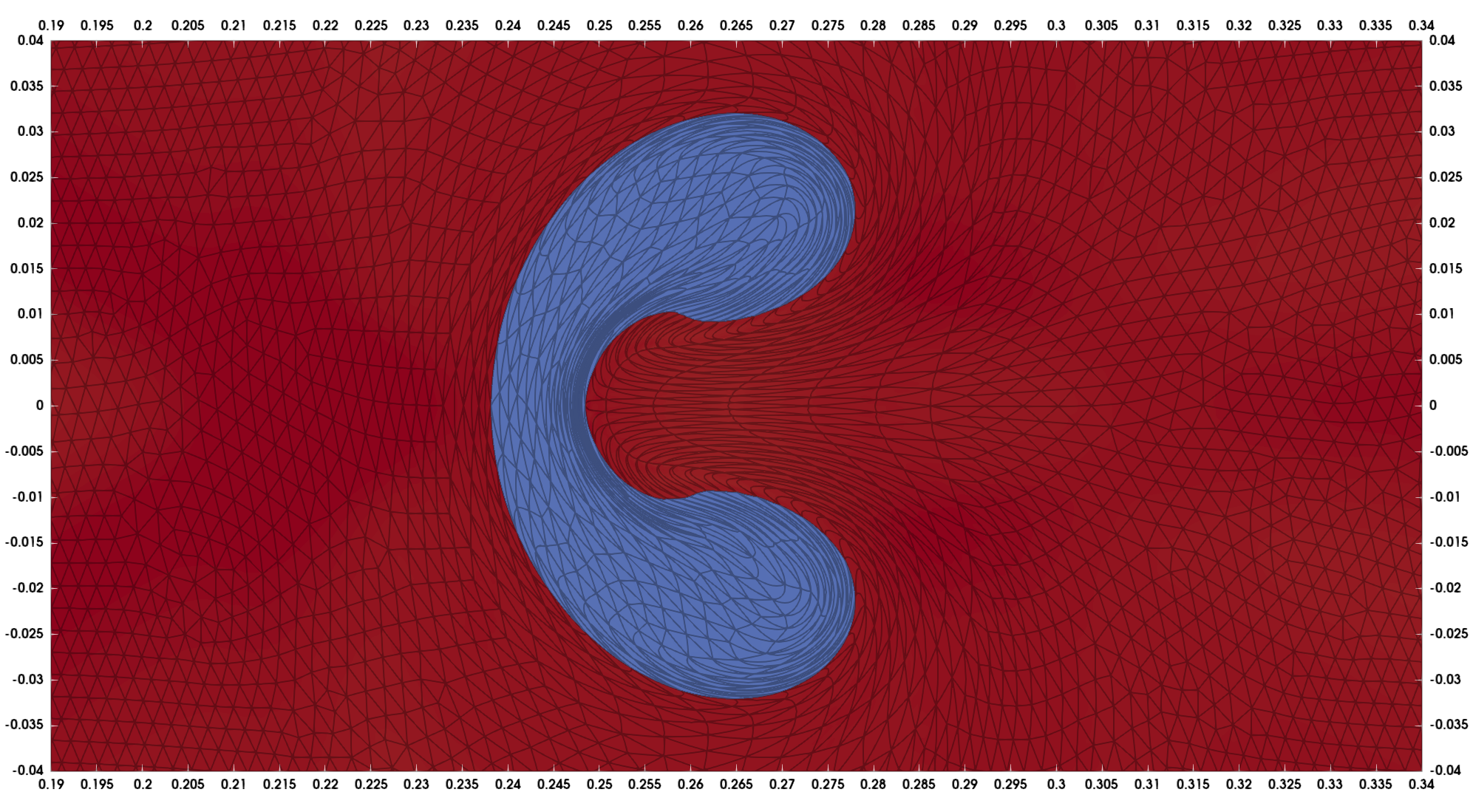}}
\subfigure[$k=4$, coarse mesh, $t$=1342.153e-6]{\includegraphics[width = .45\textwidth]{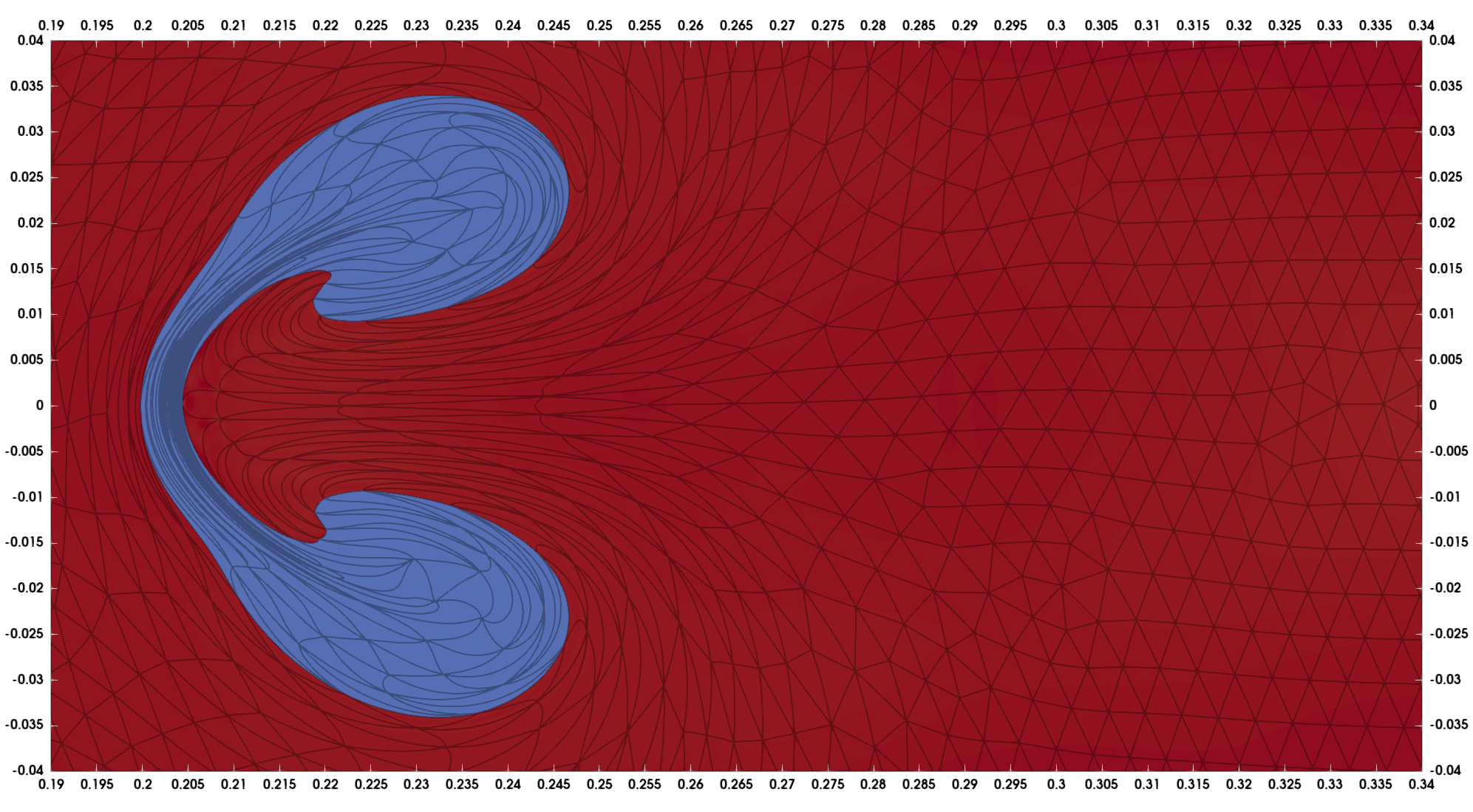}}
\subfigure[$k=4$, fine mesh, $t$=1342.153e-6]{\includegraphics[width = .45\textwidth]{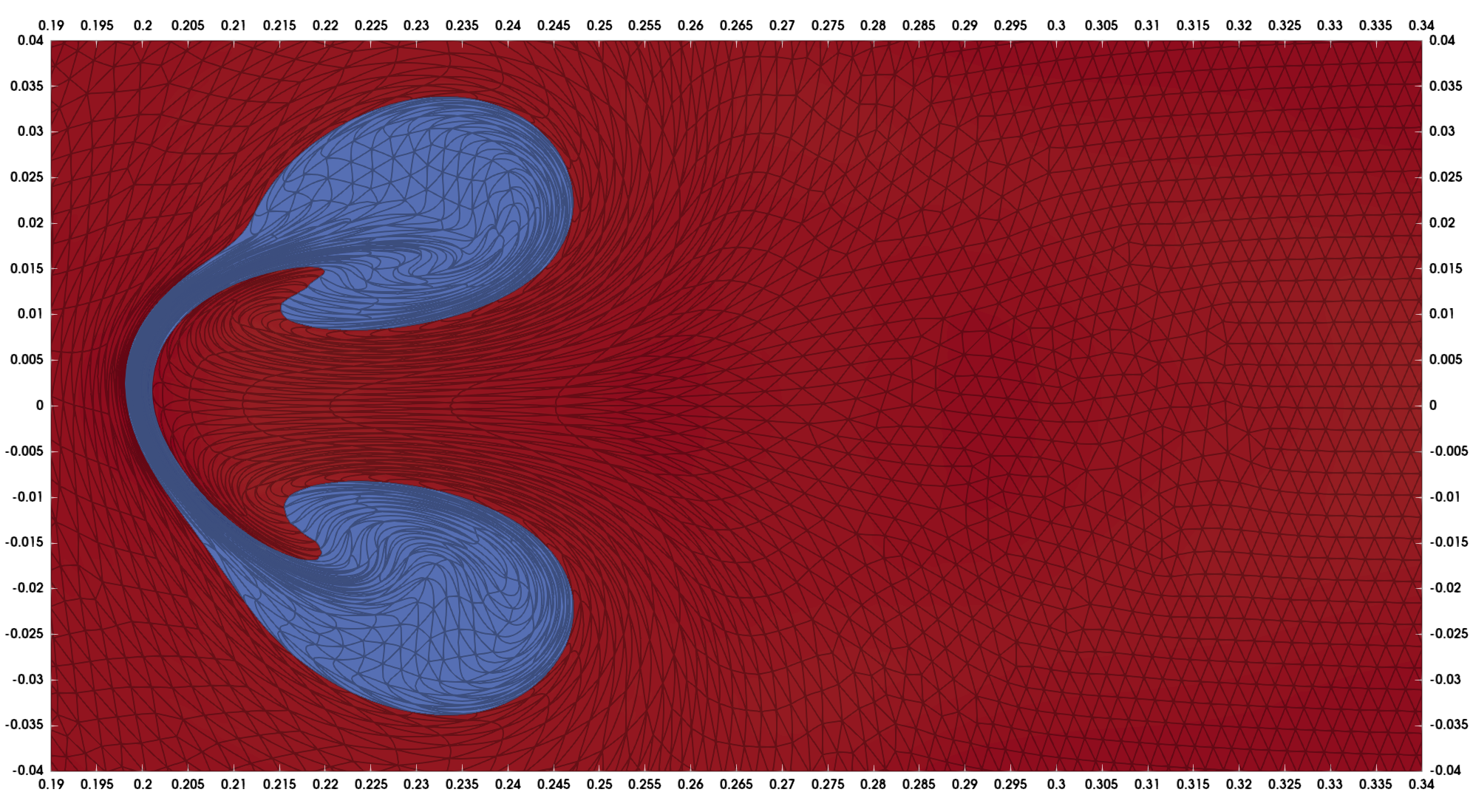}}
    \caption{Example \ref{sec:sb}: Density contour in the vicinity of the bubble at various times.
  }
    \label{fig:sb2}
\end{figure}

We display in Figure \ref{fig:sb2} the time evolution of the bubble at times $t=800\times 10^{-6}, 1100\times 10^{-6}, 1342.153\times 10^{-6}$ for $k=4$ on the two meshes. We note that the results obtained with both meshes are quite similar.

\subsection{Multimaterial implosion in cylindrical geometry}
\label{sec:imp}
\begin{figure}[ht!]
    \centering
{\includegraphics[width = .45\textwidth]{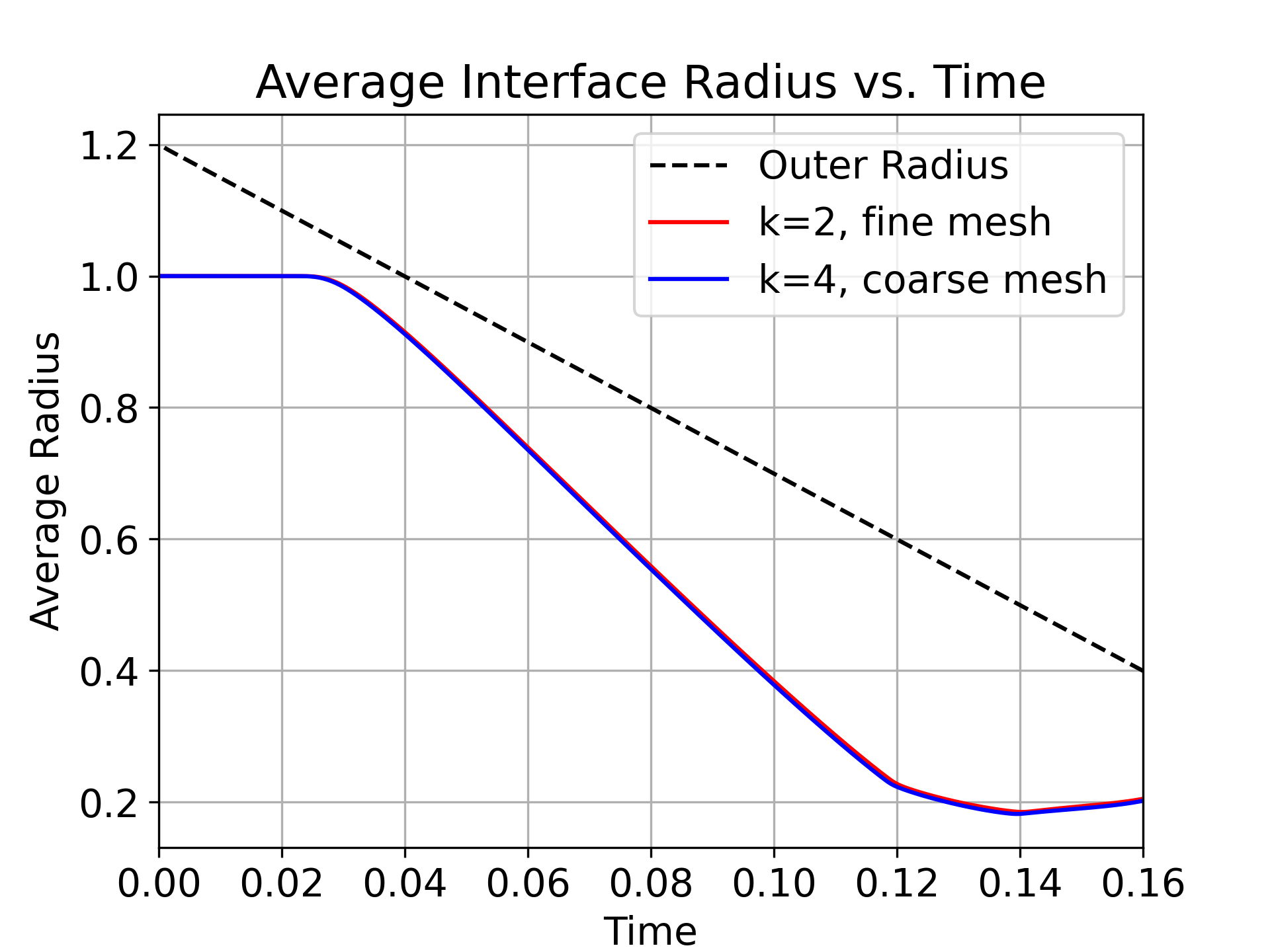}}
{\includegraphics[width = .45\textwidth]{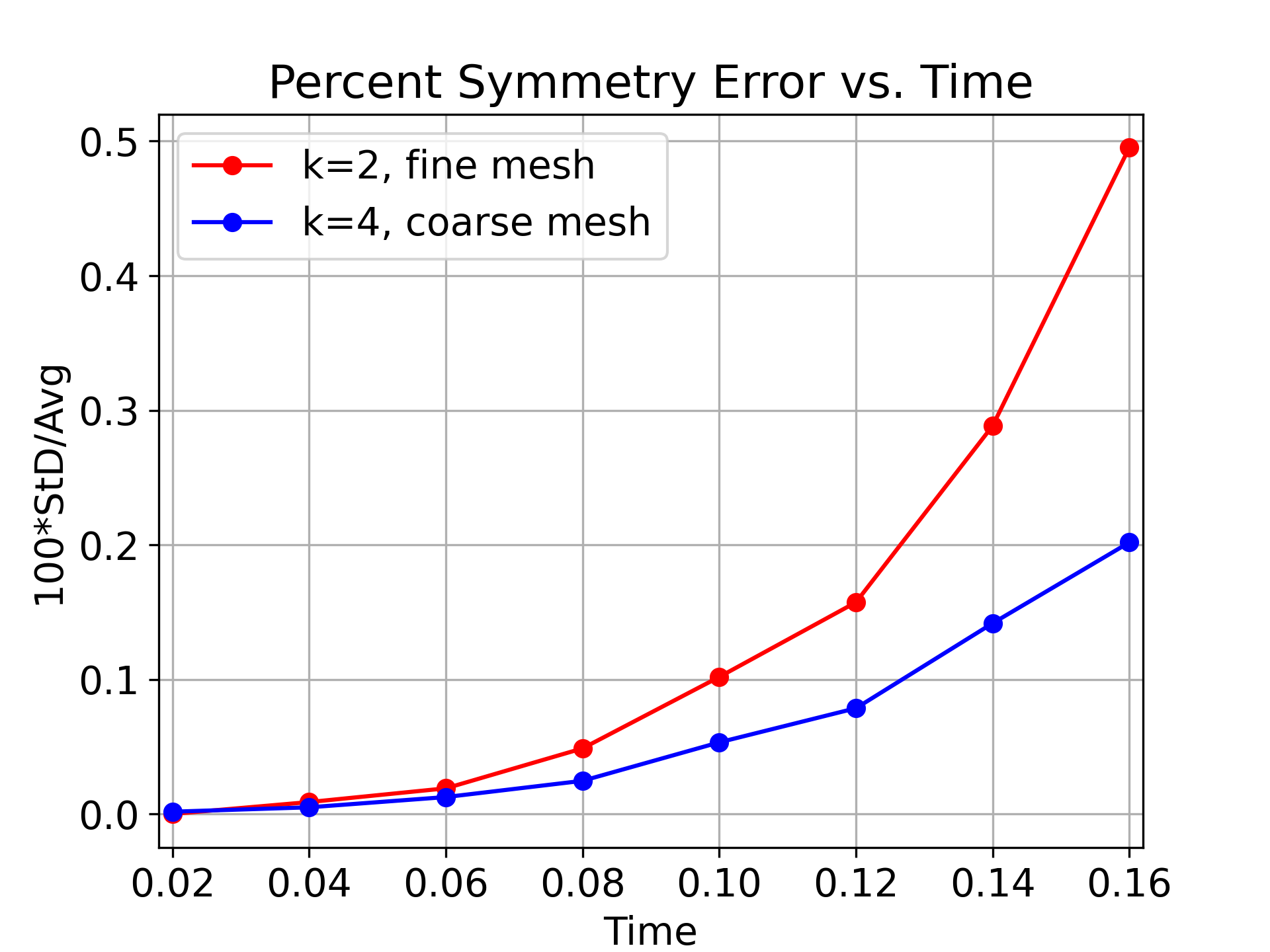}}
    \caption{Example \ref{sec:sb}: Density contour in the vicinity of the bubble at various times.
  }
    \label{fig:imp2}
\end{figure}

The aim of this example is to access the capability of the implicit high-order Lagrangian scheme to handle a multi-mode implosion in cylindrical geometry.
Here we consider a simple 1D multimaterial implosion problem on unstructured 2D triangular meshes. 
The problem consists of a low-density material with $\rho_1=0.05$ in the radial range $r\in [0,1]$
surrounded by a shell of high density material 
with $\rho_2=1.0$ in the radial range $r\in (1.0, 1.2]$. Both materials are initially at rest with pressure $p=0.1$ and adiabatic index $\gamma = 5/3$.
This problem was originally proposed in \cite{Galera10} with a time dependent pressure source on the outer radial surface $r=1.2$.
Here we use the modification in \cite{Dobrev12} that applies a constant radial velocity source of $\bld u= -5(x, y)/1.2$ on the outer boundary, which drives a cylindrical shock wave inwards.
Due to the 1D setup (symmetry) , the material interface should be a function of radius only for all time. The movement of the interface radius 
against time is shown on the left panel of Figure~\ref{fig:imp2}. 
We clearly observe the deceleration of the interface starting around $t=0.12$, and the so-called stagnation phase is reached around $t=0.14$ where the radius obtains its minimum value.
The flow becomes Rayleigh-Taylor unstable 
after this time as a small perturbation of the interface will grow exponentially as function of time due to the fact that the light fluid inside is pushing the heavy fluid outside after the stagnation phase. 
It is very challenging to preserve the interface symmetry for a Lagrangian scheme on general unstructured meshes, especially when time past the stagnation phase. 

\begin{figure}[ht!]
    \centering
\subfigure[$k=2$, fine mesh, $t$=0]{\includegraphics[width = .3\textwidth]{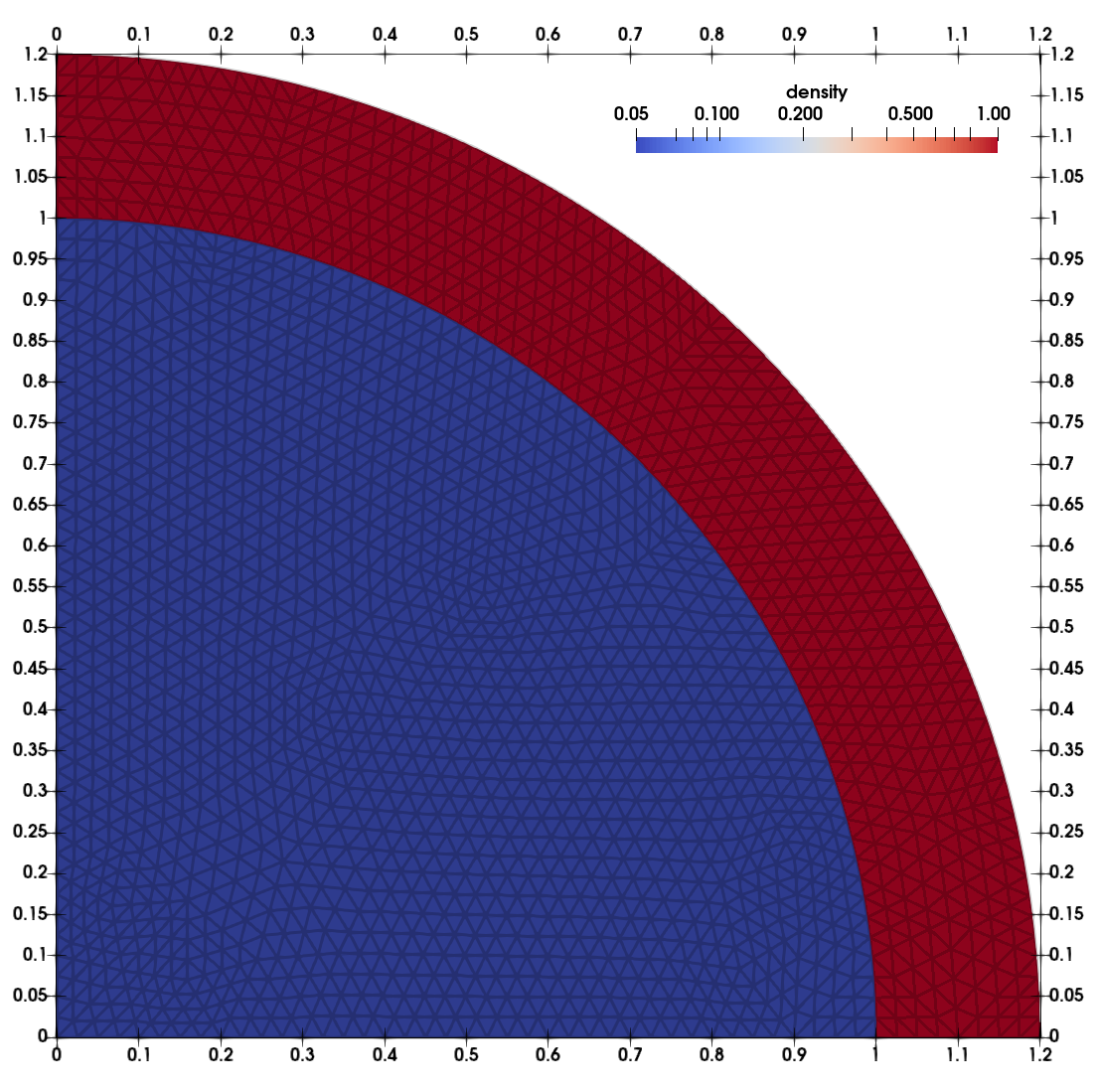}}
\subfigure[$k=2$, fine mesh, $t$=0.08]{\includegraphics[width = .3\textwidth]{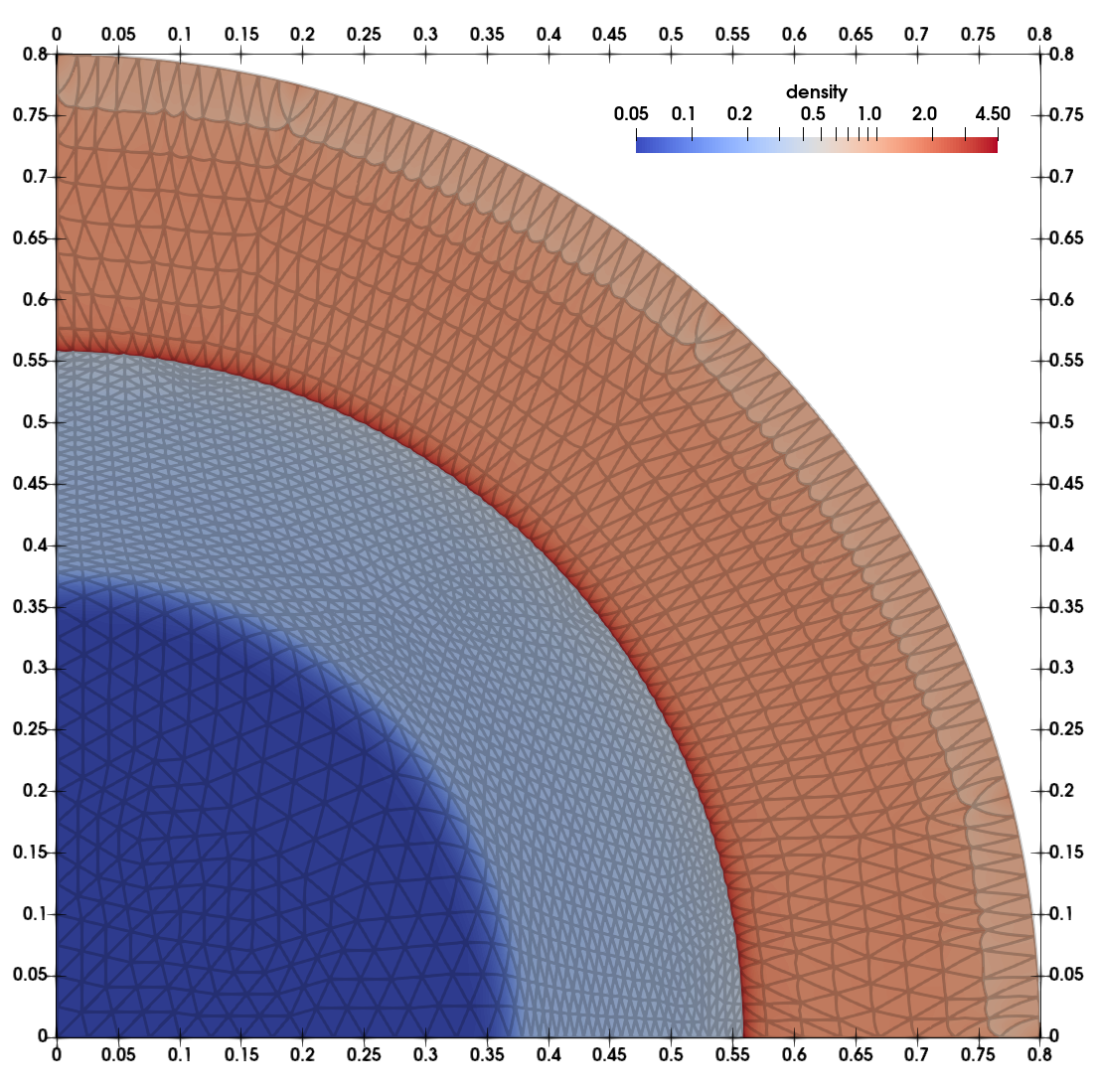}}
\subfigure[$k=2$, fine mesh, $t$=0.16]{\includegraphics[width = .3\textwidth]{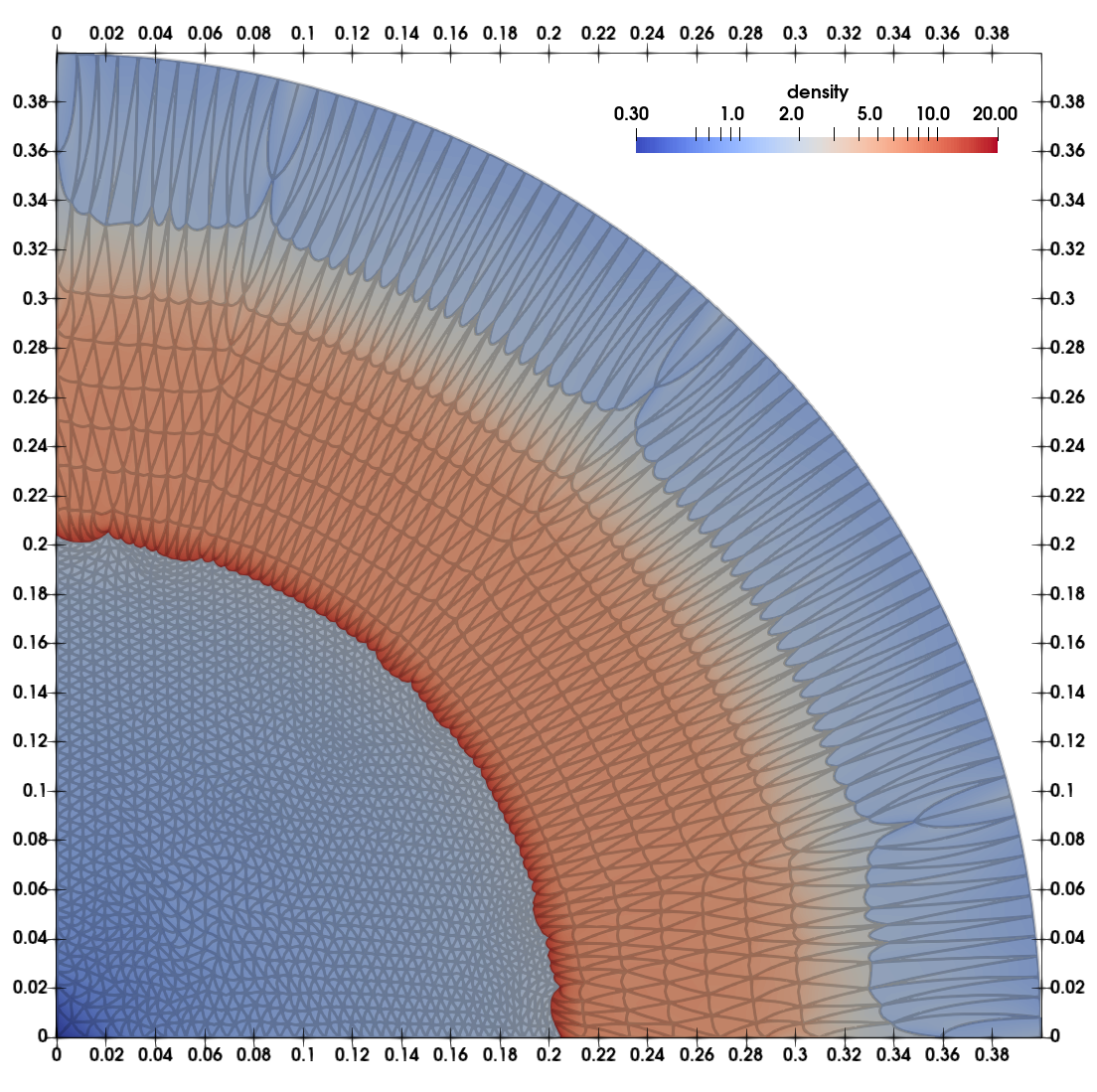}}
\subfigure[$k=4$, coarse mesh, $t$=0]{\includegraphics[width = .3\textwidth]{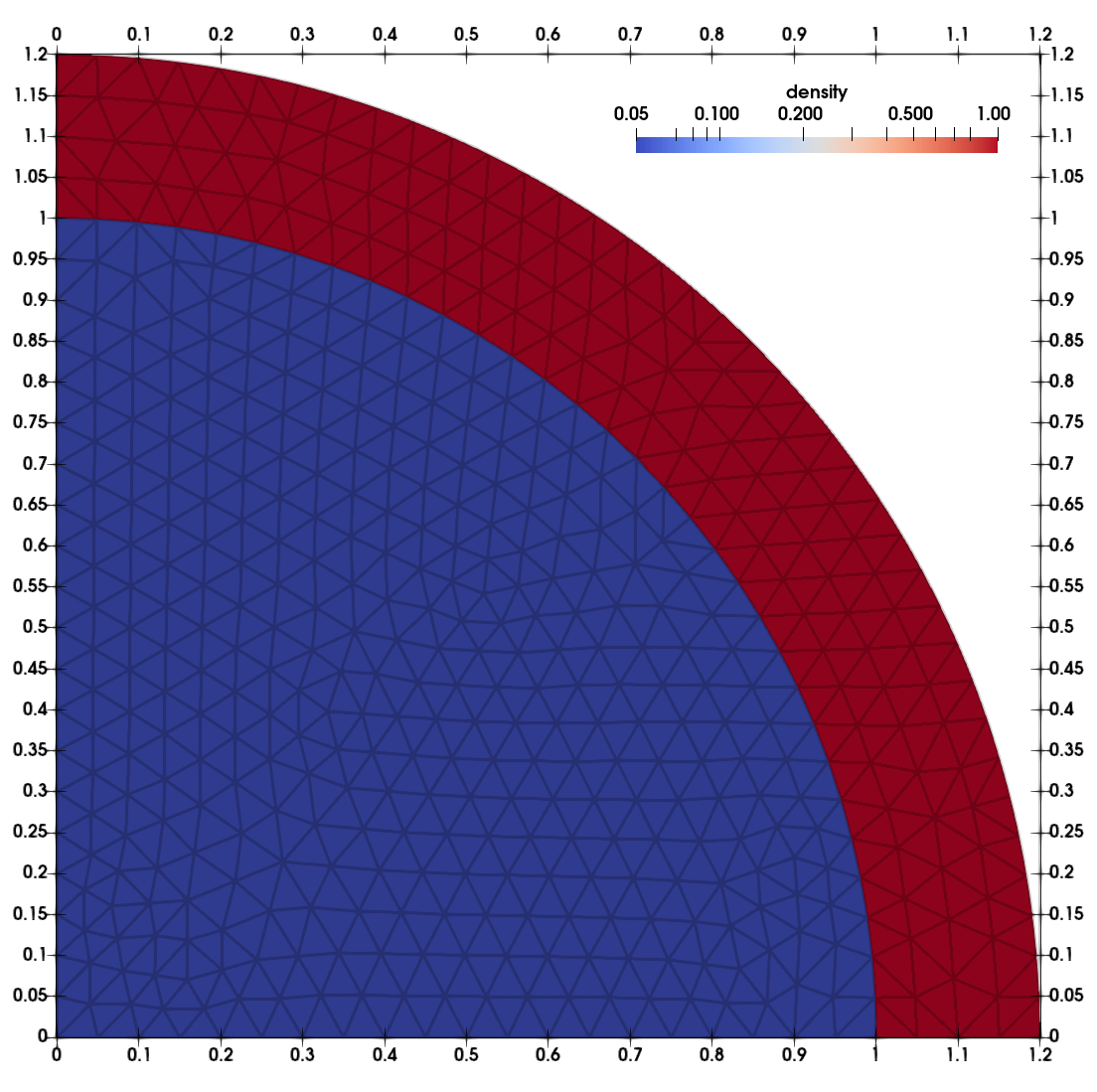}}
\subfigure[$k=4$, coarse mesh, $t$=0.08]{\includegraphics[width = .3\textwidth]{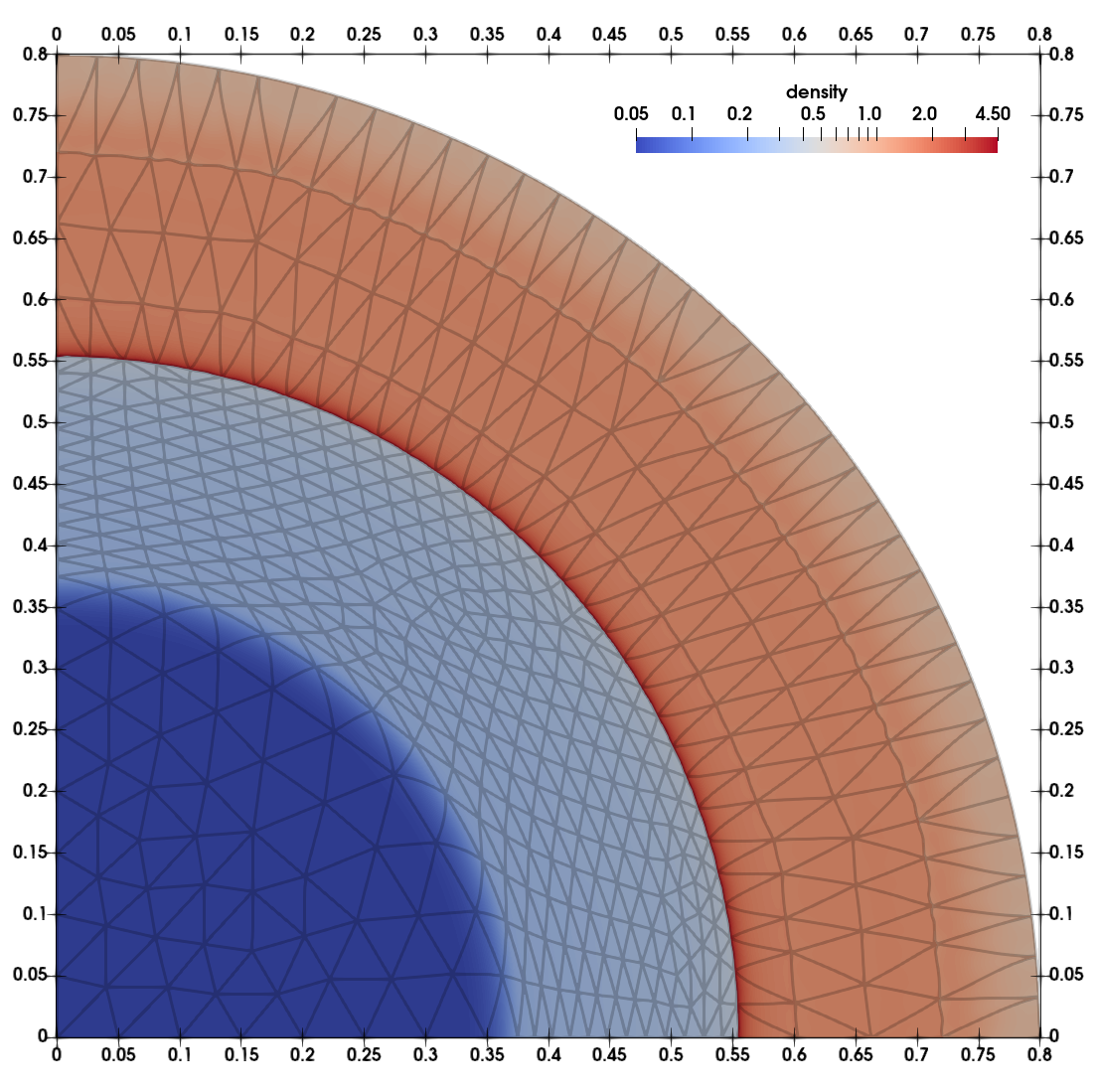}}
\subfigure[$k=4$, coarse mesh, $t$=0.16]{\includegraphics[width = .3\textwidth]{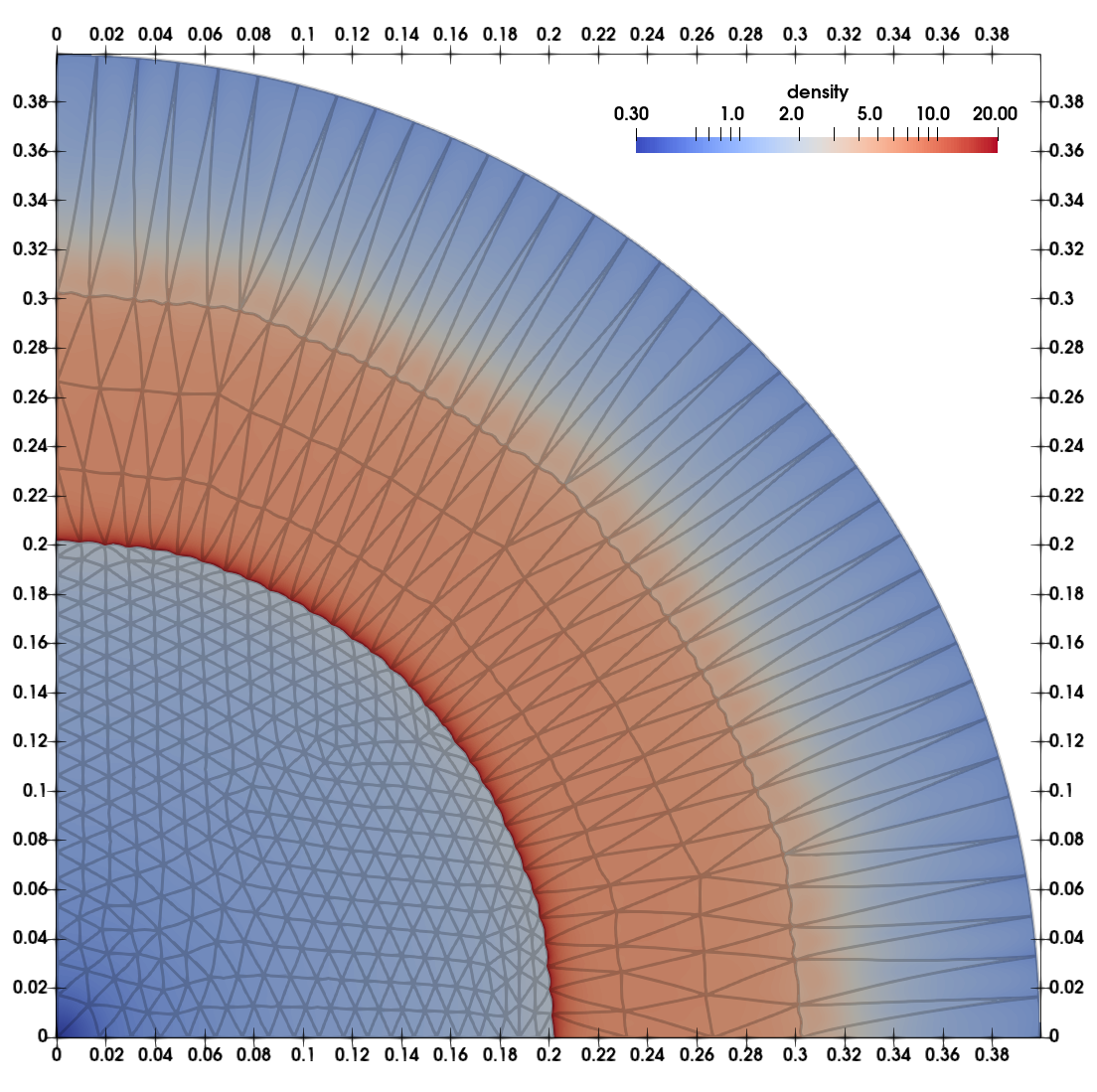}}
    \caption{Example \ref{sec:sb}: Density contour in the vicinity of the bubble at various times.
  }
    \label{fig:imp}
\end{figure}

Due to symmetry of the problem, we take the computational domain to be a quarter circle 
$\Omega=\{(x, y): x\ge 0 , y\ge 0, x^2+y^2 \le  1.2^2\}$
with wall boundary conditions on the left and bottom boundaries.
We apply the BDF2 scheme with polynomial degree $k=2, 4$ on two set of unstructured triangular meshes.
The coarse mesh has mesh size $h=0.05$ and is used for $k=4$, while the fine mesh is a uniform refinement of the coarse mesh and is used for $k=2$; 
see Figure~\ref{fig:imp}(d) for the coarse mesh and 
Figure~\ref{fig:imp}(a) for the fine mesh.
We take $q_1=1$ and $q_2=4$ as the artificial viscosity coefficients and set CFL number equals 0.25.
The two simulations have the same number of velocity DOFs, and requires a similar amount of total time steps to reach the final time $t=0.16$. 
Density plots in log-scale on the deformed meshes are shown in Figure~\ref{fig:imp} at times 
$t=0.08$ and $t=0.16$, along with the initial density at $t=0$. We observe while the results for $t=0.08$ are similar for both cases, where radial symmetry 
of the material interface is preserved quite well. The results for $k=4$ at $t=0.16$ (which past the stagnation phase) is better than that for $k=2$ in terms of the interface symmetry. 
In Figure~\ref{fig:imp2} we plot the time evolution of the average radius of the material interface, and the normalized standard deviation
of this radius
 at different times as an indication of the symmetry error over time. 
 We observe similar results for the average radius 
 for both simulations, and a smaller symmetry error for the higher order case.

\section{Conclusion}
\label{sec:conclude}
We presented a class of high-order variational Lagrangian schemes for compressible flow.
The discrete EnVarA approach is used to derive the high-order spatial finite element discretization.
Features of our spatial discretization
include mass/momentum/energy conservation and entropy stability. 
Fully implicit time stepping is then applied to the resulting ODE system. Each time step requires a nonlinear system solve for the velocity DOFs only.
Ample numerical results are shown to support the good performance of the proposed scheme.
We plan to extend this pure Lagrangian scheme to the arbitrary Eulerian Lagrangian framework, which has the potential to address the mesh distortion issue.

\section*{Appendix}
\label{sec:ap}
Here we briefly discuss the isothermal case \cite{Giga2017} where the temperature does not change over time.
In this case, the free energy \eqref{ener-den} is a function of density only: $\psi = \psi(\rho)$. 
Typical choices include $\psi(\rho)=\alpha \rho^\gamma$ or 
$\psi(\rho)=\alpha \rho\log(\rho)$.

This model only has two thermodynamic variables: density $\rho$ and pressure 
\[
p = \psi_{\rho}\rho-\psi.
\]
The EnVarA derivation leads to the 
following model equations; see also \cite{Giga2017}:
\begin{subequations}
\label{model2}
\begin{align}
\dot{\bld x} =& \bld u, \quad \quad \dot{(\rho J)} = \;0,\\
\rho\dot{\bld u} = &\;-\nabla p +\nabla\cdot\left(\eta\nabla_s\bld u+(\xi-\frac23\eta)(\nabla\cdot\bld u)\bld I\right),
\end{align}
\end{subequations}
where $p = \psi_{\rho}\rho-\psi$.

The spatial discretization for the model \eqref{model2}
is summarized below:
\begin{algorithm}[H]
\caption{Variational Lagrangian scheme for the model \eqref{model2}.}
\label{alg:X}
\begin{algorithmic}
\State 
Find $\bld x_h, \bld u_h\in \bld V_h^k$ such that 
\begin{subequations}
\label{BEX}
\begin{align}
\dot{\bld x}_h = &\;{\bld  u}_h, \\
\label{u-ThX}
(\rho_0\dot{\bld u}_h, \bld \varphi_h)_h
-\left(p_h, \nabla\cdot\bld \varphi_h\right)_h
+(\sigma_h, \nabla\bld \varphi_h)_h
= &\;0 , \quad \forall \varphi_h\in\bld V_h^k,
\end{align}
\end{subequations}
where the pressure $p_h\in W_h^k$
is a function of the density approximation with
\[
p_i^{\ell} = \psi_\rho(\rho_i^{\ell})\rho_i^{\ell}-\psi(\rho_i^{\ell}),
\]
in which
$\rho_h\in W_h^k$ satisfies mass conservation \eqref{mass-c}, 
and the 
stress
\[
\sigma_h=(\eta+\mu_{av}) J_h\nabla_s{\bld u}_h+ (\xi-\frac23\eta)J_h\nabla\cdot{\bld u}_h\bld I,
\]
in which the artificial viscosity coefficient $\mu_{av}$ is given in \eqref{av}.
\end{algorithmic}
\end{algorithm}
It is clear that the spatial discretization in Algorithm \ref{alg:X} is mass and momentum conservative. 
Next we prove Algorithm \ref{alg:X} is also entropy stable.
Using a similar derivation as in 
\eqref{potential-V} and
\eqref{var-A} and the definition of density and pressure, 
we obtain 
\[
\frac{d}{dt}(\psi(\rho_h)J_h, 1)_h
= -(p_hJ_h, \nabla\cdot \bld u_h)_h.
\]
Taking test function $\bld \varphi_h=\bld u_h$
in \eqref{u-ThX} and use the above relation, we get
\[
\frac{d}{dt}\underbrace{\left(\frac12\rho_0|\bld u_h|^2
+\psi(\rho_h)J_h, 1\right)_h}_{\text{total entropy}}
= -(\sigma_h, \nabla\bld u_h)\le 0.
\]

The ODE system in Algorithm \ref{alg:X} can be discretized using the implicit schemes discussed in 
Section \ref{sec:time}.
Here we present a variational implicit time discretization 
similar to the backward Euler scheme in Section \ref{sec:BE}.
Given data $\bld x^{n-1}, \bld u_h^{n-1}\in V_h^k$ at time $t^{n-1}$ and time step size $\delta t$, 
denote the following discrete energy functional for the flow map $\bld x_h$:
\begin{align}
E_h(\bld x_h):=
\left(\frac{\rho_0|\bld x_h-\bld x_h^{n-1}-\delta t \bld u_h^{n-1}|^2}{2\delta t^2}
+\psi\left(\rho_0/J_h\right)J_h, 1\right)_h
+\frac12(\Delta_h, 1)_h,
\end{align}
where the dissipation
\[
\Delta_h:=
(\eta+\mu_{av}^{n-1}) J_h|\nabla_s\frac{{\bld x}_h-{\bld x}_h^{n-1}}{\delta t}|^2+ (\xi-\frac23\eta)J_h|\nabla\cdot\frac{{\bld x}_h-{\bld x}_h^{n-1}}{\delta t}|^2, 
\]
with artificial viscosity $\mu_{av}^{n-1}$ explicitly evaluated at time level $t^{n-1}$, and $J_h=|\nabla_X\bld x_h|$.
The flow map at next time level is obtained by solving the following minimization problem:
\begin{align}
\label{mini}
\bld x_h^{n}:=\mathrm{argmin}_{\bld x_h\in V_h^k, J_h>0}
E_h(\bld x_h).
\end{align}
Assuming piecewise constant viscosity coefficients 
$\eta, \xi$, and $\mu_{av}^{n-1}$, the Euler-Lagrangian equation for this minimization problem is 
simply the Backward Euler scheme in \eqref{BE} applied to the ODE system \eqref{BEX}.
We note that such energy minimization interpretation is not available for the non-isothermal case discussed in Section \ref{sec:time} due to the temperature equation \eqref{theta-Th1}.

\bibliographystyle{siam}

\end{document}